\documentclass[11pt,notitlepage,twoside,a4paper]{amsart}
 \usepackage{amsfonts}

\usepackage{amsmath,amssymb,enumerate}

\usepackage{epsfig,fancyhdr,color}

\usepackage{amssymb}
\usepackage{amsmath,amsthm}
\usepackage{latexsym}
\usepackage{amscd}
\usepackage{psfrag}
\usepackage{graphicx}
\usepackage[latin1]{inputenc}
\usepackage[all]{xy}
\usepackage[mathcal]{eucal}

\definecolor{NoteColor}{rgb}{1,0,0}

\renewcommand{\textsc}{\textcolor{red}}

\newtheorem{theorem}{\rm\bf Theorem}[section]
\newtheorem{proposition}[theorem]{\rm\bf Proposition}

\newtheorem*{theorem 1}{\rm\bf Proposition 1}
\newtheorem*{theorem 2}{\rm\bf Proposition 2}

\theoremstyle{definition}

\theoremstyle{remark}

\newtheorem{problem}[theorem]{\rm\bf Problem}

\def\interieur#1{\mathord{\mathop{\kern 0pt #1}\limits^\circ}}

\title[Early history]{On the early history of moduli and Teichm\"uller spaces}

\author{Norbert A'Campo}
\address{N. A'Campo, Universit\"at Basel,  Mathematisches Institut, 
Rheinsprung 21, 4051 Basel, Switzerland.}
\email{norbert.acampo@unibas.ch}

\author{Lizhen Ji}
\address{L. Ji, Department of Mathematics, University of Michigan\\ Ann Arbor, MI 48109, USA.}
\email{lji@umich.edu}

\author{Athanase Papadopoulos}
\address{A. Papadopoulos, Institut de Recherche Math{\'e}matique Avanc\'ee,
Universit{\'e} de Strasbourg and CNRS,
7 rue Ren\'e Descartes,
 67084 Strasbourg Cedex, France
 and: CUNY;  Hunter College, 695 Park Ave, New York, NY 10065, USA.} \email{athanase.papadopoulos@math.unistra.fr}
\date{\today}

\thanks{The three authors were partially supported by the GEAR network of the NSF (GEometric structures And Representation varieties) and by the French ANR program FINSLER}

\begin{document}

\begin{abstract}  We survey some major contributions to Riemann's moduli space and Teichm\"uller space. Our report has a historical character, but the stress is on the chain of mathematical ideas. We start with the introduction of  Riemann surfaces, and we end with the discovery of some of  the basic structures of Riemann's moduli space and Teichm\"uller space. We point out several facts which seem to be unknown to many algebraic geometers and analysts working in the theory. 
 The period we are interested in starts with Riemann, in 1851, and ends in the early 1960s, when Ahlfors and Bers confirmed that Teichm\"uller's results were correct.

This paper was written for the book \emph{Lipman Bers, a life in Mathematics}, edited by Linda Keen , Irwin Kra and Rubi Rodriguez (Amercian Mathematical Society, 2015). It is dedicated to the memory of Lipman Bers who 
was above all a complex analyst and spent a large part of his life and energy working on the analytic structure of Teichm\"uller space. His work on analysis is nevertheless inseparable from geometry and topology. In this survey, we highlight the relations and the logical dependence between this work and the works of Riemann, Poincar\'e, Klein, Brouwer, Siegel, Teichm\"uller, Weil, Grothendieck and others. We explain the motivation behind the ideas. In doing so, we point out several facts which seem to be unknown to many Teichm\"uller theorists.

\bigskip


\bigskip


\end{abstract}

\maketitle

  
\tableofcontents
\addtocontents{toc}{\protect\setcounter{tocdepth}{1}}

     \section{Introduction}

Lipman Bers spent a substantial part of his activity working on Teich\-m\"uller space. On the occasion of the Centennial of his birth, we present a report on the origin and the early development of that theory. 
 The period starts with Riemann and ends with Teichm\"uller and the confirmation by Ahlfors and Bers that the statements that Teichm\"uller made (including those for which he did not provide complete proofs) were sound. The main lines of the exposition is guided by the spirit of Riemann, who defined the concept of Riemann surface,  introduced an equivalence relation between such surfaces, and claimed that the number of moduli for the set of equivalence classes of closed surfaces of genus $p\geq 2$ is $3p-3$. A central question, during several decades, was to give a precise meaning to this moduli count and to study the structure of what became known as the \emph{Riemann moduli space}. The mathematicians whose works are involved in this story include, besides Riemann, Weierstrass, Schwarz, Dedekind, Klein, Poincar\'e, Hilbert, Fricke, Koebe, Brouwer, Weyl, Torelli,  Siegel, Teichm\"uller,  Weil, Grothendieck, Ahlfors, Bers, and others. In addition to the history of moduli and Teichm\"uller spaces, this paper contains  an exposition of the birth of some of important chapters in topology which were motivated by the development of the theory of moduli spaces.

One of the multiple factors that led us to write this article is that the origin and history of Teichm\"uller theory is poorly known, probably because this history is complicated. The name of this theory suggests that it formally starts with Teichm\"uller, and indeed, he was the first to give the definition of the space that bears his name, and he studied it extensively. But the question of understanding the totality of Riemann surfaces as a space was addressed by others before him, starting with Riemann. Therefore, we wanted to provide a survey of the works done  on the subject before Teichm\"uller.

Another reason why the origin of the theory is poorly known is that  several articles by leading authors on Riemann surfaces (starting with Riemann himself) are rather arduous to read. In fact, some important results were given with only sketches of proofs.\footnote{The question of rigorous proofs has to be considered in its proper context. Klein, in his \emph{Lectures on Mathematics}, talking about the work of Alfred Clebsch (1833-1872), who was the successor of Riemann at the University of G\"ottingen and whom he describes as ``the celebrated geometer, the central figure, who was one of my principal teachers," writes: ``However great the achievement of Clebsch's in making the work of Riemann more easy of access to his contemporaries, it is in my opinion that at the present time the book of Clebsch is no longer to be considered as the standard work for an introduction to the study of Abelian functions. The chief objection to Clebsch's presentation are twofold: they can be briefly characterized as a lack of mathematical rigor on the one hand, and a loss of intuitiveness, of geometrical perspicuity, on the other [...] The apparent lack of critical spirit which we find in the works of Clebsch is characteristic of the geometrical epoch in which he lived." (\cite{Klein-Lectures} p. 4ff.)}  Furthermore, Riemann, Klein and Poincar\'e  sometimes relied on arguments from physics, and several proofs they gave were considered as unsatisfactory, e.g. those that used the Dirichlet principle or the so-called ``method of continuity" on which we shall comment. One must also bear in mind that Riemann, Poincar\'e, Brouwer and other contributors to this theory were also philosophers, and sometimes their philosophical ideas (on space, on proof, etc.) are intermingled with their mathematics. The positive effect is that besides the mathematical results, the ideas and the methods these authors introduced acted as a motivation and a basis for a great variety of theories that were developed later.
  
  Another reason why we wrote this paper is that the question of figuring out who exactly proved what, and the attribution to each of our mathematical ancestors the exact credit he or she deserves -- especially for works on which we build our own -- is our duty.  We take this opportunity to review some papers of Teichm\"uller. Most of the results they contain are known, but for some of them, the fact that they are due to him is rather unknown.

Beyond all these reasons, we wanted to read the original texts and trace back the original ideas, highlight the most important ones and the connections between them. We find this a rewarding activity. We also went through the published correspondence of Riemann, Weierstrass, Schwarz, Poincar\'e, Klein, Fricke, Brouwer, Weyl and others. Doing history of mathematics is also doing mathematics, if the historian can   penetrate into the depth of the ideas. Our paper also shows how geometry, topology and analysis are one and the same field. 

In the middle of the 1970s, Thurston came with a new vision that completely transformed the field of low-dimensional geometry and topology. His work had  many important developments and a profound impact on the theory of moduli and Teichm\"uller spaces. He introduced a topological compactification of Teichm\"uller space and he used it in the study of the mapping class group. He highlighted the relation between Teichm\"uller space and  the geometrization of  three-manifolds and with the dynamics of rational maps of the  Riemann sphere. His work motivated the later works of Ahlfors and Bers, of McMullen, more recently of Mirzakhani, and of many others.\footnote{The authors of \cite{MSV} write the following (p. 386): ``Bers and Ahlfors continued their work long into retirement, integrating their approach with his new ideas. An NSF grant application by Ahlfors famously contained only one sentence: `I will continue to study the work of Thurston'." The same episode is reported on by other authors.} Thurston conjectured a far-reaching uniformization theory for 3-manifolds, which generalizes the uniformization theory for Riemann surfaces and which culminated in the work of Perelman. All these developments constitute another story that needs a longer article. We also do not touch upon the deformation theory of higher-dimensional complex structures developed by Kodaira, Spencer, Kuranishi, Grauert and others and which was motivated by the deformation theory of surfaces.

At several places in this article, we quote \emph{verbatim} some mathematicians rather than paraphrasing them. Let us start with Poincar\'e's words, from \cite{PK-avenir}\footnote{The article \cite{PK-avenir}, titled \emph{L'avenir des math\'ematiques} (The future of mathematics) is the written version of Poincar\'e's Lectures at the 1908 ICM (Rome), and it was published at several places: \emph{Atti IV Congr. Internaz. Matematici}, Roma, 11 Aprile 1908, pp. 167-182; \emph{Bulletin des sciences math\'ematiques}, 2e s\'erie, 32, pp. 168-190; \emph{Rendiconti del Circolo matematico di Palermo} 16, pp. 162-168; \emph{Revue g\'en\'erale des sciences pures et appliqu\'ees} 19, pp. 930-939; \emph{Scientia} (Rivista di Scienza) 2, pp. 1-23, and as a chapter in the book \emph{Science et m\'ethode} published in 1908.}  (p. 930): ``In order to foresee the future of mathematics, the true method is to examine their history and their present state."\footnote{In this paper, all the translations from the French are ours.} Let us also recall a sentence of Andr\'e Weil, from his commentaries in his \emph{Collected Papers} Vol. II \cite{Weil1979} p. 545: ``There should be almost a book, or at least a nice article, on the history of moduli and moduli spaces; and for that, one should go back to the theory of elliptic functions."
  Our paper is a modest contribution to this wish.

   We tried to make the paper non-technical so that we can  communicate the chain of ideas rather than technical details. Some of the results are highlighted as theorems, when they are easy to state and without heavy notation, so that they can serve as landmarks for the reader.  The history of mathematical ideas does not progress linearly. This is why in this paper there are overlaps and flash-backs.

We are thankful to Linda Keen and Irwin Kra who, after reading our previous historical article \cite{JPH}, proposed that we write more on this history, for the Bers Centennial volume. We also thank Bill Abikoff, Vincent Alberge, Jeremy Gray, Linda Keen, Irwin Kra and Fran\c cois Laudenbach for reading  preliminary versions of this article and providing valuable remarks.

    \section{Riemann}\label{s:Riemann}

 Riemann is one of those mathematicians whose ideas constituted a program for several generations. His work encompasses several fields, in mathematics and physics.  In a letter addressed to Klein on March 30, 1882, Poincar\'e writes about him (\cite{HP-correspondance2} p. 108): 
 \begin{quote}\small 
 He was one of those  geniuses who refresh so well the face of science that they imprint their character, not only on the works of their immediate students, but on those of all their successors for many years. Riemann created a new theory of functions, and it will always be possible to find there the germ of everything that has been done and that will be done after him in mathematical analysis.\footnote{An excerpt from this letter was published in the \emph{Mathematische Annalen} 20 (1882) p. 52-53.}
 \end{quote}
 
  We also quote Klein, from an address he delivered in Vienna on September 27, 1894 whose theme was \emph{Riemann and his importance for the development of modern mathematics} \cite{Klein-Wien} (p. 484): ``Riemann was one these out-of-the way scholars who let their profound thoughts ripen, silently, in their mind."\footnote{In this footnote and in others, we provide a few biographical elements on Riemann and some other mathematicians whose works are mentioned in this paper. It is easy for the reader to find biographies of mathematicians, but we wanted to emphasize, in a few lines, facts which help understanding the background where  the ideas in which we are interested evolved.
  
  Georg Friedrich Bernhard Riemann (1826-1866) was the son of a pastor and he grew up in a pious atmosphere. He arrived in G\"ottingen in 1846 with the aim of studying theology, and indeed, on April 25 of the same year, he enrolled at the faculty of philology and theology. But he soon turned to mathematics, with his father's consent, after having attended a few lectures on that subject. It is interesting to note that this background of Riemann is quite similar to that of Leonhard Euler (of whose works Riemann was a zealous reader). It is also important to recall that besides his interest in theology and in mathematics, Riemann was (like Euler) a physicist and a philosopher. His philosophical  tendencies (in particular his ideas on space and on form) are important for understanding his motivations and the roots of his mathematical inspiration. Riemann's habilitation lecture, \emph{On the hypotheses which lie at the bases of a geometry}, is a mathematical, but also a philosophical essay. The paper \cite{Plotnitsky} is a good introduction to Riemann's philosophical ideas. In Physics, Riemann wrote articles on magnetism, acoustics, fluid dynamics, heat conduction, electrodynamics, and optics. His last written work (which was published posthumously) is an article on the physiology of sound, motivated by works of Hermann von Helmholtz. Euler also wrote on all these subjects. Riemann died of pleurisy just before he attained the age of forty. It should be clear to anyone reading his works that Riemann   wrote only a small part of his ideas. A short and careful biography of him, containing several personal details, was written by Richard Dedekind (1831-1916), his colleague and friend who also took up the task of publishing his unfinished manuscripts after his death.  About half of Riemann's \emph{Collected Mathematical Works} volume  \cite{Riemann-Gesammelte}, which was edited in 1876 by Dedekind and Weber, consists of posthumous works. Dedekind's biography is contained in this \emph{Collected works} edition. For more details, we refer the reader to the more recent biography written by Laugwitz \cite{Laugwitz}.}  
       
 In this section, we shall concentrate on one of Riemann's contributions, namely, the subject that soon after its discovery received the name ``Riemann surfaces." This theory makes a synthesis between fundamental ideas from algebra, topology, geometry and analysis. It is at the heart of modern complex geometry, complex analysis and algebraic geometry. Through elliptic and the other modular functions (Weierstrass elliptic $\frak{p}$-functions, Eisenstein series, etc.), it is also at the foundations of the analytic and the algebraic theories of numbers. Riemann surfaces are also central in mathematical physics.
 
We shall mostly refer to two writings of Riemann in which the concept of Riemann surface plays a prominent role: his doctoral dissertation, \emph{Grundlagen f\"ur eine allgemeine Theorie der Functionen einer ver\"anderlichen complexen Gr\"osse} (Foundations for a general theory of functions of a complex variable) \cite{Riemann-Grundlagen}  (1851), in which he initiated this theory, and his memoir \emph{Theorie der Abel'schen Functionen} (The theory of Abelian functions) \cite{Riemann-Abelian} (1857), which contains important developments. On Riemann's dissertation, Ahlfors writes (\cite{Ahlfors-dev}, p. 4):
\begin{quote}\small
[This paper] marks the birth of geometric function theory [...]
Very few mathematical papers have exercised an influence on the later development of mathematics which is comparable to the stimulus received from Riemann's dissertation. It contains the germ to a major part of the modern theory of analytic functions, it initiated the systematic study of topology, it revolutionized algebraic topology, and it paved the way for Riemann's own approach to differential geometry. [...] The central idea of Riemann's thesis is that of combining geometric thought with complex analysis.
\end{quote}
 
Before talking in some detail about Riemann's contributions, let us make a few remarks on the situation of complex analysis in his time.

The first name that comes to mind in this respect is Cauchy, to whom Riemann's name is intimately tied through the expression ``Cauchy-Riemann equations."    When Riemann started his research, Cauchy's techniques of path integration and the calculus of residues were known to him. He extensively used them and he contributed to their development, especially by using new topological ideas. The development of the theory of complex variables also gave rise to multi-valued functions, like the complex square root, the complex logarithm, the inverses of complex trigonometric functions, etc. and there were also multi-valued functions defined by integrals in the complex plane. To deal with all these functions, Cauchy used a method that involved deleting some curves he called ``cuts" from the complex plane, and he considered specific determinations of the multi-valued function on the resulting surface. Very often, Cauchy simply considered a single branch of the cut (and therefore simply-connected) domain. Riemann came out with a new powerful idea. He started with the cut off surfaces of Cauchy, but he then glued them together to obtain a connected surface, which becomes a ramified cover of the sphere, which he considered as the complex plane to which is added one point at infinity.\footnote{The sphere is, in Riemann's words: ``die ganze unendliche Ebene $A$."} The points in each fiber above a point represented the various determinations of the multi-valued function at that point. 
 
 Another important notion which was available to Riemann is the notion of circle of convergence of a power series, which was also formally introduced by Cauchy,\footnote{For real variables, the notion of analytic function of one variable, defined by power series, was introduced by Lagrange. See e.g. Pringsheim and Faber \cite{Pringsheim}.}  but for which a  better source -- for Riemann -- is probably Gauss' work on the hypergeometric series. The notion of analytic continuation was known to Weierstrass before Riemann, and is present in unpublished papers he wrote in 1841 through 1842.\footnote{These papers are reported on in \cite{BK}. Weyl, in \cite{Weyl1913} p. 1, also refers to an 1842 article by Weierstrass \cite{Weierstrass1842}.} Riemann was aware of this concept when he wrote his thesis in 1851, but it is unlikely that he had access to Weierstrass' notes.\footnote{According J. Gray (in private correspondence), when Riemann started his researches, the concept of analytic continuation was ``informally in the air" before Weierstrass bursts on the scene to everyone's surprise in 1854. See Chapter 6 in \cite{BG} for an exposition of  Weierstrass' work.  The question of whether Riemann discovered it or learned it from someone else is open. Dirichlet might have been aware of that principle, and Riemann followed Dirichlet's lectures in Berlin in 1847-1849. The book \cite{Laugwitz} contains a concise presentation of the notion of power series as used by Weierstrass and Riemann (see in particular p. 83). We know by the way that Weierstrass made a solid reputation in France. Mittag-Leffler recounts in his 1900 ICM talk that when he arrived to Paris in 1873, to study analysis, the first words that Hermite told him were (\cite{ML-W} p. 131): ``You are making a mistake, Sir, you should have followed Weierstrass' courses in Berlin. He is the master of all of us." It is interesting to note that this was shortly after the ravaging French-German war of 1870-1871 which ended with the defeat of France and the unification of the German States under the \emph{German Empire}. Mittag-Leffler adds: ``I knew that Hermite was French and patriotic. At that occasion I learned to what degree he was a mathematician." The esteem for the German school of function theory transformed in the twentieth century, to an esteem for the school there on several complex variables. The following story is told as several places, and in particular by Remmert in  \cite{Remmert1998} p. 222: ``[Henri] Cartan asked his students who wanted to learn several complex variables: Can you read German? If answered in the negative, his advice was to look for a different field." (This was after another war between the two countries.)}
  
Riemann's  construction by cutting and gluing back continued to be used long after him as a \emph{definition} of a Riemann surface. Jordan, for instance, in his \emph{Cours d'analyse de l'\'Ecole Polytechnique}, writes the following (\cite{Jordan-cours}, 3d edition, Vol. III, p. 626):
 \begin{quote}\small
 Let us imagine, with Riemann, a system of $n$ leaves [feuillets] lying down on the plane $P$ of [the variable] $z$. Each of these leaves, such as $P_i$, is cut along the lines $L_1,\ldots, L_v$ which are characterized by the indices $i$. Any one of these lines, $L$, having index $(ik)$, will be a cut for the two leaves $P_i, P_k$. Imagine that we  join each of the boundaries of the cut performed on $P_i$ with the opposite boundary of the cut done on $P_k$. If we perform the same thing for each of the lines $L_1,\ldots,L_v$, the result of all these joins will be to reunite our $n$ leaves in a unique surface.
  \end{quote}

A vey good and concise introduction to the birth of Riemann surfaces is contained in Remmert's paper \cite{Remmert1998} which concerns more generally the birth of complex manifolds. In the rest of this section, we give a brief account of the main contributions of Riemann that are related to the concept of Riemann surface. 
   
      \medskip

\noindent{\bf (1) Riemann surfaces and function theory.} Riemann's construction of surfaces as branched covers of the sphere is contained at the beginning of his thesis \cite{Riemann-Grundlagen}. Riemann's construction of surfaces as branched covers of the sphere was a matter of following the consequences of defining and studying complex functions on arbitrary two-dimensional surfaces. Let us quote from the beginning of his thesis: ``We restrict the variables $x,y$ to a finite domain by considering as the locus of the point $O$ no longer the plane itself but a surface $T$ spread over the plane. We admit the possibility ... that the locus of the point $O$ is covering the same part of the plane several times."\footnote{Translation in \cite{Remmert1998}, p. 205. In this paper, we are using Remmert's translations when they are available.} As a consequence of this construction, Riemann discovered surfaces on which multi-valued complex functions become single-valued.\footnote{Riemann, instead of ``single-valued," uses the expression ``perfectly determined." Weyl, in \cite{Weyl1913}  (p. 2) uses the word ``uniform."} We recall the idea, because this was one of his important achievements. We start with an analytic function, and to be precise, we shall think of the function $\sqrt{z}$. Consider a point $z_0$ in the complex plane which is different from the origin, so that the function is multivalued at that point, and take some determination of this function in a neighborhood of $z_0$. Continue defining this function along paths starting at the point $z_0$, using the principle of analytic continuation. For paths that end at $z_0$ with odd winding number with respect to the origin, the value we get at $z_0$ is different from the initial one. Here Riemann introduced the idea that the endpoint of such a path should not be considered as $z_0$, but a point on a different \emph{sheet} of a new surface on which the function $\sqrt{z}$ is defined. This introduced at the same time the notion of covering. In the example considered, the surface obtained is a two-sheeted branched cover of the complex plane (or of the sphere), and the branching locus is the origin. This construction is very general, that is, it is possible to associate to any multi-valued analytic function a Riemann surface which is a branched cover of the sphere and on which the function is defined and becomes single-valued. Let us consider an arbitrary  multivalued function $w$ of $z$ is defined by an equation 
 \begin{equation}\label{eqn:fcn}
 P(w,z)=0
 \end{equation}
 where $P$ is a polynomial (in the above example, the polynomial is  $P(w,z)=w^2-z$), or more generally an arbitrary analytic function (that is, the function  $w$ of $z$ may also be transcendental). In the polynomial case, and when $P$ is irreducible of degree $m$ in $w$, then for each generic value of $z$ we have $m$ distinct values of $w$, each of them varying continuously (in fact, holomorphically) in $z$. Thus, by the principle of analytic continuation, by following paths starting at $z$, at each generic point $z_0$ there are $m$ possible local analytic developments of $w$. A Riemann surface is constructed, which is naturally a branched covering of the plane. The degree of the covering is the number of values of $w$ for a generic value of $z$. This construction is described for the first time in Riemann's dissertation \cite{Riemann-Grundlagen} and is further developed in the section on preliminaries in his 1857 paper \cite{Riemann-Abelian}. Riemann went much further in the idea of associating a Riemann surface to a multi-valued function, as he considered Riemann surfaces as objects in and of themselves on which function theory can be developed in the same way as the theory of complex functions is developed on the complex plane. He had the idea of an abstract Riemann surface, even though his immediate followers did not. Indeed, for several decades after him, geometers still  imagined surfaces with self-intersections, first, immersed in 3-space, and, later, in $\mathbb{C}^2$.\footnote{In his \emph{Idee der Riemannschen Fl\"ache} (\cite{Weyl1913} p. 16), Weyl, writes about these spatial representations, that ``in essence, three-dimensional space has nothing to do with analytic forms, and one appeals to it not on logical-mathematical grounds, but because it is closely associated with our sense-perception. To satisfy our desire for pictures and analogies in this fashion by forcing inessential representations on objects instead of taking them as they are could be called an anthropomorphism contrary to scientific principles."} For instance, in 1909, Hadamard, in his survey on topology titled \emph{Notions \'el\'ementaires sur la g\'eom\'etrie de situation}, talking about Riemann surfaces, still considers lines along which the leaves cross themselves (cf. \cite{Hadamard-Notions} p. 204). From the modern  period, we can quote Ahlfors (\cite{Ahlfors-dev}, p. 4):
\begin{quote}\small
Among the creative ideas in Riemann's thesis none is so simple and at the same time so profound as his introduction of multiply covered regions, or Riemann surfaces. The reader is led to believe that this is a commonplace construction, but there is no record of anyone having used a similar device before. As used by Riemann, it is a skillful fusion of two distinct and equally important ideas: 

\noindent 1) a purely topological notion of covering surface, necessary to clarify the concept of mapping in the case of multiple correspondence; 

\noindent 2) an abstract conception of the space of the variable, with a local structure defined by a uniformizing parameter.
\end{quote}

Riemann also addressed the inverse question, viz., that of associating a surface to a function. More precisely, he considered the question of whether an arbitrary surface defined by a system of leaves made out of pieces of the complex plane glued along some system of lines can be obtained, by the process of  analytic continuation, from some algebraic curve defined by an equation such as (\ref{eqn:fcn}). This can be formulated as the question of finding a meromorphic function with prescribed position and nature  of singularities (poles and branch points). Riemann answered this question, and in doing so, he was the first to emphasize the fact that a meromorphic function is determined by its singularities. This result was one of his major achievements. It is sometimes called the \emph{Riemann existence theorem}, and it is  the main motivation of Riemann--Roch theorem which we discuss below. Picard, in his \emph{Cours d'Analyse}, describing this problem, writes (\cite{Picard-cours}, Tome II, p. 459): ``We enter now in the profound thought of Riemann." Riemann's existence proof was not considered as rigorous, because of a heuristic use he made of the Dirichlet principle. The proof was made rigorous by others, in particular Neumann \cite{Neumann-Vorlesungen} and Schwarz \cite{Schw}. Picard reconstructed a complete proof of this result, following the ideas of Neumann and Schwarz, in Chapter XVI of his \emph{Trait\'e d'Analyse}  \cite{Picard-cours}. The correspondence between Riemann surfaces and classes of algebraic functions is essential in Riemann's count of the number of moduli, which we shall review below.

   \medskip
   
\noindent{\bf (2) The Riemann--Roch theorem.}
 The Riemann--Roch theorem originates in Riemann's article on Abelian functions \cite{Riemann-Abelian}. It makes a fundamental link between topology and analysis. It is formulated as an equality involving the genus of a Riemann surface and the zeros and poles of meromorphic functions and of meromorphic differentials defined on it. A corollary of the theorem is that any compact Riemann surface admits a non-constant meromorphic function, that is, a non-constant map into the projective space, with control on its degree. Half of theorem (that is, one inequality) was given by Riemann in \S 5 of his article \cite{Riemann-Abelian}.\footnote{This is the inequality $d\geq m-p+1$, where $d$ is the complex dimension of the vector space of meromorphic functions having at most poles of first order at $m$ given points; cf. p. 107-108 of \cite{Riemann-Abelian}.} The other half is due to Riemann's student Roch \cite{Roch}.\footnote{Gustav Roch (1839-1866) died at the age of 26 from tuberculosis.} The theorem was given the name Riemann--Roch by Alexander von Brill and Max Noether in their 1874 article \cite{Brill-Noether}. Riemann obtained this result at the same time as the existence proof for complex functions satisfying certain conditions which we mentioned above. After the discovery of the gap we mentioned in Riemann's argument, several mathematicians continued working on the Riemann--Roch theorem. Dedekind and Weber gave a new proof in 1882 \cite{DW1882} and Landsberg gave another one in 1898 \cite{Landsberg1898a}. Hermann Weyl wrote a proof in his famous 1913 book \cite{Weyl1913}, including what he called the \emph{Ritter extension}. Teichm\"uller gave a new proof in his paper \cite{T25}, and there is also a section on that theorem in his paper \cite{T20}, where this theorem is used in an essential way.  There are several later generalizations of the Riemann--Roch theorem, and one of them is due to Andr\'e Weil, whom we shall quote thoroughly later in this paper. In his comments on his \emph{Collected papers} edition (\cite{Weil1979} Vol. I p. 544), Weil writes that for a long period of time, the Riemann--Roch theorem was one of his principal themes of reflection. In 1938, he obtained a version of a theorem he called the \emph{non-homogeneous Riemann--Roch theorem} (cf. \cite{Weil1938a}) which is valid for curves over arbitrary algebraically closed fields.\footnote{In fact, this generalized Riemann--Roch theorem was already sketched in Weil's 1935 paper \cite{Weil1935a}.} In 1953, Iwasawa published a new version of Riemann--Roch, expressed in terms of the newly discovered tools of algebraic topology \cite{Iwasawa1953}. Among the many other developments, we mention  the famous Riemann--Roch-Hirzebruch theorem which is a version of Riemann--Roch valid for complex algebraic varieties of any dimension, and its generalization by Grothendieck to the so-called Riemann--Roch-Grothendieck theorem expressed in the language of coherent cohomology. The developments and consequences of these results are fascinating, leading for example to the Atiyah-Singer index theorem which is a unification of Riemann--Roch with theorems of Gauss-Bonnet and of Hirzebruch.
A historical survey of the (classical) Riemann--Roch theorem was written by Jeremy Gray \cite{Gray-RR}.

   \medskip
   
\noindent{\bf (3) Birational equivalence between algebraic curves.}  The identification discovered by Riemann between a Riemann surface and the field (or, using more modern tools, the sheaf of germs) of meromorphic functions it carries, in the case of compact surfaces, leads to the following equivalence relation: \emph{two Riemann surfaces $S_1$ and $S_2$ associated with algebraic equations of the form (\ref{eqn:fcn}) are considered as equivalent if there exists a \emph{birational equivalence} between the curves defining them.}  

 This identification is an expression of the fact that  fields of meromorphic functions associated with a given algebraic curve are isomorphic. In the compact case, we know that this holds if and only if the associated Riemann surfaces are conformally equivalent. 
 Thus, Riemann was led to the problem of finding birational invariants of algebraic curves. This was the original form of the famous ``moduli problem." This search for invariants under birational transformations became, after Riemann's work, a subject of intense research, starting with the work of Cremona in 1863 and developed by the so-called Italian school.

   \medskip
   
\noindent{\bf (4) The number of moduli.}  Riemann, at several places, stated that there are $3p-3$ complex ``moduli" associated with a closed surface of genus $p\geq 2$. This parameter count is at the origin of the notion of \emph{moduli space}. 

   One of Riemann's counts of the moduli  is based on the association between equivalence classes of algebraic functions and equivalence classes of Riemann surfaces. This count involves the number of sheets of a branched cover associated with a curve and local information around the ramification points.\footnote{This has been made precise by Hurwitz, and it is at the origin of the so-called \emph{Riemann-Hurwitz formula}, cf.  \cite{Hurwitz39}, p. 55.} Riemann ended up with $3p-3$ moduli, a number which occurs in several places in his, works based on heuristic arguments. Several authors, in the years following Riemann's writings, gave more detailed proofs of these counts,  among them Klein \cite{Klein-Riemann} and Picard \cite{Picard-cours}. But the real meaning of the number of moduli as the complex dimension of a manifold (or an orbifold), or at least as as a number associated with some structure, was still missing, and it is hard to find a proper definition of the term ``moduli" in these earlier works.  Hyman Bass and Irwin Kra write in \cite{KB} (p. 210):
   \begin{quote}\small
   One of their important open problems, the moduli problem, left to us from the nineteenth century, was to make rigorous and precise Bernhard Riemann's claim that the complex analytic structure on an closed surface with $p\geq 2$ handles depends on $3p-3$ complex parameters.
   \end{quote}

Another one of Riemann's counts occurs in his memoir on the theory of Abelian functions  \cite{Riemann-Abelian}, when he talks about the so-called \emph{Jacobi inversion problem}. This is the problem of finding an inverse to the construction of Abelian integrals of the first kind, that is, functions of the form $u(z)=\int^z R(\zeta,w)d\zeta$ where $R$ is a rational function and where $\zeta$ and $w$ are complex variables related by a polynomial equation  $P(\zeta,w)=0$.
Riemann counted the parameters for a certain class of functions.\footnote{According to Rauch \cite{Rauch-Transcendental} p. 42, the term ``moduli" chosen by Riemann comes from this count of moduli of elliptic integrals of the first kind.} He writes (English translation \cite{Riemann--English} p. 93): 
  \begin{quote}\small
  The Jacobi inversion functions of $p$ variables are expressed using $p$-infinite theta functions.\footnote{The name ``theta functions" used by Riemann originates in Jacobi's writings, who used the letter theta to denote them.} A certain class of theta functions suffices. This class becomes special for $p>3$, and in that case we have $p(p+1)/2$ quantities, with $(p-2)(p-3)/2$ relations between them. Thus there remain only $3p-3$ which are arbitrary.
  \end{quote} 
It is interesting to note that the quantity $p(p+1)/2$ that appears in the above quote is the complex dimension of the Siegel modular space to which we shall come back in \S \ref{s:Siegel}. In modern terms, this result says that the image of the so-called Jacobi map (which is called the Jacobian locus) is of dimension $3p-3$. At the same time, Riemann states that theta functions that arise from periods of compact Riemann surfaces satisfy some special relations, and therefore, that the most general theta functions are not needed in the Jacobi inversion problem. This was later considered by Klein as a ``marvelous result" (\cite{Klein-Wien} p. 495).\footnote{In the same text, Klein asks about the role played by the general theta functions. He says that according to Hermite, Riemann already knew the theorem published later by Weierstrass and  also treated by Picard and Poincar\'e, saying that the theta series suffices for the representation of the most general periodic functions of several variables.} This consideration also suggests the so-called \emph{Schottky problem} on characterizing periods,
or the locus of the Jacobian map.\footnote{This locus is the subset of Jacobian varieties of compact
Riemann surfaces in the Siegel modular variety which parametrizes the space
of principally polarized Abelian varieties. The Schottky problem was first studied explicitly by Torelli, cf. \S \ref{s:Torelli} below.}

 In  \S 12  of the same memoir (p. 111 of the English translation \cite{Riemann--English}), Riemann ends up with the same count by computing the number of branch points of  Riemann surface coverings that he constructs.  The quantity $3p-3$ appears also in Riemann's  posthumous article,  \emph{On the theory of Abelian functions} (p. 483-504 of the English translation of the collected papers  \cite{Riemann--English}). At the end of this paper, Riemann studies the special case of genus 3, and he concludes that ``all Abelian functions, with all their algebraic relations, can be expressed via $3p-3=6$ constants, which one can regard as the moduli of the class for the case $p=3$."
 
In conclusion, in the work of Riemann, the quantity $3p-3$  appears several times as a moduli count, but never as a dimension.  The ``number of moduli" had the vague meaning of a ``minimum number of coordinates." In fact, at that time, there was no notion of dimension (and especially no notion of complex dimension) with a precise mathematical meaning, except for linear spaces. We shall discuss the question of dimension in  \S \ref{s:Brouwer}, while we survey Brouwer's work.\footnote{One can find a description of the quantity $3p-3$  that uses the word ``dimension" in Klein's writings (see \S \ref{s:Klein} below), even though there was no precise meaning of that word.} The interpretation of the number of moduli as the complex dimension of a complex space came later, when Teichm\"uller equipped Teichm\"uller space with a complex structure.\footnote{\label{F:Teich} It is interesting to note here the following from the beginning of Teichm\"uller's paper \cite{T32}: ``It has been known for a long time that the classes of conformally equivalent closed
Riemann surfaces of genus $g$ depend on $\tau$ complex constants, where (1) $\tau=0$ if $g=0$; (2) $\tau=1$ if $g=1$, and (3) $\tau=3g-3$  if $g>1$.
This number $\tau$ has been obtained using different heuristic arguments, and the result is
being passed on in the literature without thinking too much about the meaning of this statement." We shall comment on Teichm\"uler's paper in \S \ref{s:Teich} below. We note by the way that Thurston (like many others before him) makes a similar heuristic count in his Princeton Notes \cite{Thurston1976}:  A Fuchsian group representing a closed Riemann surface of genus $g$ is generated by $2g$ transformations satisfying one relation. Each generator, being a fractional linear transformation, has 3 unnormalized real parameters. Hence (subtracting 3 parameters fixed by the relation and 3 more for a normalization of the entire group) we obtain $3\times 3g-3-3= 6g-6$, which is the number of \emph{real} parameters corresponding to Riemann's count. A similar  count is also made in Bers' ICM paper \cite{Bers-Spaces1960} and in Abikoff's \cite{Ab-LNM}.} 
 
    \medskip
    
\noindent{\bf (5) Abelian integrals.} 
In his memoir on Abelian functions \cite{Riemann-Abelian}, after  studying branched coverings of the sphere, Riemann considers Abelian functions, also called Abelian integrals. In his work, these  integrals represent the transcendental approach to functions, as opposed to the algebraic. Abelian integrals  are generalization of elliptic integrals, which were so called because a special case of them gives the length of an arc of an ellipse.\footnote{Elliptic integrals of a real variable were very fashionable since the eighteenth century. They were studied in particular by Fagnano, Euler, and above all, Legendre. The name \emph{elliptic function} was coined by Legendre, who gave it to some integrals that satisfy an  additivity property discovered by Fagnano, concerning the length of the lemniscate. Legendre spent a large effort studying these functions, and he wrote a treatise on the subject, the \emph{Trait\'e des fonctions elliptiques et des int\'egrales eul\'eriennes} \cite{Legendre}. A wide generalization of indefinite integrals satisfying such an additivity property was discovered soon after by Abel. These functions were called later on \emph{Abelian}, presumably for the first time by Jacobi.} From another point of view, Abelian functions associated with Riemann surfaces of genus $p \geq 1$ are generalizations of  the elliptic functions associated with  Riemann surfaces of genus $1$. Abelian functions can also be expressed as ratios of homogeneous polynomials of theta functions. 

Riemann's work on Abelian functions is a continuation of works of Niels Henrik Abel (1802-1829) and Carl Gustav Jacob Jacobi (1804-1851). These works were done during the first half of the nineteenth century, but the theory of path integration in the complex plane, developed by Cauchy,  was unknown to these authors. An important step in the theory is Abel's discovery \cite{Abel1} that an arbitrary sum of integrals, with arbitrary limits, of a given algebraic function, can be expressed as the sum of a certain fixed number of similar integrals with a certain logarithmic expression. Riemann showed that this number is just the \emph{genus} of the curve.\footnote{For this and other notions of genus, we refer to the comprehensive paper \cite{Popescu} by Popescu-Pampu.} In 1832, Jacobi formulated the famous \emph{inversion problem for Abelian integrals}, a problem which many ramifications which generalizes the inversion problem for elliptic integrals and to which Weierstrass and Riemann contributed in essential ways. While Cauchy had failed to grasp a theory of complex integrals with multi-valued integrands, Riemann used Cauchy's path integral theory and associated with Abelian integrals Riemann surfaces, as he did for the algebraic multi-valued functions. He attacked the problem of their classification and of their moduli, and he developed the theory of periods of these integrals. For a clear exposition of the theory of Abelian integrals, with interesting historical remarks, we refer the reader to the monograph \cite{Mark}.  
    \medskip
 
\noindent{\bf (6) The Jacobian.} 
  The theory of the Jacobian asks for the characterization of Abelian integrals, and more generally, for curves, by their period matrices. It uses the integration theory of holomorphic forms on surfaces which is a generalization of Cauchy's computations of path integrals. To make a complex function on a Riemann surface defined as the integral of a differential form well defined, Riemann cuts the surface along a certain number of arcs and curves, and as a result he gets the notion of \emph{periods} of such integrals. In this setting, a period is the difference of values taken by a function when one traverses a cut.  The cuts may be arcs (as in Riemann's dissertation \cite{Riemann-Grundlagen}) or closed curves (in his Abelian functions memoir \cite{Riemann-Abelian}).\footnote{It is possible that Riemann got the idea of cuts from Gauss. Betti wrote in a letter to Tardy that Gauss, in a private conversation, gave to Riemann the idea of cuts. The letter is reproduced in Pont \cite{Pont1974} and Weil \cite{Weil1979a}.}
 Riemann claimed that the number of curves needed to cut the surface into a simply connected region (all the curves passing  through a common point) is of the form $2p$, where $p$ is the genus.\footnote{This is beginning and the heart of the classification of surfaces.} He obtained $2p$  curves which he used as a basis for the ``period parameters" for holomorphic differentials  on the surface. His argument was improved by Torelli, whose name is now attached to the theory. Riemann also introduced other kinds of parameters for these differentials, namely, the locations of poles or logarithmic poles. The idea of reconstructing the curve from its set of periods is at the origin of the embedding of the Torelli space\footnote{We recall that the Torelli space is the quotient of Teichm\"uller space by the Torelli group, that is, the subgroup of the mapping class group consisting of the elements that act trivially on the first homology group of the surface.} in the Siegel space, which was developed later, and it is also a prelude to the idea of equipping Teichm\"uller space with a complex structure using the period map. Several problems known as ``Torelli-type problems," arose from Riemann's work on periods. We shall mention some of them in \S \ref{s:Torelli}. 
  
In another direction, the problem of characterizing general Abelian varieties among all complex tori by period matrices, which was also motivated by the work of Riemann, became an important research subject in the twentieth century and involves works of Siegel, Weil, Hodge and several others.

    \medskip

\noindent{\bf (7) The Riemann theta functions.}  Given a complex vector $z\in \mathbb{C}^g$ and a complex $g\times g$ symmetric matrix $F$ whose imaginary part is symmetric  positive definite,\footnote{A matrix with these properties is usually called a \emph{Riemann matrix}.} the Riemann theta function is defined by the formula
\[\Theta(z,F) = \sum_{n\in \mathbb{Z}^{g}}e^{2\pi i (\frac{1}{2} n^T  Fn+n^T z)}\]
where $n^T$ denotes the transpose of $n$.\footnote{This way of writing the Riemann theta function is due to  Wirtinger \cite{Wirtinger} (1895). Riemann writes these functions differently, see \cite{Riemann-Abelian}, p. 93 and 119ff. of his \emph{Collected papers} edition. The fact that $\mathrm{Im}(F)$ is positive definite guarantees that the series converges for all values of $z$, and that the resulting function is holomorphic in both $z$ and $F$.} These functions are Fourier series-like functions. Riemann associated such functions with arbitrary Riemann surfaces and he studied them in the second part of his 1857 memoir on Abelian functions \cite{Riemann-Abelian}, in his solution of the generalized Jacobi inversion problem. The Jacobi theta function is a special case where $n=1$ and in this case $F$ belongs to the upper half-plane. Jacobi introduced this function in his study of elliptic functions.  Riemann's theta functions are $n$-dimensional generalizations of the Jacobi theta functions in much the same way as the Siegel half-space, which was defined much later, is a multi-dimensional generalization of the upper half-plane. This is another instance, in Riemann's work, of a generalization from the case of genus 1 to higher genus.  
 
Prior to Riemann, Gustav Adolph G\"opel (1812-1847)  and Johann Rosenhain (1816-1887) worked on theta functions.  In 1866, the Academy of Sciences of Paris announced a contest on the inverse problem of Abelian integrals in the case of curves of genus two. G\"opel and  Rosenhain independently solved the problem. G\"opel did not submit his proof to the Academy, but Rosenhain, who was a student of Jacobi, did. Both G\"opel and Rosenhain used 2-variable theta functions in their solution. G\"opel's paper was published in 1847 \cite{Goepel}. Rosenhain's paper won the prize in 1851 and his paper was published the same year \cite{Rosenhain}.

 There are many ramifications of Riemann's theta functions. They play a major role in modern Teichm\"uller theory, and also in number theory. Solutions of differential equations can often be written in terms of Abelian functions, and thus in terms of theta functions. 
 
    \medskip
    
\noindent{\bf (8) Topology of surfaces.}  
 At the beginning of the memoir on Abelian functions \cite{Riemann-Abelian}, there is a section titled \emph{Theorems of analysis situs for the theory of the integrals of a complete differential with two terms}, in which Riemann recalls what is the field of topology, and he introduces new tools in the topology of surfaces. Among these notions, we mention   \emph{connectivity} (in fact, arcwise-connectivity) and  \emph{simply-connected surface}. For non simply-connected surfaces, he introduced the \emph{order of connectivity} which gave rise to the notion of  genus for closed surfaces.\footnote{Riemann counted the so-called ``number of handles" of a surface, that is, the minimal number of curves along which one cuts it to make it planar. The word ``genus" (in German: ``Geschlecht") was introduced by Clebsch in \cite{Clebsch1} p. 43. Poincar\'e talked about  the \emph{first Betti number}, which, for closed oriented surfaces, is two times the genus, and can be generalized to higher-dimensional manifolds. Poincar\'e developed the topological tools needed to compute these numbers. } As the reader may have realized now, these topological ideas were not isolated from the rest of Riemann's work and were not developed \emph{per se}, but were part of the  development of the theory of Riemann surfaces.

The whole subject of topology of surfaces was influenced by the research on Riemann surface theory. The work of Riemann also gave a fundamental impetus for the general definition of manifold. Ahlfors writes (\cite{Ahlfors-dev} p. 4): 
\begin{quote}\small
From a modern point of view the introduction of Riemann surfaces foreshadows the use of arbitrary topological spaces, spaces with a structure, and covering spaces.
\end{quote}

   \medskip
   
\noindent{\bf (9) The use of potential theory (the Dirichlet principle).} By his substantial use of the Dirichlet principle\footnote{The Dirichlet principle has its origin in physics.  Riemann called this principle the ``Dirichlet principle" because he learned it from Dirichlet, during his stay in Berlin, in the years 1847-1849. Pierre Gustave Lejeune Dirichlet (1805-1859) was born in D\"uren, a city on the Rhine, now in Germany, but which belonged then to the French Empire. His paternal grandfather was Belgian.
During a long stay he made in Paris, between 1822 and 1827, the young Dirichlet became acquainted with several French mathematicians, and after his return to Germany, he contributed to making the link between the French and the German mathematicians. There are no existing notes from the course Dirichlet gave in Berlin, but we know that in 1854, Riemann chose, as the subject of his first lectures in G\"ottingen,  the theory of partial differential equations and their applications to physics, and that his course was modelled on the one that Dirichlet gave in Berlin (Dedekind's biography, p. 527 of \cite{Riemann--English}). Dirichlet moved to the University of G\"ottingen in 1855, where he became the successor of Gauss. He stayed there until his death in 1859. Riemann became, in 1859, the successor of Dirichlet.} and the important place he gave to the Laplace's equation, Riemann introduced potential theory as an important factor in complex function theory. To quote Ahlfors (\cite{Ahlfors-dev} p. 4): ``Riemann virtually puts equality signs between two-dimensional potential theory and complex function theory."

Brill and Noether give a mathematical and historical review of this principle in their paper \cite{BN-D}. A fact which is often repeated in the literature is that Weierstrass pointed out a gap in Riemann's proof of his existence theorem, related to his use of the Dirichlet principle. Riemann's argument was indeed poor, and in fact he was repeating unproved arguments of Gauss and Green. But he also attempted  a proof of these arguments -- a visionary extension from functions to functionals -- and the fact that he saw that these claims needed a proof is already an important step. Ahlfors writes about this issue (\cite{Ahlfors-dev} p. 5):  ``It is perhaps wrong to call Riemann uncritical, for he made definite attempts to exclude a degenerating extremal function." According to Gray \cite{Gray1994}, it was Prym (who was one of Riemann's students) and not Weierstrass who first pointed out the gap in Riemann's argument and described it in a short paper \cite{Prym1871}.  Weierstrass' objection concerned not only the use of the Dirichlet  principle by Riemann, but also its use by Green, Gauss and others. The use by Riemann of this variational principle is in the context of the definition of a holomorphic function $f(x+iy)=u+iv$ that he gives at the beginning of his memoir on Abelian functions \cite{Riemann-Abelian} (p. 79 of the English translation \cite{Riemann--English}), as a function satisfying the so-called Cauchy-Riemann equations.\footnote{Note however that at the beginning of his dissertation
\cite{Riemann-Grundlagen} (p. 2 of the English translation \cite{Riemann--English}), Riemann gives a definition in terms of conformality. He writes: ``In whatever way $w$ is determined from $z$ by a combination of simple operations, the value of the derivative $\frac{dw}{dz}$ will always be independent of the particular value of the differential $dz$," and after a small computation, he shows that this property ``yields the similarity of two corresponding infinitely small triangles."}  In brief, Weierstrass gave a counter-example that destroyed the  claim that functionals bounded below attain their lower bounds and Prym destroyed the specific claim about the Dirichlet principle, which was much more worrying for Riemann's function theory. This did not prove that Riemann's claim was false, but showed that it needed a formal proof. Riemann was not bothered by Weierstrass' criticism, and for what regards the reaction of the mathematical community, Remmert writes  (\cite{Remmert1998} p. 212)  that ``Weierstrass' criticism should have come as a shock, but it did not. People felt relieved of the duty to learn and accept Riemann's methods." Weierstrass presented his counter-example on July 14, 1870 before the Royal Academy of Sciences in Berlin, in a note titled \emph{On the so-called Dirichlet principle} \cite{W14-1870}. This was four years after Riemann's death. Neuenschwande reports in \cite{NeuenschwanderII} that ``after Riemann's death,  Weierstrass attacked his methods quite often, in part even openly."\footnote{In fact, there seems to be an incompatibility between Riemann's and Weierstrass' views on function theory. Remmert writes in \cite{Remmert1998} p. 15: ``It was Weierstrass' dogma that function theory is \emph{the theory of convergent Laurent series} ... Integrals are not permitted. The final aim is always the representation of functions. Riemann's geometric yoga with paths, cross-cuts, etc., on surfaces is excluded, because it is inaccessible to algorithmization."} Weierstrass' objection to Riemann's argument is analyzed in several papers, e.g. \cite{NeuenschwanderII}  where the author writes that after the objections by Weierstrass on Riemann's foundational methods, ``only with the work of Klein  and the rehabilitation of the Dirichlet principle by Hilbert (\cite{HD1} \cite{HD2}) could Riemannian theory again gradually recover from the blow delivered to it by Weierstrass." 

 The history of the Dirichlet principle is complicated. The title of the book \cite{Monna}, ``Dirichlet's principle: A mathematical comedy of errors and its influence on the development of analysis,"  is revealing.  After Weierstrass' criticism, the ``Dirichlet principle" turned into a ``Dirichlet problem."  The problem was solved later by works of several people, culminating with Hilbert who, in 1900-1904 (\cite{HD1} \cite{HD2}), gave a proof of a form of the Dirichlet principle which was sufficient for its use under some differentiablity conditions on the boundary of the domain that were made by Riemann. Hilbert utilized for that the so-called direct methods of the calculus of variations. The same question was also dealt with, in different ways, by Courant \cite{Courant-D} and by Weyl who, in his famous book \cite{Weyl1913} \S 19,  gave a new proof of the uniformization following a method suggested by Hilbert and others. For a concise presentation of the Dirichlet problem, we refer the reader to \cite{Garding}. We also refer the reader to \S 5.2.4 of the book \cite{BG} for the use Riemann made of the Dirichlet principle, and a discussion of where his approach breaks. 
  
  Let us conclude by quoting Hermann Weyl, from a review of mathematics in his lifetime, quoted in \cite{Century}, Part II, p. 326: ``In his oration honoring Dirichlet, Minkowski spoke of the true Dirichlet principle, to face problems with a minimum of blind calculation, a maximum of seeing thought."

   \medskip
   
\noindent{\bf (10) The Riemann mapping theorem.}  
The theorem is proved in the last three pages of Riemann's dissertation \cite{Riemann-Grundlagen}. In Riemann's words, the theorem says the following (English translation \cite{Riemann--English}, p. 36):

 \begin{theorem}\label{th:RM}
 Two given simply connected plane surfaces can always be related in such a way that each point of one surface corresponds to a point of the other, varying continuously with that point, with the corresponding smallest parts similar. One interior point, and one boundary point, can be assigned arbitrary corresponding points; however, this determines the correspondence for all points. 
  \end{theorem}
   The expression ``with the corresponding smallest parts similar" means that the map is conformal.\footnote{\label{f:cart} This notion of conformality was already used by Lambert, Gauss and others before Riemann, in their works on conformal projections of the sphere in the setting of cartography (geographic map drawing).} Riemann proves the theorem under the hypothesis of piecewise differentiability of the boundary. He proves the existence of a conformal map between a simply connected open subset of the plane which is not the entire plane and the unit disc by using the Dirichlet boundary value problem for harmonic maps. This is where he appeals to the so-called \emph{Dirichlet principle}, which characterizes the desired function among the functions with the given boundary values as the one which minimizes the energy integral. At the end of his dissertation, Riemann mentions a wide generalization of his theorem. He writes (\cite{Riemann--English} p. 39):
  \begin{quote}\small
  The complete extension of the investigation in the previous section to the general case, where a point of one surface corresponds to several points on the other, and simple connectedness is not assumed for the surface, is left aside here. Above all this is because, from the geometrical point of view, our entire study would need to be put in a more general form. Our restriction to plane surfaces, smooth except for some isolated points, is not essential: rather, the problem of mapping one arbitrary given surface onto another with similarity in the smallest parts, can be treated in a wholly analogous way. We content ourselves with a reference to two of Gauss' works: that cited in Section 3, and \emph{Disquisitiones generales circa superficies curvas}, \S 13.\footnote{This is the problem of mapping conformally a surface onto a plane, which is related to geographic map drawing, which we mentioned in Footnote \ref{f:cart}.}
    \end{quote}
   
 After the gap that was found in Riemann's argument, related to his use of the Dirichlet principle, several mathematicians tried to find new methods to prove the theorem, and this became the topic of a new research activity.  Among the developments, we mention the result of Clebsch (1865) saying that every compact Riemann surface of genus zero is biholomorphic to the Riemann sphere and every compact Riemann surface of genus one is biholomorphic to the quotient of the plane by a lattice, see \cite{Clebsch1} and \cite{Clebsch2}. 
 Hermann Amadeus Schwarz in 1870 \cite{Schwarz1870} gave a proof of the Riemann mapping theorem for simply connected plane domains whose boundary consists of a finite number of curves.\footnote{This special case is important; it applies for instance to the polygons in the hyperbolic plane which play a major role in the works of Klein and Poincar\'e.} He introduced a method called the ``alternating method," an iterative method which solves the boundary-value problem in Riemann's proof for domains which are unions of domains for which the problem is known to be solved. This method is also reported on by Klein in \cite{Klein-Vorlesungen}.\footnote{In a letter to Poincar\'e, dated July 9, 1881, Klein writes: ``There is no doubt that the Dirichlet principle must be abandoned, because it is not at all conclusive. But we can completely replace it by more rigorous methods of proof. You can see an exposition in more detail in a work by Schwarz which I precisely examined these last days (for my course)" (see \cite{HP-correspondance2}  p. 101) and Klein refers to the 1870 article by Schwarz \cite{Schwarz-Berliner}. Ten years after Schwarz's work, Klein was still referring to his approach to the problem.} In his \emph{Lectures on Riemann's theory of Abelian integrals} \cite{Neumann-Vorlesungen}, Neumann gave another proof of the Riemann mapping theorem for simply connected plane domains. Poincar\'e, in his 1890 article \cite{PK-balayage} gave a method of approaching the solution by a sequence of functions which are not harmonic but which have the required boundary values. This method became famous under the name ``m\'ethode du balayage." William Fogg Osgood (based on ideas of Poincar\'e and Harnack), in \cite{Osgood1900}  (1900), gave an existence theorem  for a Green function on bounded plane domains with smoothness conditions on the boundary. Schwarz \cite{Sch} (1869) and Elwin Bruno Christoffel \cite{Chr} (1867 and 1871) independently, gave explicit Riemann mappings of open sets whose boundaries are polygons.\footnote{The uniformization of polygonal regions plays a major role in the works of Klein and Poincar\'e. In their work on Fuchsian groups, the fundamental domain of the action is a polygon in the hyperbolic plane. Such a polygon corresponds naturally to the cut up surface of Riemann. One gets a Riemann surface as the quotient of the fundamental plane by the group action  that identifies the boundary components. One form of the uniformization problem is whether any Riemann surface can be obtained in this way. This is where the so-called \emph{continuity method}, which we mention several times in this paper, was used. In his letter to Poincar\'e dated July 2, 1881, Klein writes, considering the uniformization of polygonal regions (\cite{HP-correspondance2} p. 99): ``Weierstrass can determine effectively the constants involved by convergent processes. If we want to use Riemann's methods, then we can establish the following very general theorem." He then states a general uniformization theorem for polygons, and he says: ``The proof is exactly similar to the one given by Riemann in \S 12 of the first part of his \emph{Abelian functions} article for the particular polygon constituted by $p$ parallelograms stacked above each other and connected by $p$ ramification points. It seems to me that this theorem, which I shaped only these last days, included, as easy corollaries, all the existence proofs that you mention in your notes."} Such a map is known under the name \emph{Schwarz-Christoffel map}.\footnote{The name of Schwarz is attached the famous \emph{Schwarzian derivative} which later became one of the basic tools used in Poincar\'e's approach to uniformization via differential equations. Riemann used the Schwarzian derivative before Schwarz, namely, in his posthumous paper of 1867 on minimal surfaces. The idea of the Schwarzian derivative is also contained in Kummer's paper on the hypergeometric series. (We learned this from Jeremy Gray.) According to Klein (\cite{Klein-Lectures} p. 35), the name ``Schwarzian derivative" was coined by Cayley. The Schwarzian derivative plays a major role in the modern theory of Teichm\"uller spaces and moduli. Let us also recall that the Schwarzian derivative is a basic tool in the work of Bers. It appears in particular in the \emph{Bers embedding} of Teichm\"uller space, a holomorphic embedding of that space into a complex vector space of the same dimension, namely, the space of integrable holomorphic quadratic differentials on a Riemann surface.} 
  The work on the Riemann mapping theorem culminated in the proof in 1907 by Koebe and Poincar\'e of a very general uniformization theorem. We shall report on it in more detail in \S \ref{s:PK} below.
  Comparing Riemann and Weierstrass, Poincar\'e,  writes, in  a survey on the work of the latter \cite{Poincare-Weierstrass} (p. 7):
  \begin{quote}\small
  In one word, the method of Riemann is primarily a method of discovery; that of Weierstrass is primarily a method of proof.  
  \end{quote}
   Klein obviously stood up for Riemann. In the introduction of his book on the presentation of Riemann's ideas \cite{Klein-Riemann}, he writes:
\begin{quote}\small
Riemann, as we know, used Dirichlet's principle in their place in his writings. But I have no doubt that he started from precisely those physical problems, and then, in order to give what was physically evident the support of mathematical reasoning, he afterwards substituted Dirichlet's principle. Anyone who clearly understands the conditions under which Riemann worked in G\"ottingen, anyone who had followed Riemann's speculations as they have come down to us, partly in fragments, will, I think, share my opinion. -- However that may be, the physical method seemed the true one for my purpose. For it is well known that Dirichlet's principle is not sufficient for the actual foundation of the theorems to be established; moreover, the heuristic element, which to me was all-important, is brought out far more prominently by the physical method. Hence the constant introduction of intuitive considerations, where a proof by analysis would not have been difficult and might have been simpler, hence also the repeated illustration of general results by examples and figures.
\end{quote}
        \section{Klein}\label{s:Klein}
 
 Klein was  passionately devoted to mathematical research. He was among the first to acknowledge the works of Lobachevsky on hyperbolic geometry. He gave a formula for the disc model of the hyperbolic metric (which was already discovered by Beltrami) using the cross ratio, making connections with works of Cayley and giving the first interpretation of hyperbolic geometry in the setting of projective geometry; cf. his articles \cite{Klein-Ueber} and \cite{Klein-Ueber1}, and the commentary \cite{ACPK1}. Hyperbolic geometry eventually played a major role in the theory of Riemann surfaces, and this is an aspect Riemann did not touch. 
  Klein is also well known for the formulation of the \emph{Erlangen program} \cite{Klein-Erlangen}, in which geometry is considered as transformation groups rather than  spaces, and the study of invariants of these groups. The text of this program was set out at the occasion of Klein's inaugural lecture when he became professor at the University of Erlangen. The program served as a guide for mathematical and theoretical physical research, during several decades. We refer the reader to the papers in the book \cite{JP} for a modern view on the Erlangen program.

After the collapse of his health in 1882, Klein's research slowed down but he remained a great organizer and editor,\footnote{Klein was, for several years, the leading editor of the \emph{Mathematische Annalen}. He also co-edited Gauss' \emph{Collected Works} and the famous \emph{Encyklop\"adie der mathematischen Wissenschaften mit Einschluss iher Anwendungen} (Encyclopedia of Pure and Applied Mathematics). He was awarded the official title of ``Geheimrat" (Privy Councilor) by the German state. The history of this title goes back to the Holy Roman Empire, and it was awarded to very few people. Among them were Leibniz, Goethe, Gauss, Planck, Fricke  and Hilbert. Klein also  became the representative of the University of G\"ottingen at the upper chamber of the Prussian parliament.} and his influence on  the German mathematical school of the last quarter of the ninteenth century is enormous. At the turn of the century, he was considered to be the leader of German mathematics. Under his guidance, the University of G\"ottingen became the German center for mathematics and theoretical physics. The main factor was that Klein attracted there some of the most talented mathematicians, including Hilbert, Minkowski, Koebe, Bieberbach, Courant and Weyl.   Of major importance is also Klein's relationship with Poincar\'e, part of which survives in the form of a rich correspondence.
  
   Hermann Weyl credits Klein with the concept of an abstract Riemann surface. He writes, in the Preface of \cite{Weyl1913} (p. {\sc vii} of the 1955 edition):
\begin{quote}\small
    Klein had been the first to develop the freer concept of a Riemann surface, in which the surface is no longer a covering of the complex plane; thereby he endowed Riemann's basic ideas with their full power. It was my fortune to discuss this thoroughly with Klein in diverse  conversations. I shared his conviction that Riemann surfaces are not merely a device for visualizing the many-valuedness of analytic functions, but rather an indispensable essential component of the theory; not a supplement, more or less artificially distilled from the functions, but their native land, the only soil in which the functions grow and thrive. 
    \end{quote}

 Between the years 1880 and 1886, Klein lectured at the University of Leipzig on Riemann's work on algebraic functions and on Riemann surfaces. His book \cite{Klein-Riemann}, in which he explains the ideas of Riemann and titled \emph{On Riemann's theory of algebraic functions and their integrals; a supplement to the usual treatises}, is connected with these lectures. It contains a section on  Riemann's moduli space. Klein states the following result and he explains how to obtain it (p. 82):

\begin{theorem}
All algebraical equations with a given $p$ form a single continuous manifoldness in which all equations derivable from one another by a uniform transformation constitute an individual element.
\end{theorem}
 
Klein writes that ``the totality of the $m$-sheeted surfaces with $w$ branch-points form a \emph{continuum}," and for this property he refers to theorems of L\"uroth and Clebsch in \emph{Mathematische Annalen} t. {\sc iv, v}.\footnote{The papers of L\"uroth and Clebsch on the subject include \cite{Luroth1872} \cite{Clebsch1872a} \cite{Clebsch1872b} \cite{Clebsch1872c}. The word ``continuum"  appears very often in the writings of Klein and his contemporaries. In this context, a continuum is an arcwise connected open subset of the complex plane. 
In 1851 \cite[p. 129]{bo},  Bolzano had a  different definition, to which 
 Cantor made objections. We refer to the paper by \cite[p. 226]{mo} for several historical notes on this subject. Bolzano's definition was:
``a continuum is present when, and only when, we have an aggregate of simple entities (instants or points or substances)
 so arranged that each individual member of the aggregate has, at each individual and sufficiently small distance from 
 itself, at least one other member of the aggregate for a neighbor." The notion adopted by Klein was introduced by Weierstrass in his Berlin lectures, cf. Weierstrass, \emph{Monatsb. Akad. Berlin} 1880 p. 719-723 and Weierstrass' \emph{Werke} 2, Berlin 1895, p. 201-205. This was repeated by Weierstrass  in 1886 \cite[p. 65]{w2}. Weierstrass' definition was also adopted by G. Mittag-Leffler  in \cite[p. 2]{mit}. This definition and the more general one are recalled in the article by Osgood in Klein's Encyclopedia \cite{Osgood-Encyclopaedia} p. 9. In his 1928 edition of  textbook on complex analysis \cite[p. 162]{os2}, Osgood called a set in the complex plane a 
``two-dimensional continuum" if it is path-connected and open. In the book by Picard and Simart \cite{Picard-Simart}, (Vol. I, p. 24, in the chapter on \emph{Analysis situs}) the authors write: ``Instead of using the expression \emph{space}, we shall often use the words \emph{variety} (vari\'et\'e), \emph{multiplicity}, and \emph{continuum}. 
In 1910 \cite[p. 138]{Denjoy1910}, Denjoy writes: ``Let us call a continuum (Gebiet) a set which only contains interior points and which is connected, and let us call a domain (Bereich) the sum of a continuum and its boundary."}
 
Besides the word continuum, the word ``dimension" appears  for the first time in the description of the quantity $3p-3$. Klein writes in \cite{Klein-Riemann} p. 82: 
\begin{quote}\small
Hence we conclude that \emph{all algebraic equations with a given $p$ form a single manifoldness}. [...] For the first time, a precise meaning attaches itself to the number of the moduli; \emph{it determines the dimensions of this continuous manifoldness}.\footnote{The italics are Klein's.}
\end{quote}

The set that Klein describes is the moduli space $\mathcal{M}_p$ of Riemann surfaces of genus $p$.  A few comments are in order.

1) The letter $p$ denotes, as before, the genus. The notion of \emph{genus} is associated here with an algebraic equation, an idea due to Clebsch, as is the word ``genus." This is the genus of the Riemann surface that is associated with the multi-valued function defined by the equation.

2) The term \emph{manifoldness} (``Mannigfaltigkeiten") was introduced by Riemann,\footnote{In the French edition of Riemann's works (1898), the term is  translated by ``multiplicit\'e" (multiplicity). This is also the term used in Hermite's Preface to that edition. The term was also used by Poincar\'e and other French geometers before the introduction of the word ``vari\'et\'e," which is today's French word for manifold and which had a more precise meaning.} originally as a mathematico-philosophical notion.\footnote{It seems that Riemann, who was first educated in theology, had in mind, a word which would generalize the word ``Dreifaltigkeit" which designates the holy Trinity. See the interesting article \cite{Plotnitsky}.} It is usually translated by ``manifoldness," and it appears in several places in the work of Klein. This word was never defined precisely, and the notion of a general ``manifold" with charts and coordinate changes was still inexistent. In Riemann's context,  the word manifoldness means (vaguely) an object which can be described using a certain number of complex parameters.\footnote{Notice that if the word ``Mannigfaltigkeiten" denoted a topological manifold (as it has been sometimes suggested), then Riemann and Klein would have written $6p-6$ and not $3p-3$.} 

3) The word ``single" refers to the fact that the space is connected.

Klein then adds some more precise facts (p. 85):
\begin{quote}\small
 To determine a point in a space of $3p-3$ dimensions we do not generally confine ourselves to $3p-3$ coordinates; more are employed connected by algebraical, or transcendental relations. But moreover it is occasionally convenient to introduce parameters, in which different series denote the same point of the manifoldness. The relations which then hold among the $3p-3$ moduli necessarily existing for $p>1$ have been but little investigated. On the other hand the theory of elliptic functions had given us an exact knowledge of the subject for the case $p=1$. 
\end{quote}

 Klein tried to make Riemann's count of moduli more precise.
 In the same booklet on Riemann's work, he recalls (\S 19) that ``Riemann speaks of all algebraic functions of $z$ belonging to the same class when by means of $z$ they can be rationally expressed in terms of one another," and then he writes: ``The number in question is the number of different classes of algebraic functions which, with respect to $z$, have the given branch-values." He then adds in a note: ``If I may be allowed to refer once more to my own writings, let me do so with respect to a passage in \emph{Mathematische Annalen} t. {\sc xii} (p. 173), which establishes the result that certain rational functions are fully determined by the number of their branchings, and again to ib. {\sc xv} p. 533, where a detailed discussion shows that there are ten rational functions of the eleventh degree with certain branch points." In fact, Klein, instead of  the quantity $3p-3$, has $3p-3+\rho$, where $\rho$ is ``number of degrees of freedom in any one-to-one transformation of a surface into itself," that is, the dimension of the isometry group of the surface. (This takes care of the case of genera 0 and 1 in the formula.) The form $3p-3+\rho$ is the one which appears later in the work of Teichm\"uller. In his count of the number of moduli, Klein gives an argument which is basically the same as Riemann's: Each Riemann surface can be realized as an $m$-covering over the Riemann sphere, with $m > 2p-2$. There are $w=2m+2p-2$ branch points, and for each Riemann surface of genus $p$ there are $2m-p+1$ ways to cover   the sphere with $m$ sheets. Once the branches are fixed, there are finitely many ways to glue the various sheets to get a Riemann surface.\footnote{We know that the finite number of choices cause a singularity of moduli space (or, its orbifold structure).}
 From this number of parameters, one needs to subtract the $2m-p+1$ ways of realizing
a Riemann surface of genus p as an $m$-sheet cover. Thus, the total count for moduli is
$w-(2m-p+1)=2m+2p-2-(2m-p+1)=3p-3$. Klein then writes on p. 82: ``Thus for the first time, a precise meaning attaches itself to the number of moduli: it determines the dimensions of this continuous manifoldness."\footnote{The reader may notice that several problems arise if one wants to make these statements precise. Namely, taking quotients might destroy any manifold property of the space. Singularities might also appear from the collision of branch points. Thus, it isn not always true that we get a manifold structure on the resulting quotient space. This is the sort of problem that led Teichm\"uller later to say that people made counts without knowing what they were talking about (see Footnote \ref{F:Teich} and the discussion in \S \ref{s:Teich} below).} Klein could not define a manifold structure on the moduli space $\mathcal{M}_p$ (and in fact, this space is not a manifold), and the moduli count does not provide that space with a topology. It is therefore not a surprise that his ``method of continuity," which he used in his attempt to prove uniformization in 1882, and on which we shall comment later, was considered as problematic. One should also add that Klein's objective, in going through Riemann's ideas, was to transmit them, and not to make them more rigorous, since he never doubted their validity.\footnote{Other mathematicians were suspicious. In his essay on \emph{Riemann and his importance for the development of modern mathematics} (\cite{Klein-Wien} p. 70), Klein writes:  ``Riemann's methods were kind of a secret for his students and were regarded almost with distrust by other mathematicians." (Remmert's  translation \cite{Remmert1998} p. 206).} Like Riemann, Klein was more interested in general ideas than in details.\footnote{\label{f:Brunel} We mention here a letter to Poincar\'e, dated July 7, 1881, by Georges Brunel, a young French mathematician following Klein's lectures in Leipzig (\cite{HP-correspondance1} p. 94): ``I am extracting from his lectures two or three lines which show clearly the fundamental idea and which answer your question about Riemann: `Riemann stated the theorem that on any Riemann surface there exist functions, but this theorem cannot be considered as proved.' This is written by the hand of the student who did the writing. And then, by the hand of Klein: \emph{Despite this fact, I will use without hesitation this theorem. I think it is possible to give a rigorous proof of the general proposition.}"}  Like his famous predecessor, Klein also resorted to physics. Right at the beginning of his expository monograph \cite{Klein-Riemann}, he  writes (p. 1-2): ``The physical interpretation of those functions of $x+iy$ which are dealt with in the following pages is well known. [...] For the purpose of this interpretation it is of course indifferent of what nature we may imagine the fluid to be, but for many reasons it will be convenient to identify it here with the \emph{electric fluid}." Klein refers to Maxwell's \emph{Treatise on electricity and magnetism} (1873).\footnote{In the introduction of his book \cite{Klein-Riemann}, Klein writes: ``I have not hesitated to take these physical conceptions as the starting point for my presentation." Let us note however that in the paper \cite{Bottazzini}, Bottazzini argues that Prym and others denied that Riemann had been pursuing a physical analogy.}  It is interesting to note that in his essay, Klein draws equipotential curves which are analogous to the measured foliations that are used in the Thurston theory of surfaces. The figure below is extracted from Klein's essay \cite{Klein-Riemann}.
\[
\includegraphics[width=1\linewidth]{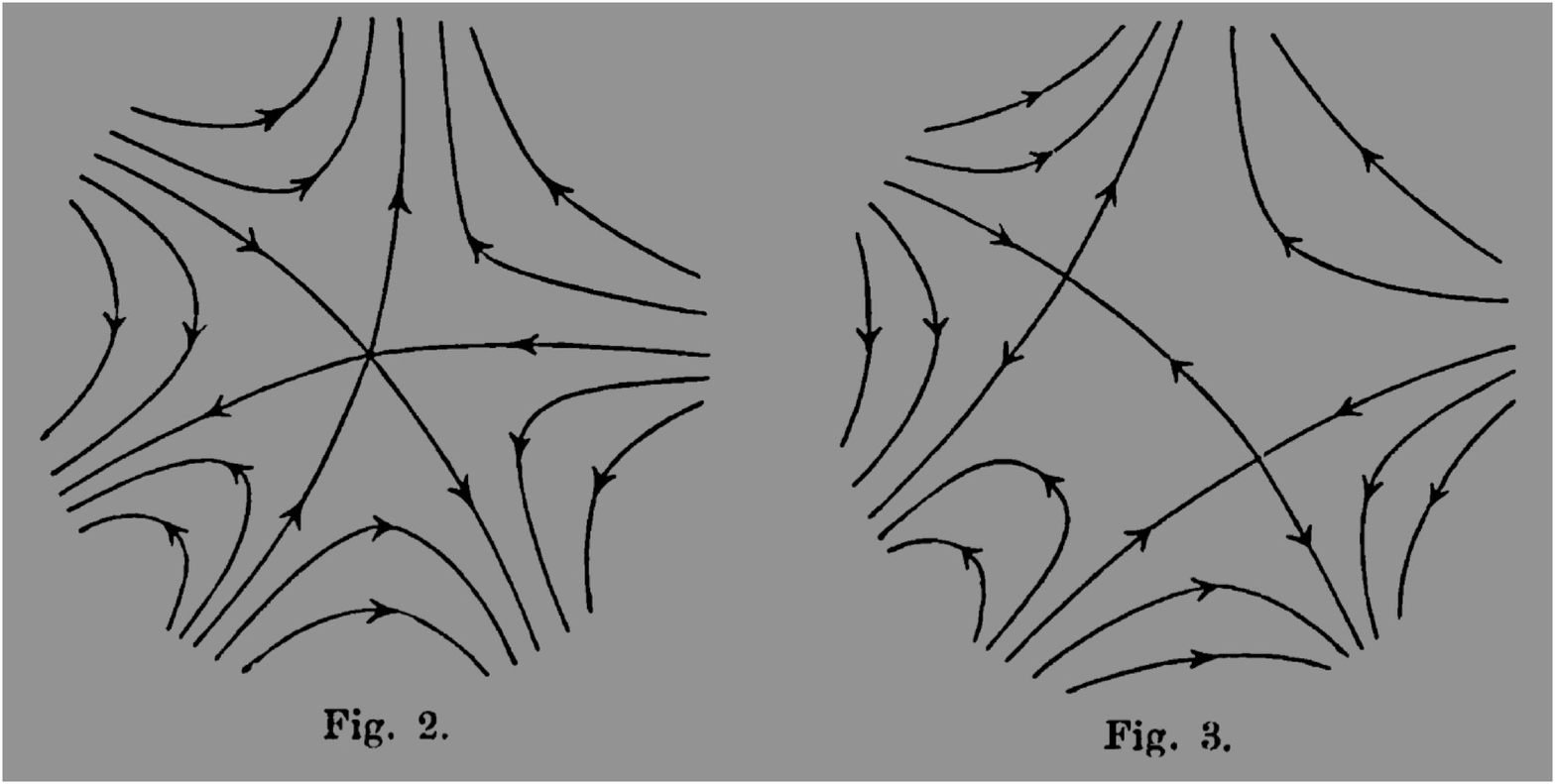}
\]

 Let us now quote Poincar\'e, from his \emph{La valeur de la science} (1905):\footnote{Our translation from the French.} 
   \begin{quote}\small
   Look at the example of Mr. Klein: He is studying one of the most abstract questions in the theory of functions; namely, to know whether on a given Riemann surface there always exists a function with given singularities. What does the famous German geometer do? He replaces his Riemann surface by a metal surface whose electric conductivity varies according to a certain rule. He puts two of its points in contact with the two poles of a battery. The electric current must necessarily pass, and the way this current is distributed on the surface defines a function whose singularities are the ones prescribed by the statement.  
      \end{quote}

 He several important discoveries in function theory, although he himself considered that he was merely explaining the ideas of Riemann. Indeed, as we mentioned several times, by his writings and his achievements, Klein was a true perpetuator of the tradition of Riemann. 
 In his biography of Arnold Sommerfeld, who was Klein's assistant in G\"ottingen,  M. Eckert writes (\cite{Eckert} p. 50): ``Riemann was the great model for Klein. For Riemann, the proximity of mathematics to physics had been axiomatic. Klein advised his students to acquire a personal sense of Riemann's work from the primary source literature." Klein also contributed to invariant theory, number theory,  algebra and differential equations. For a concise biography of Klein, the reader may consult \cite{Lizhen-Klein}.   

 Besides his work on discontinuous group actions on the hyperbolic plane, Klein made several contributions to Riemann surfaces. He  addressed the problem of whether every compact Riemann is equivalent to a surface smoothly embedded in a Euclidean space, on which the holomorphic structure is defined by the angle structure induced from that of the ambient space.\footnote{The problem was also addressed by Gauss. It was solved by Garsia \cite{Garsia} for closed surfaces, and by R\"uedy \cite{Ruedy} in the general case.} See \cite{Klein-Gesammelte} Vol. 3 p. 635. He  knew that the automorphism group of a compact Riemann surface of genus $\geq 2$ is finite. He states in \cite{Klein-Riemann}: ``It is always possible to transform into themselves in an infinite number of ways by a representation of the first kind surfaces for which $p=0$, $p=1$, but never surfaces for which $p>1$."

 \section{Poincar\'e}\label{s:PK}

 As a mathematician, Poincar\'e (1854-1912) was a geometer, an analyst, an algebraist, a number theorist, a dynamical system theorist, and above all a topologist. His contribution to the theory of Riemann surfaces involves  all these fields. 
 In fact, Poincar\'e considered topology as part of geometry. In his \emph{Derni\`eres pens\'ees} (Last thoughts), published posthumously (1913) \cite{PK-dernieres}, he writes: 
 \begin{quote}\small
 Geometers ordinarily classify geometry into two types, calling the first one metric and the second one projective. [...] But there is a third geometry, where quantity is completely banished and which is purely qualitative: this is \emph{Analysis Situs}. [...]  \emph{Analysis Situs} is a very important science for the geometer. It gives rise to a series of theorems which are equally  linked as the ones of Euclid. It is on this set of propositions that Riemann constructed the most remarkable and abstract theories of pure analysis. [...]  This is what makes this \emph{Analysis Situs} interesting: the fact that it is there that geometric intuition is really involved.
  \end{quote}
  
  P. S. Alexandrov,\footnote{Pavel Sergue\"ievitch Aleksandrov (1896-1982) studied in Moscow under Egorov and Luzin (the latter had studied in G\"ottingen between 1910 to 1914). Aleksandrov made major contributions to topology, where his name is attached to the Alexandrov compactification and the Alexandrov topology. Alexandrov visited G\"ottingen several times and taught there.  In her biography of Courant, Constance Reid writes (\cite{Reid} p. 106): ``Since 1923 Alexandrov had returned each year, either alone or accompanied by countrymen. From 1926 through 1930, Courant always arranged for him to give courses in topology, each for a quite different audience of mathematicians. The summer that Courant was trying to keep up some semblance of attendance at Wiener's lectures, Alexandrov's were crowded."} in a talk he gave at a celebration of the centenary of Poincar\'e's birth \cite{Alex}, says the following: ``To the question of what is Poincar\'e's relationship to topology, one can reply in a single sentence: he created it [...]." Although this may be, strictly speaking, an overstatement, there is a lot of support for it.\footnote{In his \emph{Analysis of his own works} (\cite{Poin-Acta} p. 100), Poincar\'e considers that his predecessors, in topology, are Riemann and Betti. The latter was a friend of Riemann and was influenced by him. The book \cite{Pont1974} and the article \cite{Weil1979a} contain two letters from Betti to his colleague and friend Tardy in which he reports on conversations he had with Riemann on topology.} We owe to Poincar\'e such basic notions as the fundamental group, homology, Betti numbers, torsion coefficients and duality.  He is the founder of combinatorial topology and of the qualitative approach to differential equations from the point of view of the properties of the vector fields they generate: singularities, existence of periodic orbits, behavior of integral curves, etc. This is also part of topology. Poincar\'e's study of discrete group actions is also a chapter in topology. He is the author of the celebrated \emph{Poincar\'e conjecture} stating that simply connected, closed 3-manifold is homeomorphic to the 3-sphere and which was, until its recent proof, one of the most central open problems in mathematics. For a history of algebraic topology, and in particular on Poincar\'e's contribution, we refer the reader to the books \cite{Pont1974} by Pont and \cite{Scholz} by Scholz.
Poincar\'e was also a physicist. His works on special relativity and on celestial mechanics, which includes the three body problem and the question of the stability of the solar system and where his qualitative study of differential equations finds other beautiful applications, are decisive in the fields. It is not possible to describe in a few pages his achievements. We are interested in his work on Riemann surfaces and in his relationship to Klein. We refer to the recent scientific biography by Gray \cite{Gray-PK}.

In the 1880s, Poincar\'e gradually became the leading mathematician in France, as Klein was in Germany.  Hadamard, in his commentary on Poin\-car\'e's mathematical work \cite{Hadamard-PH}, writes that unlike other scientists,  Poincar\'e, in the choice of the subjects of his investigations, did not follow his personal inclinations, but the need of science (p. 214). He adds: ``He was present everywhere a serious gap had to be filled or a big obstacle had to be overcome." Poincar\'e's work on Riemann surfaces involves, like Riemann's and Klein's, function theory, topology, and potential theory. He also introduced in that theory techniques of group theory and the newly discovered methods of hyperbolic geometry. Unlike Klein, Poincar\'e did not consider a Riemann surface as an object obtained by making cuts in the complex plane and gluing them again, but from the beginning, he thought of it in terms of polygons in the plane, whose edges are circular arcs, with discontinuous groups acting on them, the Riemann surface appearing as a quotient of this action. He later identified the disk with the hyperbolic plane, and the discontinuous groups with groups of hyperbolic motions.  He called these groups Fuchsian, a name he used after he read a paper by Lazarus Fuchs on second-order differential equations \cite{Fuchs1866}. He also introduced automorphic functions as a major ingredient in the theory of Riemann surfaces. 
 On Poincar\'e's style and his relation to Riemann, Alexandrov writes in \cite{Alex}: 
 \begin{quote}\small
 The close connection of the theory of functions of a complex variable, which Riemann has observed in embryonic form, was first understood in all its depth by Poincar\'e [...] Poincar\'e was a master of topological intuition as no other mathematician of his time and preceding eras; may be only Riemann could be compared to him in this respect, but he did not succeed in developing it with such breadth and diversity of applications as Poincar\'e. 
 Topological intuition penetrates the majority of Poincar\'e's most signifiant works -- the theory of automorphic functions and uniformization (this is the supreme triumph of the `Riemann' approach to the theory of functions of a complex variable).
 
 The force of Poincar\'e's geometrical intuition sometimes led him to ignore the pedantic strictness of proofs. Here there is still another side; finding himself under constant influx of a set of ideas in the most diverse fields of mathematics, Poincar\'e `did not have time to be rigorous', he was often satisfied when his intuition gave him the confidence that the proof of such and such a theorem could be carried through to complete logical rigor and then assigned the completion of the proof to others. Among the `others' were mathematicians of the highest rank.
 \end{quote}

 Poincar\'e came to Riemann surfaces and uniformization because he was interested in the global behavior of solutions of differential equations, more precisely, of second order linear equations whose coefficients are meromorphic functions.\footnote{The subject of Poincar\'e's doctoral thesis (defended in 1879) was differential equations. This subject prepared him for his subsequent research in mathematics and physics. In 1880, Poincar\'e obtained the ``mention Tr\`es Honorable" for the annual competition organized by the Acad\'emie des sciences. The subject was ``To improve in some important way the theory of linear differential equations with one independent variable." Poincar\'e submitted a memoir titled \emph{The integration of all linear differential equations with algebraic coefficients}. Hadamard writes in his survey of Poincar\'e's work (\cite{Hadamard-PH} p. 206): ``The integration of differential equations and of partial differential equations remains until now the central problem of modern mathematics. It will presumably stay one of its major problems, even if physics continues the path it is following presently towards the discontinuous."} His aim was to prove that one can express these solutions by a uniform function.  Thus, like Riemann, but through another path, Poincar\'e was led to the question of uniformization. It seems that Poincar\'e, when he started working on the subject, was not aware of Riemann's work, and that first learned about it through his correspondence with Klein.\footnote{Dieudonn\'e writes in his biography of Poincar\'e \cite{Dieu-Poin}: ``Poincar\'e's ignorance of the mathematical literature, when he started his research, is almost unbelievable. He hardly knew anything on the subject beyond Hermite's work on modular functions; he certainly had never read Riemann, and by his own account had not even heard of the Dirichlet principle."}

 Several of Poincar\'e's early results on differential equations are announced in a series of \emph{Comptes Rendus} notes written  in 1881 and 1882, and more details are published in the 1882 \emph{Acta Mathematica} papers of the same year \cite{P1882a} \cite{P1882b} and in later ones.  The 1881 \emph{Comptes Rendus} note \cite{P1881} starts as follows: 
 \begin{quote}\small
 The goal which I propose in this work is to search if there are analytic functions which are analogous to the elliptic functions and which allow the integration of various linear differential equations with algebraic coefficients. I managed to prove that there exists a very wide class of functions which satisfy these conditions and to which I gave the name \emph{Fuchsian functions}, in honor of Mr. Fuchs, whose works have been very useful to me in this research.
  \end{quote}

  Poincar\'e started a correspondence with Fuchs in 1880, after he noticed the latter's 1866 paper on linear differential equations with complex coefficients. Several letters from Poincar\'e to Fuchs are published in \cite{PK-Fuchs1}, and two letters from Fuchs, in French translation, are published in \cite{HP-correspondance1}. In these letters, Fuchs responds to Poincar\'e who had asked for clarifications and details on his works.  \[
\includegraphics[width=1\linewidth]{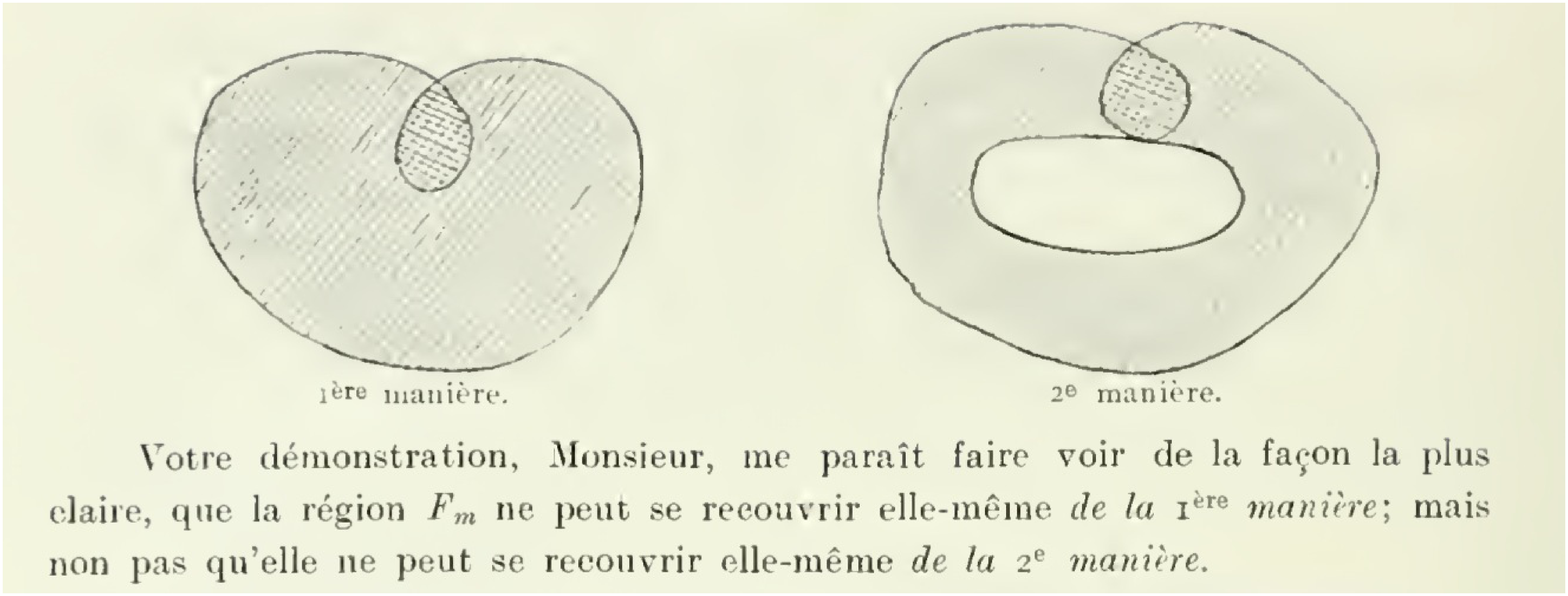}
\]
 The above figure is extracted from a letter from Poincar\'e to Fuchs, where the former constructs a region in the plane obtained by analytical extension of a meromorphic function which covers itself. He explains that in a proof, Fuchs missed one of the two cases. In fact, Poincar\'e saw that there were imprecisions, inconsistencies and other problems in the writings of Fuchs, but he nevertheless remained faithful to the terminology ``Fuchsian equation," ``Fuchsian function" and ``Fuchsian group," in gratefulness to the ideas that he got by reading Fuchs' papers. He eventually gave new proofs and generalized the results of Fuchs. Dieudonn\'e (whose opinion is sometimes extreme), in a letter he wrote to the editors of \cite{HP-correspondance1}, says that ``it seems evident that Fuchs never understood the objections of Poincar\'e" (quoted on p. 152).

 Poincar\'e's theory of Fuchsian groups and Fuchsian functions provided a lot of examples of Riemann surfaces which are equipped by non-constant meromorphic functions, which is in the direct lineage of Riemann's program. As a by-product of this work, Poincar\'e obtained a version of the uniformization theorem, which is also one of the important questions addressed by Riemann. Poincar\'e stated this result in his 1883 article \cite{P1883}. We reproduce the statement, which is amazing in its generality:
 \begin{theorem} \label{unifomr:PK}
 Let $y$ be an arbitrary non-uniform analytic function of $x$. Then, we can always find a variable $z$ such that $x$ and $y$ are uniform functions of $z$.
 \end{theorem}
In this theorem, the function $y$ is not necessarily  algebraic. 
 
 At about the same time, Klein announced a proof of a similar theorem, in his note \cite{Klein-MA1882}. Klein's statement is closer to the spirit of Riemann. It is summarized in the following:
  \begin{theorem} \label{unifomr:Klein}
 For every Riemann surface, there exists a single-valued function $\eta$, which is uniquely defined up to composition by a fractional linear transformation, which, after cutting the surface into a simply-connected region, maps it conformally onto a simply-connected region of the sphere, and such that along the cuts the function changes by fractional linear transformations.
 \end{theorem}

  In general, the works of Klein and Poincar\'e were close, although their approaches to several problems were different. Poincar\'e kept thinking in terms of solutions of differential equations while Klein was following the function-theoretic path traced by Riemann. 
Klein's vision of geometry as transformation groups expressed in the \emph{Erlangen program} (1871) was completely shared by Poincar\'e. Klein was also among the first to recognize the importance of Lobachevsky's geometry and to transmit it through his writings; this geometry played later on a major role in Poincar\'e's writings.

The uniformization theorem is discussed in  the correspondence between the two men.\footnote{This correspondence was edited by various people. One set of letters was published by N. E. N\"orlund in  \emph{Acta Mathematica} \cite{Noerlund} (1923). A second set is included in Klein's \emph{Collected Works} \cite{Klein-Gesammelte}. A third set is contained in the collection of mathematical letters  of Poincar\'e, published with French translations, in two parts, \cite{HP-correspondance1} and \cite{HP-correspondance2}. Some of the letters are reproduced in the book \cite{uni}. 
The second volume of the correspondence \cite{HP-correspondance2} contains 26 letters exchanged with Klein between June 1881 and September 1882, 3 letters exchanged in 1895, and  one letter for each of the years 1901, 1902, 1906.} Most of these letters concern Riemann surfaces and they also give us an idea of the period when the two men were in close relationship with one another.

When this correspondence started, in 1881, Poincar\'e was 27 years old, and Klein was 32.
Klein wrote in German and Poincar\'e responded in French.  Poincar\'e, like Klein, was investigating discrete groups of linear fractional transformations and their automorphic functions, and he discovered the relation between these groups and the group of  isometries of the hyperbolic plane. 
 Both had developed methods  for constructing Riemann surfaces by assembling polygons in the plane, and they both constructed associated discontinuous groups. Hadamard writes in \cite{Hadamard-PH} that Poincar\'e's discovery of Fuchsian groups made in 1880 ``attracted the attention and the admiration of all geometers." Hadamard also recalls the result on Fuchsian functions came as a surprise to Poincar\'e, whose initial aim was to show that these functions, as generalizations of the modular function and the inverse of the hypergeometric series, do not exist (\cite{Hadamard-PH} p. 208).\footnote{This is stated by Poincar\'e in \emph{Science et M\'ethode} \cite{PK-Science} p. 51. The passage is reproduced below.} The two years 1881 and 1882 were a period of intense mathematical research for Klein.  After that, his research activity slowed down considerably because of a nervous breakdown, and his correspondence with Poincar\'e slowed to a trickle.

The first letter is from Klein, dated June 12, 1881.  It starts by: 
``Sir, Your three \emph{Comptes Rendus} Notes \emph{Sur les fonctions fuchsiennes}, of which I became aware only yesterday, and only rapidly, are so closely related to the thoughts and the efforts that occupied me during the last years that I feel obliged to write to you."
Then Klein reports on several of his articles, and also on some unpublished results concerning triangles, and (more generally) polygons, whose boundaries are pieces of  ``circular arcs."  The next two figures, representing such polygons, are extracted from a letter from Klein to Poincar\'e \cite{HP-correspondance2}.
 \[
\includegraphics[width=.4\linewidth]{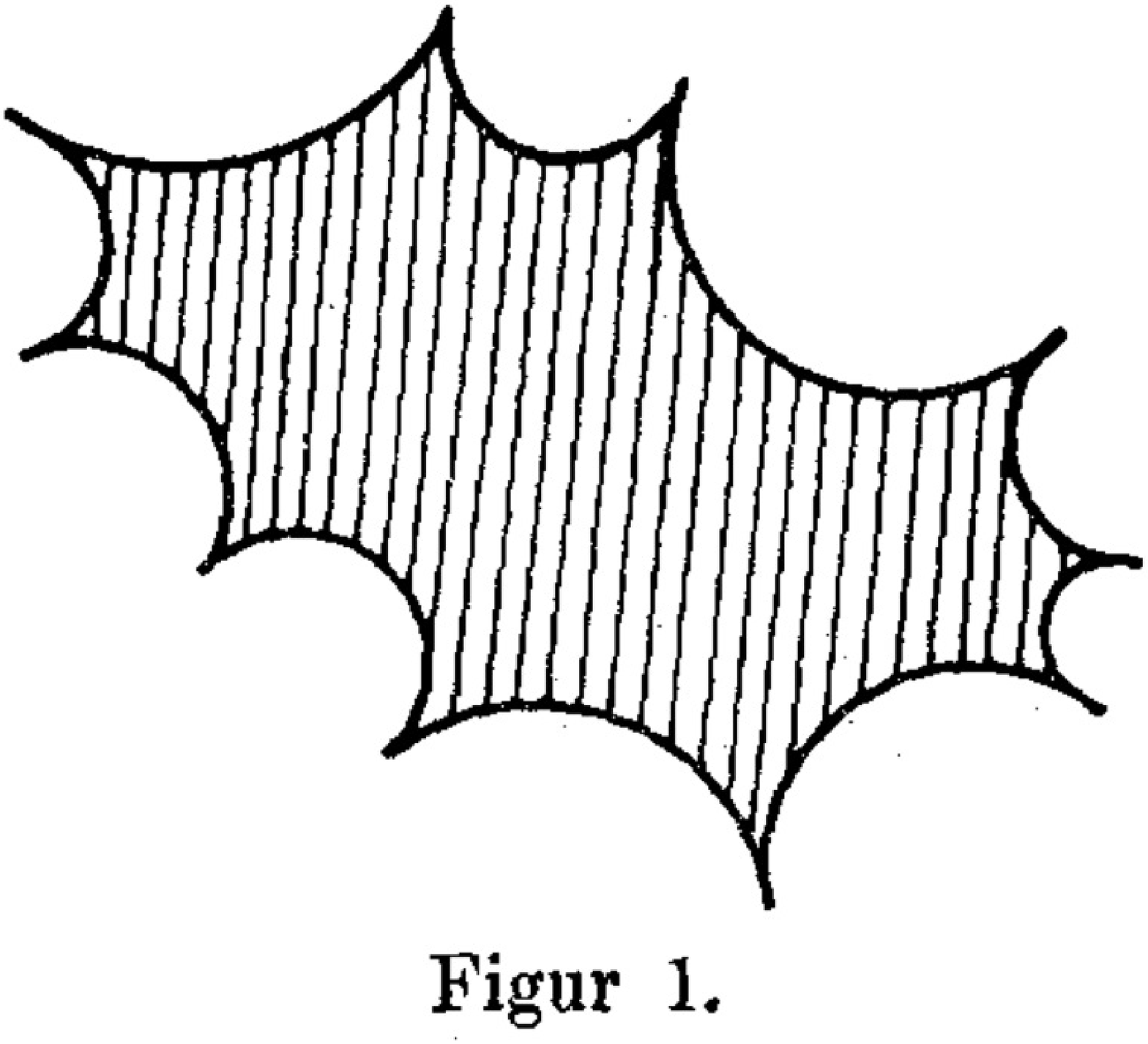}
 \includegraphics[width=.6\linewidth]{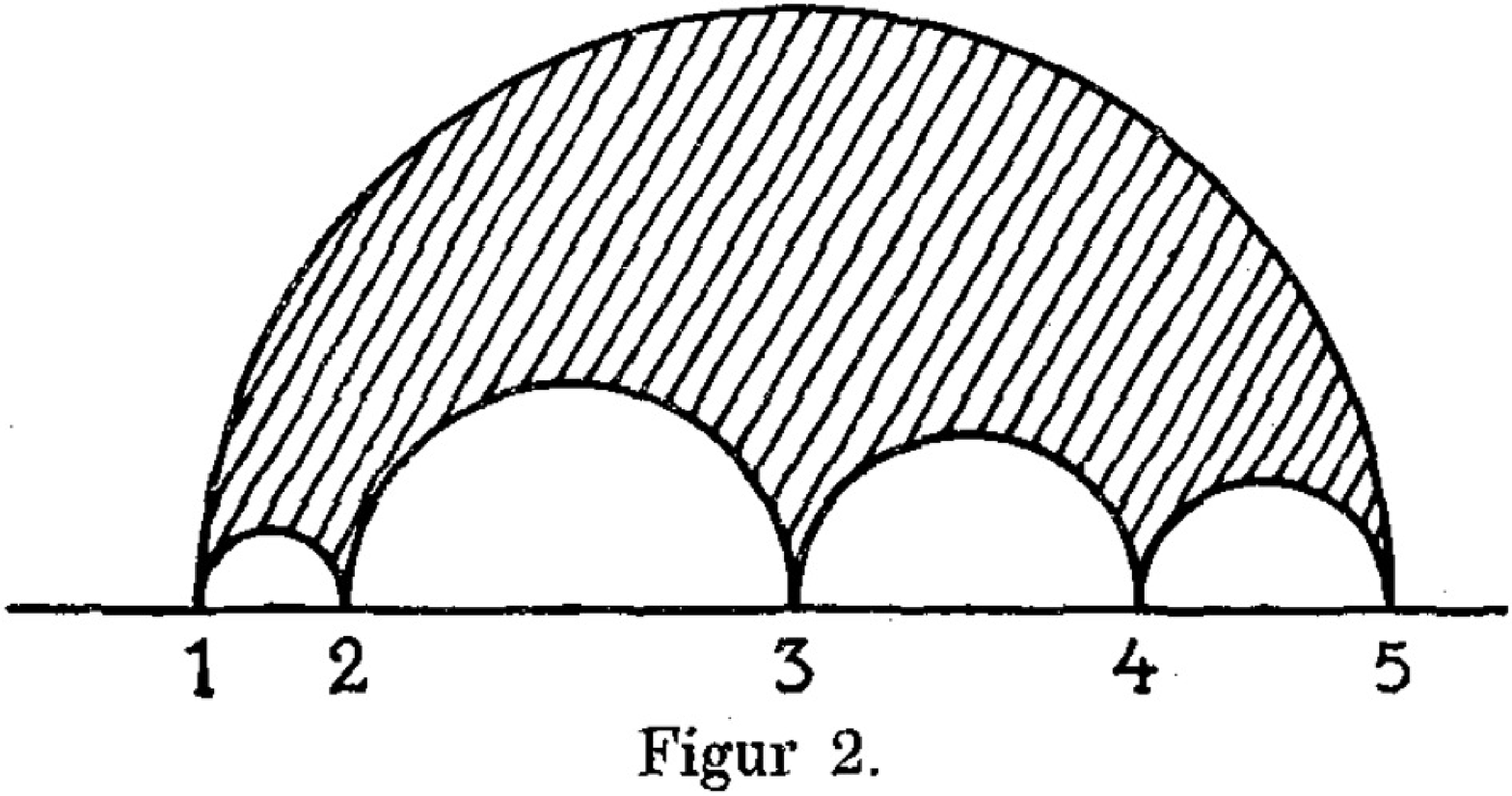}
\]
Klein writes: ``I never published anything about these concepts, but I presented them, in the summer of 1879, in a course at the General Technical School of Munich." He concludes his letter by the sentences: ``Your approach, the one which guides your works since 2-3 years, is in fact very close to mine. I will be very pleased if this first letter gives rise to a regular correspondence. It is true that at this moment other commitments take me away from these works, but I am all the more excited of  resuming them that I have to give, next winter, a course on differential equations."

 Poincar\'e responds on June 15, 1881:  ``Your letter proves that you saw before me some of my results on Fuchsian functions. I am not at all surprised, because I know how much you are knowlegeable in non-Euclidean geometry, which is the true  key to the problem with which we are concerned." Poincar\'e promises to quote Klein's work in his future publications. He then asks him several questions: ``What is the \emph{Theory of the fundamental polygon} (Fundamentalpolygone)? Did you manage to find all the \emph{polygons with circular arcs}  (Kreisbogenpolygone) which give rise to a discontinuous group? Did you prove the existence of functions that correspond to every discontinuous group?" To these questions, Klein responds in his next letter (June 19, 1881): ``I was not able yet to establish the existence of all the discontinuous groups. I only noticed that there were several for which there is no fixed fundamental circle\footnote{The \emph{fundamental circle} is the circle invariant by a Fuchsian group. The terminology seems to be due to Poincar\'e. This remark made by Klein shows that he knew the existence of quasi-Fuchsian groups.} and therefore to which we cannot apply an analogy with non-Euclidean geometry (with which, by the way, I am very familiar)." We do not wish to enter here into the mathematical details contained in these letters, and the interested reader should go through them. The results discovered and the problems raised there are at the basis of the theory of Fuchsian groups and its relation to Teichm\"uller theory.

In several letters, the choice of the name \emph{Fuchsian function} is discussed. Klein considers that this is a misnomer and that Fuchs' results were known to others before him. Klein seems to have good reasons to complain, but Poincar\'e, who had already decided about this name and who had used it in his published articles, is unrelenting.\footnote{There are several places in the mathematical literature where one can read on Klein's discontent about Poincar\'e's choice of the adjective \emph{Fuchsian}. We mention a  source which is rarely quoted, viz., a letter from Brunel to Poincar\'e, written from Leipzig, on June 1881 (\cite{HP-correspondance1} p. 92). Brunel reports that at on that day, Klein, at his seminar, presenting works of Poincar\'e,  declared to the audience: ``I protest about the choice of the name \emph{Fuchsian function}. The fundamental idea is due to Riemann, and the merit of applying this idea of Riemann goes to Schwarz. Later, I myself worked in that direction and in my lectures at the Munich \emph{Polytechnikum} I presented a few results which were at the basis of the work of Mr. Poincar\'e. As to Mr. Fuchs, who once wanted to deal with similar questions, he only managed to do the following: to show us that in this matter he understands strictly nothing." Brunel adds: ``I may conclude: your Fuchsian functions belong to you and to nobody else, you can give them the name you want, and no German can criticize this." In his letter dated June 27, 1881, Poincar\'e acknowledges that he might as well have chosen another name for these functions, namely, Schwarzian functions, but he says that out of respect for Fuchs, he cannot change the name. Today, the terminology ``Fuchsian function" is no more used; these functions are called ``automorphic," a word that was suggested by Klein in \cite{Klein1890} p. 549. But the terminology ``Fuchsian group" survives.} It is also worth noting that in a letter dated January 13, 1882, Klein protests Poincar\'e's choice of the word \emph{Kleinian}. He considers that the credit is rather deserved by Schottky, and that in any case, the origin of everything is in Riemann's work. 

In his letter dated June 27, 1881, Poincar\'e asks Klein for details about the classification  of groups according to \emph{genera}. In his letter dated July 2, 1881,  Klein informs Poincar\'e that some of his work on  theta functions is not new, and that a large number of young mathematicians are working on the subject,  viz., to find conditions that distinguish Riemann's theta functions from the more general ones. He notes that he is surprised by a remark of Poincar\'e saying that the constant associated with Riemann's theta functions is equal to $4p+2$, whereas it should be equal to $3p-3$, and he asks him: ``Didn't you read Riemann's explanations on that, or the discussion of Brill and Noether in  (\"Uber die algebraischen Functionen und ihre Anwendung in der Geometrie, \emph{Mathematische Annalen} 7 (1874), 269-310)?" On July 5, 1881, Poincar\'e writes that the $4p+2$ constants he mentioned are not the moduli but rather ``parameters," and that $4p+2$ is only an upper-bound for the number of moduli, and that this is sufficient for his needs. From here, we conclude that the term \emph{moduli} had the meaning of ``minimal" number of parameters.

In a letter dated  December 4, 1881, Klein asks Poincar\'e to send him an article (even in the form of a letter) which he proposes to publish in the \emph{Mathematische Annalen}. It is interesting to note that Klein asks Poincar\'e for an \emph{an outline} of his views and results. Klein was not interested in technical details. He also says that he could publish, in the same issue, a note which presents his point of view on the matter, and explain how Poincar\'e's actual program served as a principle for the orientation of his own work on modular functions. 

The letter from Klein to Poincar\'e dated January 13, 1882 concerns uniformization: 
\begin{quote}\small
I just wrote a small article \cite{Klein-MA1882}  which will be published next to yours. It presents, again without proof, some of the results in this domain, and above all the following:
\\
Any algebraic equation $f(w,z)=0$ can be solved in a unique manner by $w=\phi(\eta)$, $z=\psi(\eta)$ as long as we have drawn on the corresponding Riemann surface $p$ independent  cuspidal cuts [R\"uckkehrscnitte], where $\eta$ is a discontinuous group, such as those which you have told me about after my letter. This theorem is all the more beautiful that this group has exactly $3p-3$ essential parameters [wesentliche Parameter], that is, as many parameters as the equations of the given $p$ possess moduli.

[..] But the proof is difficult. I always work with the ideas of Riemann, respectively with the \emph{geometria situs}. It is very difficult to write this clearly. I will try to do it later on. In the meanwhile, I would like very much to correspond with you on this matter, and also on your proofs. Be sure that I will study with the greatest interest your letters that you allow me to hope on this subject and that I will respond rapidly.  If you want to publish them, under one form or another, the \emph{Annalen} are naturally at your disposal."
\end{quote}

 Poincar\'e responded on April 4, 1882, after he received Klein's article on uniformization (\cite{HP-correspondance2} p. 111): ``Thank you very much for your note \cite{Klein-MA1882a} which you kindly sent me. The results which you state are of great interest to me for the following reason: I already found them some time ago, but I did not publish them, because I first wanted to clarify the proof. This is why I would like to know yours, when you will clarify them from your side."

The approach to uniformization, both by Poincar\'e and Klein, raised a difficulty, due to their use of the co-called \emph{continuity method}. This method was problematic, for several reasons. First of all, it involved a map between two spaces which are required to be manifolds (using modern terminology). For the spaces to which Poincar\'e and Klein applied it, this was not proved. In fact, this was even wrong, because one of the two spaces was Riemann's moduli space, which is not a manifold.
   Another problem was that the ``number of moduli" had not been made precise, and in order to be useful, this number had to be interpreted as a dimension (assuming the spaces involved are known to be manifolds ...), specially after Cantor found his famous bijection between the real line and the plane.  We shall talk about this in more detail in \S \ref{s:Brouwer}.   
   
In his 1884 memoir on Kleinian groups \cite{poin},  Poincar\'e's surveys the method of continuity as it was used by Klein  (\S VIII), and then, by  himself (\S XIV). He writes:
``Mr. Klein and myself have been led independently of one another to a method which allows one to prove that any Fuchsian type contains a Fuchsian equation, and that we can call the \emph{method of continuity}.  We made various applications of this method (see \emph{Comptes Rendus}, t. 92, p. 1200 and 1486; Klein, \emph{Mathematische Annalen} Bd. 19, p. 565; Bd. 20, p. 49 and Bd. 21 p. 211; \emph{Comptes Rendus}, t. 84, p. 1038)."
He then gives a brief presentation of this method. It is interesting to note that Poincar\'e describes this method as a general principle, not restricted to the case in hand. He starts with examples, namely, with maps between one-dimensional objects and then surfaces, making the distinction between ``closed" and ``open" surfaces. He writes: an \emph{open surface has a boundary, or frontier}.\footnote{This terminology of ``open" surface, as a surface with boundary, is used in Picard's \emph{Cours d'Analyse} (\cite{Picard-cours} p. 465). The reader should note there that this word does not have the same meaning today.} He then passes to the higher-dimensional case:
\begin{quote}\small
The same thing will happen, when $S$ and $S'$ are regarded as surfaces situated in the space of more than two dimensions, as \emph{Mannigfaltigkeiten} (multiplicities) of more than two dimensions, as the Germans say.

Suppose now that to each point $m$ of $S$ corresponds one and only one point $m'$ of $S'$  in such a way that the coordinates of $S'$ are analytic functions of those of $m$, unless $m$ goes to the boundary of the multiplicity $S$, in the case where this multiplicity is one. 
Suppose that there are no points of $S'$ which correspond to more than one point of $S$. If $S$ is a \emph{closed multiplicity}, then we can be sure that to any point of $S'$ corresponds a point of $S$. If, on the contrary, $S$ is an \emph{open multiplicity} having a boundary, or a frontier, then we cannot affirm anything.

\end{quote}
Poincar\'e then applies this result to a space of equivalence classes of Fuchsian groups. The equivalence relation is the following:
\begin{quote}\small Two Fuchsian groups belong to the same class when the generating polygon has the same number of sides and when the vertices are dispatched in the same number of cycles, in such a way that the cycles correspond one to one, and likewise for the angle sum of two corresponding cycles. 
\end{quote}
He then describes a map between a multiplicity $S$ and a multiplicity $S'$. The main question is whether  it is surjective.
\begin{quote}\small
For that, and by the preceding, it suffices that $S$ be a \emph{closed} multiplicity, and that {it has no boundary. This is not at all evident a priori.}

Indeed, among the groups of the same class, infinitely many could be \emph{limit groups}, corresponding to boundary points of the multiplicity $S$. These are groups whose generating polygon presents one or several infinitesimal sides. We can see, in effect, that we can always construct a generating polygon having a side which is as small as we want. In fact, it will not suffice for $S$ to be an open multiplicity, since, from \S IX of the \emph{Th\'eorie des groupes fuchsiens}, the same group may be generated by an infinity of equivalent polygons and it is possible that among these polygons we can always choose one whose sides are all greater than a given limit.

Thus, it is not evident that $S$ is a closed multiplicity and it is necessary to prove it, by a discussion which is special to each particular case, before asserting that to every point of $S'$ corresponds a point of $S$. This is what Mr. Klein omitted to do. \emph{There is here a difficulty which we cannot overcome in a few lines.}
\end{quote}
It is worth noting that about thirty years later, in a letter he wrote to Brouwer in January 1912, which we quote in \S \ref{s:Brouwer}, Poincar\'e refers to the Klein approach, to deal with the problem of singular points of Riemann's moduli space.

Poincar\'e then describes his own use of the method of continuity, in 
\S XIV,  after having established several lemmas. 
He starts with a linear order-two differential equation -- one of his favorite -- of the form 
\begin{equation}\label{e1}
\frac{d^2v}{dx^2}=\phi(x,y)v
\end{equation}
such that $x$ and $y$ satisfy an algebraic equation
\begin{equation}\label{e2}
\psi(x,y)=0.
\end{equation}
He argues that the space of Fuchsian groups which it leads to has no boundary. The reason is that if this space had a boundary, then this boundary would be of codimension one. But then he shows that any degeneration in that space falls on a subspace of codimension two and not one. He writes:
\begin{quote}\small
The parameters representing the moduli form a \emph{multiplicity} $M$ of dimension $q$ which is closed and has no frontier. Indeed, Riemann showed that if a multiplicity of dimension $q$ were limited by another multiplicity, then this new multiplicity would have dimension $q-1$. [...] The fact is that by varying the parameters of the class, we can only reach the frontier of $M$ by attaining certain singular points of this multiplicity, corresponding to the case where the type\footnote{The word ``type" refers to the classification of Fuchsian group.} $T$ is reduced to a simpler type $T'$. But this can happen only in two manners:

1) Two singular points collide; but this gives a complex condition, that is, two real conditions;

2) The genus of the algebraic relation (\ref{e2})  is lowered by one unit, that is, the curve represented by this relation. But this also gives two real conditions.
\end{quote}
Poincar\'e concludes that in both cases, ``the singular points  would form a multiplicity of dimension $q-2$, and therefore cannot form a frontier for the multiplicity $M$, whose dimension is $q$." Reasoning with limits of polygons, he argues that the degeneration would give a multiplicity of dimension $q-1$. Thus, the multiplicity has no boundary. 
This, according to him, settles the continuity method. 

  It is of highest interest that Poincar\'e examined the possibility of a boundary structure for his space of equivalence classes of Fuchsian groups, and that he also noticed that the dimensions of all the spaces involved (including the boundaries, if they exist) should be even. The reader will recall that there was no higher-dimensional complex structure yet present in the discussion.

The rest of the correspondence between Poincar\'e and Klein on the method of continuity is also rich. On May 7, 1882, Klein writes (\cite{HP-correspondance2} p. 113):\begin{quote} \small
I recently read your \emph{Comptes Rendus} Note of April 10 \cite{P1882}. I was all the more interested that I think that your present considerations are close to mine, as for what regards the method. I prove my theorems using the \emph{continuity}, based on the following two lemmas: 
 
1) With any discontinuous group is associated one Riemann surface; 
 
2) Only \emph{one} of these groups belongs to a Riemann surface which is conveniently cut [zerchnittene] (as long as there is an associated group).
\end{quote}

Poincar\'e responds, in a letter dated May 12, 1882 (\cite{HP-correspondance2} p. 114): 
\begin{quote}\small
I think, as you do, that our methods are very close and they differ less by the general principle than by the details.  \end{quote}
In a letter dated May 14, 1882 (\cite{HP-correspondance2} p. 116), Klein explains his method (\cite{HP-correspondance2} p. 114): 
\begin{quote}\small  I would like to explain to you, in two words, how I use the ``continuity." [...] We have to prove, above all, that two manifolds [Mannigfaltigkeiten] that we compare -- the set of systems of substitution that we consider, and on the other hand, the set of Riemann surfaces that exist effectively -- not only have the same number of dimensions ($6g-6$ real dimensions), but are \emph{analytic} with \emph{analytic} frontiers [Grenze] (in the sense introduced by Weierstrass) [...] But now it turns out this relation is \emph{analytic} and even, as it follows from the two propositions, an analytic relation \emph{whose functional determinant is never zero}.
\end{quote}  
Poincar\'e's response is dated May 18, 1882 (\cite{HP-correspondance2} p. 117):
\begin{quote} \small 
It is probable that we do not establish by the same method the analytic character of the relation which relates our two \emph{Mannigfaltigkeiten} which you talk about. For me, I relate this to the convergence of my series; but it is clear that we can obtain the same result without such a consideration. 
\end{quote}
We refer the interested reader to the complete set of letters.

 There is a historical section in Poincar\'e's 1882 article \cite{P1882a} in which he reports that results on discontinuous groups of fractional linear transformations follow from work of Hermite, and he then mentions works on this subject by Dedekind, Fuchs, Klein, Hurwitz and Schwarz.    In his \emph{Analysis of his own works} \cite{Poin-Acta}, p. 43, he describes elliptic functions as univalent functions which are periodic.\footnote{The word \emph{univalent} does not mean injective, as the function is periodic. One has to take into account the fact that the function is automorphic.} The periodicity here is the important property, and it is a property of group invariance. Poincar\'e always thought in terms of groups. He notes that the periods form a discrete group, and he mentions the fundamental domain for this action, a parallelogram in the complex plane. ``The knowledge of the function in one of the parallelograms implies its knowledge everywhere in the plane." He then considers the situation where the fundamental region is a ``curvilinear polygon." He addresses the question of finding the conditions under which the group generated by reflections along sides on such a polygon is discontinuous, and he relates this question to non-Euclidean geometry. This   leads him to consider theta functions which are not periodic, but ``are multiplied by an exponential when the variable is augmented by a period." (p. 46).  This is the origin of the subject of automorphic forms.        

Poincar\'e's theory of Fuchsian groups and Fuchsian functions, as well as his more general theory of automorphic forms, led to many developments in connection with the study of moduli in works of Bers \cite{Bers-auto}, Kra \cite{Kra1972}, Sullivan \cite{Sullivan},  McMullen \cite{McM} and others, and it is not possible for us to develop this subject here.

 We saw how much Poincar\'e and Klein had close interests. Two events recounted by the two men can also be put in parallel. The first one, described by Poincar\'e, is well known, and it occurred in the Summer of 1880, at the time of his correspondence with Klein, when he made his lightning discovery that the discrete groups that arise from the solutions of differential equations he was studying are the  groups of isometries of the hyperbolic plane. Poincar\'e recounted this discovery, which he considered as one of his major ones, in a text written in 1908, that is, 28 years after the episode, at a talk he gave at the Soci\'et\'e de Psychologie de Paris. It is published in his book \cite{PK-Science}, and it is worth recalling, and put in parallel with a text written by Klein. Poincar\'e writes:

\begin{quote}\small
During fifteen weeks, I was working hard to prove that there was no function analogous to the one which I had called since then \emph{Fuchsian function}. At that time, I was very ignorant. Every day, I sat at my work table. I used to spend one hour or two. I tried a large number of combinations, but I obtained no result. One evening, I drank black coffee, contrary to my habit. I could not sleep, because of ideas flocking up. I felt them colliding, until I hooked two of them, and they formed a sort of stable combination. In the next morning, I established the existence of a class of Fuchsian functions which arise from the hypergeometric series. I had only to write the results and this took me a few hours.

Then, I wanted to represent these functions as quotients of two series. The idea was perfectly conscious and deliberate. I was guided by the analogy with elliptic functions. I was wondering what should be the properties of these series, if they existed, and I arrived without difficulty to construct the series I called thetafuchsian.

At that moment, I was leaving Caen, where I used to live at that time, to take part in a geological competition organized by the \emph{\'Ecole des Mines}. The adventures of the trip made me forget about my mathematics. When we reached Coutance, we got in a local train for some tour. At the moment I set my foot on the footboard, the idea came to me, without any preparation from previous thoughts, that the transformations that I used for the Fuchsian functions were identical to those of non-Euclidean geometry. I did not check it, and in fact, I wouldn't have got time to do it since as soon as I sat in the train, I continued a conversation which I had started before; but I immediately got a complete certitude. When I came back to Caen, I checked the result at my leisure, although I was sure of it.

I started to study some questions in number theory, without real visible result and without suspecting that this could have the slightest relation with my research. I became sick of my failure, and I went for a few days to the seaside, where I began thinking about something else. One day, while I was wandering over the cliff, an idea stoke me, again, with the same immediate suddenness and certainty, that the arithmetic transformations of the indefinite ternary quadratic forms were identical to those of non-Euclidean geometry.

When I returned to Caen, I thought again of that result and I deduced from it some consequences. The example of quadratic forms showed me that there were other Fuchsian groups than those which correspond to the hypergeometric series, and I saw that I could apply to them the theory of thetafuchsian functions, and, therefore, that there exist other thetafuchsian functions than those that derive from the hypergeometric series, which were the only ones I knew until that time. Naturally, I  tried to form all these functions. I made  a systematic siege of them. [...]

After that, I went to the Mont-Val\'erien, where I was going to do my military service. I therefore had very different concerns. 
One day, while I was crossing the boulevard, the solution of the difficulty appeared to me suddenly. I did not try to go immediately deeper into it and it was only after my military service that I came back to the problem. All the elements were in my possession; I had only to assemble them and to put some order in them. I then wrote my memoir, without a break and without any difficulty.

\end{quote}

 Klein's episode is related by N\"orlund, in his edition of the correspondence between Klein and Poincar\'e \cite{Noerlund}. It concerns the discovery of the uniformization theorem. Klein writes:

\begin{quote}\small In order to follow the advice that the doctors gave me at that epoch, I decided to come back, at Easter 1882, on the shores of the North Sea, this time at Norderney. I wanted to write up, in this quiet place, the second part of my book on Riemann, that is, to have under a new written form, the existence theorem for algebraic functions on a Rieman surface. In fact, I could not stand this stay more than eight days; life there was too morose, since it was impossible to go out because of the violent storms. Thus, I decided to come back as soon as possible to my home in D\"usseldorf. During the last night of my stay there, the night of March 22 to 23, which I spent sitting on a couch because of my asthma attacks, suddenly, around 2:30, the central theorem [Zentraltheorem] appeared to me, in the way it is sketched in the figure of the 14-gon [14-Eck] in Vol. XIV of the \emph{Mathematische Annalen} p. 126.\footnote{This is Klein's article \cite{Klein1879}.} The next morning, in  the stagecoach which at that time circulated between Norden and Emden, I pondered about my finding, examining again every detail. I knew I had found an important theorem. As soon as I arrived to D\"usseldorf, I wrote a memoir, dated March 27, sent it to Teubner\footnote{This is Klein's article \cite{Klein-MA1882a}.} and I asked that the proofs be sent to Poincar\'e and to Schwarz.
 
I already recounted how Poincar\'e reacted in the \emph{Comptes rendus} \cite{P1882} on April 10. To me, he wrote: \emph{The results you state are of great interest for me, for the following reason: I already found them, some long time ago} [...] He never said precisely how and since when he knew them. One can easily understand that our relations became tense. 
Schwarz, on  his side, first thought, after an insufficient counting of constants, that the theorem must be false. But later on, he brought a major contribution to the new methods of proof.
 
But in reality the proof was very difficult. I used the so-called continuity method, which assigns to the set of Riemann surfaces of a given [genus] $p$ the corresponding automorphic group of the unit circle. I have never doubted of the method of proof, but I always stroke against my lack of knowledge in the theory of functions, or against the theory of functions itself, of which I could only tentatively assume the existence, and which in fact, was obtained only thirty years later (1912) by Koebe \cite{Koebe1912}.
\end{quote}

Klein's systematic presentation of his work on uniformization is contained in the book he authored with Fricke \cite{fk} on which we comment in \S \ref{s:book} below. In 1907, Poincar\'e \cite{Poincare1907} and Koebe\footnote{Koebe (1882-1945) wrote his dissertation under Schwarz in 1905. After he proved the uniformization theorem, Koebe practically spent the rest of his career in trying to improve and expand his proof, and in doing so he discovered the notion of planarity in both the topological and complex analytic contexts.  Koebe also spent a lot of energy in trying to convince others that the proof of the uniformization which Klein gave in 1882 cannot be made rigorous. We shall mention him again in \S \ref{s:Brouwer} dedicated to Brouwer. His work is described by Bieberbach  in the article \cite{Bieber} and by Freudenthal in the  \emph{Dictionary of Scientific Biography}.}  \cite{Koebe1907} \cite{Koebe1909} gave, independently, a proof of the general uniformization theorem.  Both proofs contain very new ideas and do not use the continuity method. We refer the reader to the exposition in \cite{uni}. Other proofs were obtained later on, by several authors.

It is good to conclude this section by Bers' comments on uniformization. He writes in his survey on this question (\cite{Bers-Hilbert} p. 509): 
  \begin{quote}\small
 A significant mathematical problem, like the uniformization problem which appears as No. 22 on Hilbert's list, is never solved only once. Each generation of mathematicians, as if obeying Goethe's dictum,\footnote{Was du ererbt von deinen V\"atern hast, erwirb es, um es zu besitzen.} rethinks and reworks solutions discovered by their predecessors, and fits these solutions into the current conceptual and notational framework. Because of this, proofs of important theorems become, and if by themselves, simpler and easier as time goes by -- as Ahlfors observed in his 1938 lecture on uniformization. Also, and this is more important, one discovers that solved problems present further questions.
 
 [...]  The personal element is also fascinating. It involves some of the most illustrious mathematicians of that time: Schottky, about to conjecture, in 1875, a fairly general uniformization theorem, but deflected by the authority of Weierstrass (according to Klein), the rivalry between Klein, then at the height of fame and productivity, and the yet unknown Poincar\'e (see their correspondence from the years 1881--1882). Schwarz, suggesting, in a private communication, two methods for proving the main uniformization theorem (one using the universal covering surfaces, the other involving the partial differential equation $\Delta u = e^{2u}$), Hilbert, reviving the interest in the problem by his Paris lecture, and several years later creating a new tool for uniformization by ``rehabilitating" the Dirichlet principle, Brouwer, embarking on  his epoch making topological investigations in order to put the original ``continuity method" of Klein and Poincar\'e on firm foundations, Poincar\'e, returning to the uniformization problem after a quarter of a century and finally achieving a full solution, but having to share this honor with Koebe. 
  Koebe ``went on to explore, with the most varied methods, all facets of the uniformization problem." The (slightly rephrased) quotation is from the 1955 edition of Weyl's celebrated  \emph{Idee der Riemannschen Fl\"ache}. [...]
   
 The modern developments in uniformization, which began after a period of hibernation, utilize quasiconformal mappings and the recent advances in the theory of Kleinian groups. Quasiconformal mappings give new proofs of classical uniformization theorems, akin in spirit though not in technique, to the old ``continuity method," and also proofs of new theorems on simultaneous uniformization. The theory of Kleinian groups permits a partially successful attack on the problem of describing all uniformizations of a given algebraic curve. An unexpected application of simultaneous uniformization is Griffith's uniformization theorem of $n$-dimensional algebraic varieties, which answers a question also raised by Hilbert in the 22nd problem. (For non-Archimedean valued complete fields there are uniformization theorems due to Tate, for elliptic curves, and to Mumford, but the present author is not competent to report on this work).
  \end{quote} 
 An exposition of these later uniformization theories would need a longer article. Let conclude by recalling that the solution by Perelman of the most general uniformization problem so far, namely, Thurston's geometrization,  uses as an essential ingredient, an evolution equation and the theory of partial differential equation, which is completely in the tradition of Riemann and Poincar\'e in their use of the Laplace equation and the Dirichlet principle. By the way, this is a field where Bers worked at the beginning of his career; cf. the article by Abikoff and Sibner in this volume \cite{AS}.

\section{Fricke and Klein's book}\label{s:book}

   The two-volume book \emph{Vorlesungen \"uber die Theorie der automorphen Funktionen} \cite{fk} by Fricke\footnote{Karl Emmanuel Robert Fricke (1861-1930) studied mathematics, philosophy and physics at the universities of G\"ottingen, Z\"urich, Berlin and Strasbourg. He obtained his doctorate under Klein at the University of Leipzig in 1885. After this, he taught at the Braunschweig gymnasium and then worked as a private tutor; During this period he kept his relations with Klein, with whom he wrote the two-volume treatise \emph{Vorlesungen \"uber die Theorie der elliptischen Modulfunctionen}, published in 1890 and 1892 respectively. It was  between the publication of these two volumes that Fricke obtained an academic position, first at the University of Kiel (1891), then at the University of G\"ottingen (1882) where Klein was teaching, and finally, in 1894, at the Braunschweig Polytechnikum, succeeding Dedekind. Fricke married in 1894 a niece of Klein. We shall talk more thoroughly  about another two-volume treatise which he wrote with Klein, \emph{Vorlesungen \"uber die Theorie der automorphen Funktionen} \cite{fk}.} and Klein is in the stream of the study made by Klein and Poincar\'e of polygons and spaces of polygons equipped with group actions. The book may be considered as the first comprehensive study of groups of discrete motions of the Euclidean, elliptic and hyperbolic planes. The hyperbolic case constitutes the largest part of the book, since the classification of group actions in that case is much more involved.  
 The authors quote several articles by Klein, in particular \cite{Klein1885} and \cite{Klein1886}. They describe families of  polygons equipped with group actions using Riemann's term \emph{Mannigfaltigkeit}. The space of groups includes those that act with fixed points, and thus, spaces of orbifold surfaces are also considered.

Chapter 2 of Vol. 1 carries the title \emph{The canonical polygons and the parameters for rotation groups}.\footnote{Fricke and Klein use the word  ``rotation" for an  isometry of any one of the three constant-curvature planes. This is consistent with the other writings of Klein. Motivated by projective geometry, Klein considered an isometry as a rotation whose center might be in the space, or at infinity, or beyond infinity. The authors of \cite{fk} make use of Klein's projective model, with a conic at infinity, a point of view which Klein had  adopted in his articles \cite{Klein-Ueber} and \cite{Klein-Ueber1}. Different choices of conics lead to different geometries.} 
  For each type of polygon with geodesic sides in the hyperbolic plane, the authors consider a space of polygons. They call a ``canonical polygon" a polygon obtained through a system of cuts that forms a basis of the fundamental group of the surface. Canonical polygons and their moduli are developed systematically on pages 284-394 of Vol. 1. In a note on p. 284, it is specified that isomorphic groups should not be considered as different, which shows that the space  constructed is Riemann's moduli space and not Teichm\"uller space. The authors claim that the moduli of the canonical polygons are \emph{de facto} the moduli of the discrete groups.  They introduce traces of matrices as new parameters. They describe in detail  several special cases, in particular the surfaces of type $(0,3)$ (pair of pants) and $(1,1)$ (torus with one hole) and they show that in these cases the moduli space is in bijection with the product $[0,1]^3$ with some faces removed, which they call a ``cube" (p. 290 and 302 of Vol. 2). Other surfaces of low genus and a small number of boundary components are also considered in detail. After the special cases, the authors develop a general theory by gluing the small surfaces along their boundaries and adding one parameter for each gluing. This theory is a precursor of the later rigorous theories of gluing hyperbolic surfaces developed Fenchel and Nielsen.\footnote{ After topological tools and quasiconformal mapping theory were made available, a more accurate construction of the canonical polygons considered by Fricke and Klein in \cite{fk} was made by Linda Keen in \cite{Keen1}. The trace coordinates that are used in \cite{fk} were also made more precise by Keen in \cite{Keen2}. }  In the book by Fricke and Klein, one also finds the idea of a marking of a surface, considering that the polygons are marked by the choice of the elements of the fundamental group that correspond to the curves along which it is cut.  A count for the number of moduli is given, and the result is the same as Riemann's.  There is a ``modular group" acting on spaces of polygons by identifying parameters of hyperbolic polygons which lead to the same hyperbolic surface (cf. p. 308 of Vol.  2 of \cite{fk}). The authors claim that this action is given by rational functions. The quotient by this modular group is mapped injectively into the moduli space of hyperbolic surfaces.  
The authors also aim to give a systematic presentation of the method of continuity for its use in the uniformization theorem (although
 without substantial
improvement). Brouwer, Koebe, Bieberbach and  others were very suspicious about these attempts.\footnote{For instance, in a letter sent on May 22, 1912 to Brouwer (\cite{Brouwer-correspondence} p. 148), Bieberbach writes: ``In Fricke, the proof in the boundary circle case seems to rely completely on the fact that only one fundamental domain belongs to a system of generators. But this is not satisfied."}  This method is exposited on pages 402ff. of Vol. 2.
The authors  consider two spaces, namely, a space of unmarked Fuchsian groups of a fixed type,  and a space of Riemann surfaces of a fixed topological
type (the Riemann moduli space). They define a natural map between them, by assigning to each Fuchsian group  
the quotient Riemann surface. Such a surface is identified, following Riemann, with the field of meromorphic functions it carries, and a Fuchsian group is defined accordingly as the space of automorphic functions on it. Several complications arise in making such a theory rigorous. We already mentioned Poincar\'e's remarks on degenerations of polygons, on
 the difference between open and closed multiplicities, etc. These problems are not solved in the book of Fricke and Klein.

 We have seen that Klein, like Riemann, was more interested in communicating ideas than writing proofs, and his articles and books have a conversational style. His book with Fricke is an exception: it is much more technical and difficult to read than the rest of his writings. The notions and results are sometimes presented in a confusing way, and everybody who has tried to figure out what this book contains precisely knows that this task is almost impossible.\footnote{We quote for a letter from Blumenthal to Brouwer, dated August 26, 2011 (\cite{Brouwer-correspondence} p. 99): ``I would like now to deal quickly with the automorphic functions [...] In Klein's article in \emph{Mathematische Annalen} 21 the problem is completely and clearly formulated from a set theoretic point of view,\footnote{We recall that at that time the expression  ``point-set theory" denoted general topology.}  even though the answer given there does not satisfy the standard for rigor. I strongly advise you to go through the matter there, and \emph{not} in the fat Fricke and Klein, where one trips again and again over details that obscure the general idea." Abikoff, in a correspondence with one of the authors of the present article, recalls the following: ``In 1971, Ahlfors and I were at the Mittag-Leffler Institute. I was trying to read Fricke-Klein and mentioned to him that my German wasn't up to dealing with that text. He said, `Your German is maybe 10 \% of the problem'."} This is certainly due to Fricke's style, and in fact, this is one the reasons why the results in this book are attributed to Fricke.\footnote{People talk more thoroughly about the ``Fricke space" and not the ``Fricke-Klein" space, which is somehow unfair to Klein because, as Fricke acknowledges at several points, most of the ideas are due to Klein. In the two volumes, the foreword is signed by Fricke alone, but Fricke emphasizes Klein's contribution.}  
 
 Bers claimed that Fricke and Klein,  as well as Poincar\'e, were dealing with Teichm\"uller space. For instance, in his paper \cite{bers-uni}, he writes: ``The space $\mathcal{T}_p$ appears implicit in the early continuity arguments by Klein and Poincar\'e; it has been constructed as a manifold of $6p-p$ dimensions by Fricke and Klein \cite{fk} (who proved it to be a cell) and by Fenchel-Nielsen \cite{FN}."
A statement that Fricke defined the space of marked hyperbolic or Riemann surfaces is also made in the paper \cite{bg} by Bers and Gardiner. Ahlfors is somehow more dubious, but on the other hand, he also talks about the complex structure of Teichm\"uller space, and not only the topological definition of that space. He says in his \emph{Collected papers} edition,  Vol. 2 p. 122, commenting his paper \emph{The complex analytic structure of the space of closed Riemann surfaces} (1960): ``An incredibly patient reader of Fricke-Klein, a two-volume 1300-page book, might have been able to discern that Fricke had anticipated Teichm\"uller's idea of a space $\mathcal{T}_g$ of marked Riemann surfaces of genus $g>1$ with a $3g-3$ dimensional complex structure. Although in 1959 nobody had yet reconstructed Fricke's proof, it was generally believed that such a structure exists and that the elements of the Riemann matrix are holomorphic functions with respect to that structure."
 The reader will compare this with the precise results of Teichm\"uller that are surveyed in \S \ref{s:Teich} of this paper.

    There are many  modern  developments of the theory developed in Fricke and Klein's book, but it is not possible to review them here. They include discrete subgroups of Lie groups, number theory, automorphic forms, and higher Teichmuller theory.

  \section{Brouwer and the method of continuity}\label{s:Brouwer}

  When Poincar\'e died in 1912. Brouwer, who was 31, became the greatest living topologist. His doctoral dissertation, titled \emph{On the foundations of mathematics}   and his inaugural dissertation, \emph{The nature of geometry}, which he had defended in 1907 and 1909 respectively, like Riemann's doctoral dissertation and inaugural lecture, are huge programs. Brouwer's doctoral dissertation has three parts: I. The construction of mathematics; II. Mathematics and experience; III. Mathematics and logic. In his inaugural lecture, Brouwer expresses the fact that topology is the most fundamental part of mathematics. Both dissertations, have a marked logico-philosophical side. The reader can read the English versions of the two dissertations in \cite{Brouwer-collected}.\footnote{The two dissertations are included in Vol. 1 of Brouwer's collected works \cite{Brouwer-collected}, that is, the volume which contains the philosophical writings, and not in Vol. II which contains the mathematical writings.} Chapter 3 of Vol. 1 of \cite{Van-Dalen} contains commentaries on these works.

  Brouwer had an uncommon personality. He was extremely self-disciplined, with a rigorous life, a schedule for every day and an extensive reading program. He also had high moral standards and an enormous will to avoid all sort of mistakes. Furthermore, his view of mathematics was unusual. Knowing this is important for understanding his contentious relations with mathematicians.  In a letter to Poincar\'e in 1912, asking him to write a recommendation letter for his student Brouwer, Korteweg writes: ``The goal is to find for him an appropriate place in one of our universities; this will not be easy, given his distinguished but very peculiar personality."\footnote{The same year, Brouwer was appointed professor at the University of Amsterdam and he was elected member of the Academic of Sciences of Netherlands.}
  
A review of Brouwer's view on mathematics is contained in \cite{Van-Dalen}, which is an extremely interesting  two-volume mathematical, philosophical and personal biography.  We also recommend to the reader the two-volume collected works of Brouwer \cite{Brouwer-collected} and the valuable volume of correspondence of this mathematical giant \cite{Brouwer-correspondence}.
 
   General topology and the topology of manifolds were undergoing birth at the time Brouwer began his work, and his work was crucial in making this field grow.   Among his many contributions to the field are a precise definition of the notion of dimension and its topological invariance, the notion of degree of a map, several fixed point theorems and a rigorous proof of the Invariance of Domain Theorem.\footnote{There were several attempts prior to Brouwer's work to prove the invariance of dimension (with a notion of dimension which was still not fully satisfying).  The need to have rigorous results in this domain was particularly felt after Cantor's example (1878) of his bijection between a segment and a square. Several studies of dimension and dimension invariance were made, in particular cases, by Jordan and by Schoenflies.  L\"uroth proved, for some special values of $m$ and $n$, that if we have a bicontinuous bijection between two continua of dimensions $m$ and $n$, then $n\geq m$. He obtained this result for $m=1, n\geq 2$, and then for $m=2, n=3$ and finally for $m=3, n=4$. \cite{Luroth1899}  \cite{Luroth1907}.}  His new tools included approximation by simplicial mappings and an extensive use of homotopy and the degree of a map.
    
 For several of his discoveries, Brouwer was influenced by Poincar\'e,  of whose works (including the philosophical essays) he was a devoted reader.\footnote{Poincar\'e is quoted ten times in Brouwer's doctoral dissertation (especially referring to his philosophical writings).} Like Poincar\'e, Brouwer was interested in automorphic functions and uniformization, but his stress was on the topological aspects. His style was different though, since Poincar\'e was satisfied by intuitive arguments and by sketches of proofs whereas Brouwer was obsessed by precise definitions and by rigorous proofs. In trying to make rigorous the method of continuity as used by Poincar\'e and Klein, Brouwer insisted on the fact that the integral parameter that was involved should be interpreted as a topological dimension. This required a precise definition of dimension, the proof of its invariance, and, of course, a precise definition of the topologies of the two spaces involved: the space of equivalence classes of polygons equipped with group actions, and Riemann's moduli space. 

Brouwer proved the first version of the Invariance of Dimension Theorem in his paper 1911 \cite{B1911}, and he obtained stronger results later on. Also implicit in his paper \cite{B1911} is the notion of degree of a mapping, an idea that he made more explicit in \cite{B1912}.\footnote{Brouwer obtained the idea of a degree of a map by reading a work of Poincar\'e.  In a letter to Hadamard on December 24, 1909 (\cite{Brouwer-correspondence} p. 62):, he writes:  ``Cher Monsieur, Thank you very much for having pointed out Mr. Poincar\'e's memoirs on algebraic vector distributions. [...] I have had another idea. First we remark that if we adapt the concept of an `index' (quoted from the first memoir, p. 400) to general continuous vector distributions, Corollary I of p. 405 becomes the following: \emph{If the singular points are finite in number, each of them has a finite index and the algebraic sum of all the indices is equal to 2}."} We refer to Koetsier \cite{Koetsier-Hand} and to the comments in Brouwer's Collected Works edition by Freudenthal for summaries of Brouwer's works on topology.

   In the intense period 1911-1913, an extensive and interesting exchange of letters took place between Brouwer, Koebe, Klein, Blumenthal, Baire, Poincar\'e, Fricke, Hilbert, Scheonflies, Hadamard and others, concerning the birth of the topological notion of dimension and the method of continuity.

Brouwer gave a memorable talk at the automorphic functions-session of the annual meeting of the German Mathematical Society (DMV) which was held in Karlsruhe on September 27 to 29, 1911. He presented there Poincar\'e's original proof of uniformization, in which he filled two gaps, and a version of his Invariance of Dimension Theorem, which rescued the method of continuity. This talk was a source of disagreement and fighting with several mathematicians, that lasted several years.  In particular, a conflict between Koebe and Brouwer started just after the meeting. The debates caused by these conflicts are interesting for us, because they concern spaces of Riemann surfaces, and because they contributed making the methods and problems that concern us more precise.
 The technical issues involved are summarized in letter from Brouwer to Fricke, dated December 22, 1911, and we quote it in its entirety because it is interesting (\cite{Brouwer-correspondence} p. 116): 
\begin{quote}\small
Dear Geheimrat,
With reference to our last conversation I inform you about some remarks related to the topological difficulties of the continuity proof, which I have presented at the meeting of Naturforscher in Karlsruhe.
 
Let $\kappa$ be a class of discontinuous linear groups of genus $p$ with $n$ singular points and with a certain characteristic signature; for this class the \emph{fundamental theorem of Klein} holds if to every Riemann surface of genus $p$ that is canonically cut and marked with $n$ points there belongs \emph{one and only one} canonical system of fundamental substitutions of a group of class $\kappa$.
 
In the continuity method, which Klein uses to deduce his fundamental theorem, the six following theorems are applied. 
 
\begin{enumerate}
\item The class $\kappa$ contains for every canonical system of fundamental substitutions that belongs to it \emph{without exception} a neighborhood that can be represented one to one and continuously by $6p-6+2n$ real parameters.

\item During continuous change of fundamental substitutions within the class $\kappa$ the corresponding canonically cut Riemann surface likewise changes continuously.
\item Two different canonical systems of fundamental substitutions of the class $\kappa$ cannot correspond to the same cut Riemann surface.
\item When a sequence $\alpha$ of canonically cut Riemann surfaces with $n$ designated points and genus $p$ converges to a canonically cut Riemann surface with $n$ distinguished points and genus $p$, and when each surface in the sequence $\alpha$ corresponds to a canonical system of fundamental substitutions of the class $\kappa$; then the limit surface likewise corresponds to a canonical system of fundamental substitutions of the class $\kappa$.

\item The manifold of cut Riemann surfaces contains for every surface belonging to it \emph{without exception} a neighborhood that can be one to one and continuously represented by $6g-6 +2n$ parameters.

\item In the $(6g-6+2n)$-dimensional space the one to one continuous image of a $(6g-6+2n)$-dimensional domain is also a domain.
\end{enumerate}

I am ignoring here Theorems 1, 2, 3, 4. For the case of the boundary circle they have been completely treated by Poincar\'e in Vol. 4 of \emph{Acta Mathematica}. For the most general case only Theorems 3 and 4 await an exhaustive proof. In this matter also this gap will be filled in by Mr. Koebe in papers that are to appear soon.

  Theorems 5 and 6  are those which constitute the topological difficulties of the continuity proof that are emphasized in your book about automorphic functions \cite{fk}.  However, of these, Theorem 6 is settled by my recent article  \emph{`Beweis der Invarianz des $n$-dimensional Gebiets} \cite{Brouwer-gebiets}, whereas the application of Theorem 5 can be avoided by carrying out the continuity proof in the following modified form:

We choose $m>2p-2$ and consider on the one hand the set $M_g$ of automorphic functions belonging to the class $\kappa$ that only have simple branching points and with $m$ simple poles in the fundamental domain, and on the other hand the set $M_f$ of Riemann surfaces covering the surface, of genus $p$, with $n$ signed points, and with $m$ numbered leaves and with $2m+2p-2$ numbered simple branching points not at infinity, for which the sequential order of the leaves and the branching points correspond to the canonical relations in the sense of L\"uroth-Clebsch.\footnote{The editor of \cite{Brouwer-correspondence} writes, in a footnote: Crossed out footnote of Brouwer: `Two of these surfaces are considered identical if and only if the corresponding not-covered surfaces can be mapped so much similarly onto each other that corresponding return cuts and stigmata behave the same with respect to the construction of the covering surface."}

The set $M_f$ constitutes a continuum, and possesses for each of its corresponding surfaces \emph{without exception} a neighborhood which is one-to-one and continuously representable by $4p-8+2n+4m$ real parameters.

For an arbitrary automorphic function $\phi$ belonging to $M_g$ there exists in $M_g$ a neighborhood $u_{\phi}$ which can be determined by $4p-8+2n+4m$ real parameters; these parameters are the $m$ complex places of the poles in the fundamental domain, the $m-p-1$ complex behaviors of the $m-p$ arbitrary pole residues, and the $6p-6+2n$ parameters of the canonical systems of fundamental substitutions. The value domain of the parameters belonging to $U_\phi$ constitutes a $(4p-8+2n+4m)$-dimensional domain $w_\phi$.

With the function $\phi$ there corresponds a finite number of surfaces belonging to $M_\phi$. Furthermore we conclude from Theorems 1, 2, 3 and the remark that possible birational transformations into itself not only for the single Riemann surface, but also for the totality of Riemann surfaces belonging to $u_\phi$, cannot become arbitrarily small, that with a  sufficiently small $w_\phi$ in $M_f$ there corresponds a finite number of one to one and continuous images, and hence because of Theorem 6 a domain set. However, \emph{then the total set $M_g$ in $M_f$ corresponds with a domain set $G_f$ too.}

Now we formulate Theorem 4 in the following form:

When a sequence of canonically cut surfaces of $M_f$ converges to a canonically cut surface of $M_f$ and when each surface of the sequence corresponds to a canonical system of fundamental substitutions of the class $\kappa$, then the limit surface also corresponds to a canonical system of fundamental substitutions of the class $\kappa$.

This property immediately entails that the domain set $G_f$ cannot be bounded in $M_f$, and \emph{hence it must fill the whole manifold $M_f$.} This proves the fundamental theorem for every Riemann surface of genus $p$ on which there exist algebraic functions with more than $2p-2$ simple poles and with exclusively simple branching points, i.e. just for any Riemann surface of genus $p$.
\end{quote}

In this letter, both spaces $M_f$ and $M_g$ are \emph{not} the Riemann moduli space; their dimension is greater since $m>2p-2$. The space $M_f$  is a space of coverings of the Riemann sphere with some specified branching above  the moduli space. The space $M_g$ also sits above the Riemann moduli space, and it corresponds to spaces of  hyperbolic surfaces (or Fuchsian groups)
with certain additional structure obtained by imposing conditions on the automorphic functions of the group. The idea of introducing these spaces is to enhance the structures of Riemann surfaces and hyperbolic surfaces in order to get manifold structures on the resulting  moduli spaces.  This is a precursor to the idea of \emph{marking}.

 In fact, Koebe was not prepared to publish the paper which is announced in Brouwer's letter. On February 2, 1912, he wrote to Brouwer, expressing various thoughts, including a critique of Fricke and Klein's book \cite{fk}:
``Poincar\'e rather represents the interpretation of \emph{closed} continua by adding \emph{limits of polygons} [...] Poincar\'e recently informed me in conversation that the continuity method cannot be used \emph{at all} if one wants to prove the no-boundary-circle theorem, because these manifolds are \emph{not closed}. [...] Also \emph{Fricke-Klein's} `Vorlesungen \"uber automorphe Funktionen' has adopted extensively the view of Poincar\'e about \emph{closedness}." He then sent a text which he asked Brouwer to insert in his paper, in which he says the following: ``[...] The proofs found by Mr. Koebe extend to the case of boundary circle uniformization, the only one considered by Poincar\'e, and imply a \emph{life giving} advance, because of the liberation from the thoughts introduced by Poincar\'e and copied by Klein-Fricke about polygonal limits and closed continua, an advance which is at the same time a return to Klein's old standpoint of non-closed continua which was vigorously attacked by Poincar\'e."

Brouwer responded to Koebe, on February 14, 1912: ``Fortunately I am still in possession of the abridged text of my Karlsruhe talk, which I enclose, so that you can no longer maintain that I used in Karlsruhe in the talk or in discussion the `closed'  manifolds of Poincar\'e! That you could make such a statement only proves that modern set theory must be absolutely unfamiliar to you. For, the elaborations of Poincar\'e who works with the so-called `closed manifolds' are pure balderdash, and can only be excused by the fact that at the time of their formulation there was not yet any set theory." In his notes on Brouwer's collected Works, Vol. II, Freudenthal writes: ``The essential point is that Poincar\'e, in 1884, uses the closed `manifold' of Riemann surfaces, which Brouwer calls \emph{Bl\"odsinn} [non-sense] in his letter to Koebe; whereas Klein [Neue Beitr\"age zur Riemann'schen Functionentheorie, \emph{Mathematische Annalen}, 21 (1882), 41-218] does this with the open variety of cut Riemann surfaces, the way Brouwer did it later."

In a letter to Hilbert, written on March 3, 1912, Brouwer writes: ``Koebe apparently does not understand that a not one-to-one but continuous specification of a set by $r$ real parameters does not guarantee at all that this set is an $r$-dimensional manifold without singularities. [...] Koebe moves in a vicious circle, because on the one hand he demands from me that I extensively praise his paper which hasn't appeared yet, on the other hand he tries to prevent me from seeing this article [...] That the planned note of Koebe doesn't contain any falsehoods or insinuations concerning me, is, by the way, more in Koebe's interest than in mine, because in my eventual refutation I will probably not be able to avoid to disgrace him irreparably."

At the same period, a conflict started between Brouwer and Lebesgue, regarding the notion of dimension.\footnote{One may recall here that in 1877, Cantor gave his famous bijection between a line at a plane \cite{Cantor1878}. Before that, it was assumed that the dimension of a space means the number of coordinates necessary to parametrize a point in that space. Cantor's result showed that this is problematic. Schoenflies writes in \cite{Schoeflies1908} that when the first example of a space-filling curve was given, ``the geometers felt that the ground on which their doctrine was founded was shaking."  

At a bout the same time, Peano gave of a space-filling curve. All this showed that geometers and analysts had to revise the definition of the notion of ``curve" that they were using, in particular in the theories of Cauchy and of Riemann. This also confirmed the fact that the method of continuity, which was already suspicious, was not correct if not applied precisely. The response was given by Brouwer in his 1911 paper \cite{B1911}.  For the historical development of the notion of dimension, the interested reader can refer to the papers \cite{Johnson1} and \cite{Johnson2}.} It is related in Freudenthal's paper \cite{Freudenthal1975}. Blumenthal (who was the managing editor of \emph{Mathematische Annalen} at that time) met Lebesgue in Paris, during the Christmas holidays of 1910, and told him about Brouwer's discovery of the invariance of dimension. Lebesgue told Blumenthal that he already knew this result, and that he had several proofs. Blumenthal asked him to send him a paper containing one of his proof so he could publish it in the \emph{Annalen}.  Lebesgue provided his paper, and the proof, based on Lebesgue's so called \emph{paving property}, was published next to Brouwer's paper.  The proof in that paper was insufficient. According to Brouwer, the paper contained no proof at all. Brouwer complained to Blumenthal, and wrote to Lebesgue, who remained evasive, and never provided a proof, until the year 1921, when he published a proof of the paving principle, which was based on Brouwer's ideas.

  We must mention now Brouwer's correspondence with Poincar\'e,\footnote{In \cite{Brouwer-correspondence}, there are two letters from Poincar\'e to Brouwer and one from Brouwer to Poincar\'e, but the content shows that there were more. The set \cite{HP-correspondance1} contains one letter from Poincar\'e and one from Brouwer.} which started in 1911. Brouwer addresses the difficulties he mentioned in his Karlsruhe lecture, and he talks about the possible singularities of Riemann's moduli space. It seems that this is the first time that the latter problem is addressed.\footnote{Although Brouwer says in his first letter to Poincar\'e that this issue can be easily resolved, he will come back to it later, and he will consider it as a real difficulty. We refer the interested reader to Vol. II of  Brouwer's Collected works.} We reproduce part of this correspondence because it is extremely interesting, and several problems with which we are concerned here are addressed there in a very precise way.  In his letter to Poincar\'e dated December 10, 1911, Brouwer writes (\cite{Brouwer-correspondence} p. 113):
 \begin{quote}\small
 My `Beweis der Invarianz des $n$-dimensionalen Gebiets'\footnote{Proof of invariance of the $n$-dimensional domain, \cite{Brouwer-gebiets}.} has been inspired by reading your ``m\'ethode de continuit\'e" in Vol. 4 of \emph{Acta Mathematica}.\footnote{Poincar\'e's paper \cite{PK1887}.} It was in the course of this reading that I had the impression that on the one hand one did not know in the  general case if the one-one and continuous correspondence between the two $6g-6+2n$-dimensional manifolds concerned, is analytic, and on the other hand, that in order to be able to apply the method of continuity, one has to start by proving the absence of singular points in the variety of modules of Riemann surfaces of genus $p$; this last demonstration, incidentally, turns out to be fairly easy. Now after having read somewhere in an article by your hand (I believe about the equation $\Delta u=e^u$ in the Journal de Liouville) that you considered your exposition of the method of continuity as perfectly rigorous and complete, I started to fear that I had poorly understood your memoirs in \emph{Acta}, and I have published my memoir `Beweis der Invarianz des $n$-dimensionalen Gebiets' without indicating there the application to the method of continuity, restricting myself to an oral communication on this subject on September 27, 1911 at the Congress of the German Mathematicians in Karlsruhe, of which communication I join the text in this letter. At the occasion of this talk, Mr. Fricke expressed to me his doubts that at the start I had formulated exactly the result of your arguments of pages 250-276 of the \emph{Acta}. Meanwhile, I continue to believe that I have interpreted you exactly.
 
 In fact, if the conditions of this statement, in which the word ``uniformly" (uniform\'ement) is the key word, are satisfied, the \emph{reduced polygons of the sequence of groups converge also uniformly to the boundary of the $(2n+6p-6)$-dimensional cube}, and because of your arguments there exists at least a \emph{reduced limit polygon} that only has parabolic angles on the fundamental circle, corresponding for that reason to a limit Riemann surface, for which either the genus is decreased, or the singular points have become coincident.
 
 Would I ask you too much of your benevolence and your precious time, asking you to be so kind as to convey briefly to me your opinion about the disputed points [...]"
 \end{quote}

 On December 10, 1911, Poincar\'e responds (\cite{Brouwer-correspondence} p. 114):
  \begin{quote}\small
 [...] I do not see why you doubt that the correspondence between the two manifolds would be analytic; the moduli of Riemann surfaces can be expressed analytically as functions of the constants of Fuchsian groups. It is true that certain variables only can have real values, but the functions of those real values nonetheless preserve their analytic character.
 
Now in your eyes the difficulty arises from the fact that one of these manifolds does not depend on the constants of the group but does depend on the invariants. If I remember correctly, I considered a manifold depending on the constants of the fundamental substitutions of the group;\footnote{The term ``substitutions" denotes the fractional linear transformations.} so to a group there will correspond a discrete infinity of points in this manifold; next I subdivided this manifold into partial manifolds in such a fashion that to a group [element] corresponds a single point of each partial manifold (in the same way as one decomposes the plane into parallelograms of the periods, or the fundamental circle into Fuchsian polygons). The analytic correspondence does not seem to be altered for me.
\\ 
 With regard to the manifold of Riemann surfaces one can get into problems if one considers those surfaces in the way Riemann did: one may for example wonder if the totality of these surfaces form \emph{two} separate manifolds. This difficulty vanishes if one views   these surfaces \emph{from the point of view of Mr. Klein}: the continuity, the absence of singularities, the possibility to go from one surface to another in a continuous way become then almost intuitive truths.

  I apologize for the disjointed fashion and the disorder of my explanations, I have no hope they are satisfying to you, because I have presented them very poorly to you; but I think they will lead you to make the points that bother you more precise, so I can subsequently give you complete satisfaction. I am happy to have this opportunity to be in contact with a man of your merit.
      \end{quote}

A few remarks are due, and some consequences can be drawn from this letter:
\begin{enumerate}
\item In this letter, Poincar\'e uses the usual French word for  ``manifold," whereas in the papers quoted, he was using the terms ``multiplicity" which was motivated by the German word \emph{Mannigfaltigkeit} which Riemann had introduced and which did not have at that time a precise mathematical meaning.
\item  The letter shows that Poincar\'e has a \emph{paving} of moduli space. This cellular decomposition appears clearly to Poincar\'e, thirty years after he wrote the memoirs concerned, and this does not seem to have been highlighted before.

\item Poincar\'e says that if we stick to the methods of Riemann, it is not clear whether the manifold we obtain is connected but that the problem can be settled if one adopts the ``point of view of Klein," which also shows the ``absence of singularities."  The reference to the ``point of view of Klein" is also made in the next letter.
\end{enumerate}
The letter from Poincar\'e to Brouwer, sent in January\footnote{According to the editor of \cite{Brouwer-correspondence}, the exact date is not certain.} 1912,  also contains valuable information. It appears from that letter that in the meantime,
 Brouwer had sent several letters to Poincar\'e, asking him for details on what he knew, in particular on the question of the singularities of Riemann's moduli space 
  (\cite{Brouwer-correspondence} p. 120). Concerning the latter point, Poincar\'e repeats that the solution comes from adopting 
  Klein's point of view instead of Riemann's. The letter says the following:\begin{quote}\small
Thank you very much for your successive letters; I will study the matter in detail as soon as I have time. I still believe that the simplest way to prove the absence of a singular point would be to not use the Riemann surfaces in the form given by Riemann, in other words with stacked flat leaves and cuts, but in the form given by Klein; an arbitrary surface with a convenient connection and some law (with a representation that is or is not conformal) for the correspondence of the points of this surface with the imaginary\footnote{The word ``imaginary" means here ``complex."} points of the curve $f(x,y)=0$.
\\
Already many years ago I have expounded my ideas about this point during a session of the Soci\'et\'e Math\'ematique de France; but I did not publish them because Mr. Burckhardt, who was present at that session, told me then that Mr. Klein had already published them in his \emph{autographically prepared} lecture notes; may be you can avail yourself of these.
\\
It all amounts to this. Let $f(x,y)=0$ be a curve of genus $p$; to this curve I let correspond a Riemann-Klein surface $S$ and a law $L$ of correspondence between the real points of this surface and the complex points of the curve $f(x,y)=0$. Next I consider surfaces $S'$ and laws $L'$ that differ infinitesimally little from $S$ and $L$. One first must prove that there are $\infty^{6p-6}$ such surfaces $S'$ (which are not considered distinct if they can be transformed into each other by birational transformations); and then one can always pass from an arbitrary $S',L'$ to another arbitrary $S',L'$, without moving too far from $S,L$ and without passing by $S,L$.
\end{quote}

The last passage is extremely interesting. The space of equivalence classes under birational transformations   is Riemann's moduli space, and Poincar\'e says that in order to get rid of the singularities of that space, one has to consider a space of pairs $(S,L)$ of a real surface equipped with a ``law of correspondence" (he probably means a homeomorphism) between the real points of the surface $S$ and a complex curve. He then says a few words about the topology of that space of pairs. Although it is said in very concise words, the idea that one has to consider marked surfaces instead of surfaces in order to overcome the singularities of the moduli space is clearly stated. It seems that Poincar\'e attributes this idea to Klein. Unfortunately, Poincar\'e passed away a few months later and the correspondence with Brouwer ended there. It is unlikely that Teichm\"uller, who was the first to make explicit the idea of marking, was aware of this letter.

       \section{Weyl}\label{s:Weyl}

  Hermann Weyl (1885-1955) is another one of the mathematical giants who worked in several fields at the same time, both in mathematics and in theoretical physics.\footnote{The third author of this paper remembers a conversation with Pierre Cartier, in Strasbourg, where he asked him why Andr\'e Weil used the French-Parisian pronunciation ``Vey" for his name, and not the Alemannic pronunciation ``Wayl" (the family Weil was Alsacian). Cartier responded spontaneously that Weil had too much respect for Hermann Weyl and did not want himself to be confused with him.} In the proceedings of a conference celebrating Weyl's legacy \cite{Weyl-AMS}, Bott writes: ``Reverence and gratitude for the breadth and beauty of his work I believe we all share." Weyl, in his turn, had a great devotion for Riemann.\footnote{Weyl also used to express his gratefulness to Klein as a source of inspiration. His approach to mathematics is often described as a synthesis of Riemann's geometric ideas with Klein's ideas on group theory. The reference to Klein's \emph{Erlangen program} is recurrent in his works.}
 In 1919, he edited his \emph{Habilitationvortrag} \cite{Weyl1919}, with a commentary making the relation with relativity theory. His point of view on Riemann surfaces is a synthesis of the sometimes diverging points of view of Riemann, Weierstrass, Klein, Hilbert and Brouwer. In the winter semester of the academic year 1911/1912, Weyl gave at the University of G\"ottingen \cite{W1911} a course on Riemann surfaces. The written notes start with a historical survey:
 \begin{quote}\small
 Complex variables have an essentially real purport. Functions of complex variables have a geometric interpretation as maps and a physical interpretation as images of electric currents. Their development took place on the basis of such interpretations. One can also consider function theory as part of mathematical physics. The bases of this field were laid down by Cauchy, Riemann and Weierstrass.

Riemann took as a starting point mathematical physics. He borrowed from there the problems that characterize the development he made of function theory. Weierstrass followed another direction, starting with power series. They essentially reached the same point, and we must take into account the method of Cauchy-Riemann and that of Weierstrass.
 \end{quote}

  In this section, we shall mostly be concerned with Weyl's book \emph{Die Idee der Riemannschen Fl\"ache}  \cite{Weyl1913} (1913). The book is a summary of his 1911/1912 lectures \cite{W1911}. Bers writes in \cite{Bers-Hilbert}  about this book: ``The appearance of the first edition in 1913 may be thought of as concluding the classical period. The \emph{Idee} contained, among other things, the concepts, though not the present names, of cohomology and Hausdorff space, and the proof of the uniformization theorem based on an idea of Hilbert." In his Preface to the 1955 edition, Weyl writes:\footnote{We are quoting from the English translation of the 1955 version \cite{Weyl1913}.} 
  \begin{quote}\small 
    Three events had a decisive influence on the form of my book: the fundamental papers of Brouwer on topology, commencing in 1909; the recent proofs by my G\"ottingen colleague P. Koebe of the fundamental uniformization theorems; and Hilbert's establishment of the foundation on which Riemann had built his structure and which was now available for uniformization theory, the Dirichlet principle.
  \end{quote}
  In the  Preface to the 1913 edition, Weyl was more explicit about topology, because the subject was still nascent; he writes (p. v):
   \begin{quote}\small 
   Until now, such a rigorous presentation, building the foundations of the concepts and the propositions of \emph{Analysis situs} to which the theory of functions refers not on intuitive plausibility, but on exact proofs from set theory, was missing.
        \end{quote}
    
  As the title indicates, this book concerns chiefly Riemann surfaces.  The book is generally referred to as the reference for the abstract definition of a Riemann surface, that is, as an object independent of any algebraic function. On p. 17, in introducing surfaces, Weyl first defines the notion of manifold as a topological space, then he restricts to a characterization of two-dimensionality.  On p. 32, Weyl gives the definition of a Riemann surface as an atlas with analytic coordinate-changes. This explicit definition of a Riemann surface is at the same time the first abstract definition of a manifold defined by coordinate charts. Weyl notes, on p. 33, that the definition of an $n$-dimensional manifold was first given in Riemann's \emph{Habilitationsvortrag},  \emph{\"Uber die Hypothesen welche der Geometrie zu Grunde liegen} \cite{Riemann-Ueber} (1854). He recalls, on the same page,  the following: `` This formulation of the concept of a Riemann surface, first developed in intuitive form by F. Klein in his monograph \cite{Klein-Riemann} is more general than the formulation which Riemann himself used in his fundamental work on the theory of analytic functions. There can be no doubt that the full simplicity and power of Riemann's ideas become apparent only with this general formulation." Klein's ideas on the notion of abstract Riemann surface, confirming Weyl's comments, are presented by Remmert in \cite{Remmert1998} \S 1.2.

    Let us quote the following from the article on Weyl by Chevalley and Weil. The text is interesting because it discusses, besides Weyl, several mathematicians whose works are related to our subject matter, and it summarizes several important steps  (\cite{Chevalley-Weil}  p. 170):
      \begin{quote}\small
  It is as a student of Hilbert, and as an analyst, that Weyl had to tackle the subject of one of the first courses that he taught at G\"ottingen as a young Privatdozent: function theory after Riemann. As soon as the course was finished and written up, he ended up as a geometer, and the author of a book which was to exert a profound influence on the mathematical thought of his century. It might be that he only proposed to put on the agenda, using the ideas of Hilbert on the Dirichlet principle, the traditional  expositions of which the classical book by C. Neumann was a model. But it must  have appeared to him soon that in order to substitute correct and ``rigorous" (as these were called at that time) reasonings to his predecessors' constant appeals to intuition -- and in Hilbert's circle, it was not admitted that somebody cheats on that -- it was above all the topological foundations that had to be renewed. It seems that Weyl was not prepared for that by his previous works. It was possible for him, in this task, to rely on Poincar\'e's work, but he hardly talks about  that. He mentions, as having been, regarding this subject, profoundly influenced by the researches of Brouwer, which at that time were very new. But in reality, he does not make any use of them. Frequent contacts with Koebe, since the latter was fully dedicated to the uniformization of complex variables, became surely very useful to him, in particular for clearing up his own ideas. The first edition of the book is dedicated to Felix Klein who, of course, as Weyl writes in his preface, was fairly interested in a work which was so close to his youth's concerns, and because he provided the author with pieces of advise that were inspired from his intuitive temper and his profound knowledge of the work of Riemann. Even though he never met him, Klein embodied, in G\"ottingen, the Riemann tradition. Lastly, in one of his memoirs on the foundation of geometry, Hilbert had formulated a system of axioms based on the notion of neighborhood, underlining the fact that one can find there the best starting point for a ``rigorous axiomatic treatment of analysis situs." Of all these elements which are so different that were provided him by tradition and his environment, Weyl extracted a book which is profoundly original and which made history.
  \end{quote}

   \section{Torelli}\label{s:Torelli}

In 1957, one year before he gave a Bourbaki seminar talk on Teichm\"uller's work, Andr\'e Weil gave a Bourbaki seminar  \cite{Weil1957} on a result of Torelli\footnote{Ruggiero Torelli (1884-1915) had a very short life. He worked as an assistant to Severi in Parma, starting in 1904, and in Padova, starting in 1907. In 1915, during World War I, he joined the Italian army and died while he  was a soldier, from a heart attack,  at the age of 31. He was a very promising mathematician and his early death was felt as a great loss by the Italian mathematical community. For an exposition of Torelli's life and works, see Castelnuovo \cite{Castelnuovo-Torelli}  and Severi \cite{Severi-Torelli}.} on which recent progress had been made. The theorem says that we can recover a Riemann surface from its polarized Jacobian. 

 Torelli's theorem roughly says that the surface can be reconstructed from its Jacobian. One point which makes the theorem interesting is the fact that the Jacobian of the surface, although a higher-dimensional complex manifold, is in principle is simpler than the surface, due to of its linear character. As stated by Weil, Torelli's theorem says the following:

\begin{theorem}
Let $C$ and $C'$ be two algebraic curves of genus $g$ with polarized Jacobian varieties $J$ and $J'$, and let $\phi:C\to J$ and $\phi':C'\to J'$ be their canonical maps. If we have an isomorphism $\alpha:J\to J'$ (respecting the polarizations), then there exists a mapping $f:J\to J'$ which is of the form $\pm \alpha+c$ (where $c$ is a constant) which maps $\phi(C)$ to $\phi'(C')$. In particular, since $\phi$ is an isomorphism from $C$ to $\phi(C)$ and  since $\phi'$ is an isomorphism from $C'$ to $\phi'(C')$, this implies that $C$ and $C'$ are naturally isomorphic.
\end{theorem}

The proof of this theorem consists in constructing, from the knowledge that the polarized variety $J$ is the Jacobian of a curve, the curve $C$ and the map $\phi$ from $C$ to $J$. Stated in other words, the proof involves the construction of  a space of Riemann surfaces in terms of periods of holomorphic differentials. Abelian integrals appear in this theory, and we are back to the ideas of Riemann. Torelli's main article on this subject is \cite{T1913}. There is also a local Torelli theorem which amounts to the injectivity of the differential of the period map at non-hyperelliptic
curves, see  (\cite{Kodaira1964}, Theorem 17). This result is attributed to Andreotti and Weil.

The Torelli theorem and its proof are important for our subject for several reasons. First of all, the theorem concerns our problem of moduli of algebraic curves. Secondly, the theorem involves a map from an algebraic curve into its Jacobian variety, and the collection of such maps embeds the space of equivalence classes of what Weil calls ``Torelli surfaces" (cf. \S \ref{s:Weil}) into a complex vector space. This embedding can be used to show that the Torelli space has a complex-analytic structure. The result is in the direct lineage of Riemann's work on periods of Abelian integrals, and it is also a forerunner of Teichm\"uller's work on  the complex structure of Teichm\"uller space using period matrices. Thirdly, this theorem  is a rigidity theorem. It says that an isomorphism between the Jacobian varieties associated with two algebraic curves is induced by an isomorphism between the curves. Several rigidity theorems which have the same flavor were proved later in Teichm\"uller theory.
 
 In the introduction of the Bourbaki seminar paper on Torelli's theorem,  Weil writes the following, which gives a good idea of the whole setting:
 \begin{quote}\small
 It is said that an algebraic curve of genus $g>1$ depends on $3g-3$ ``moduli" (whatever meaning one should give to this sentence!).\footnote{ This is reminiscent of a sentence in Teichm\"uller's paper \cite{T32} which we already quoted, cf.  Footnote (\ref{F:Teich}).} But its integrals of the first kind, normalized as it is usually done in the classical theory using a system of ``retrosections," admit a period matrix which is symmetric and contains therefore $g(g+1)/2$ independent coefficients. Nothing prevents us from conjecturing that there are $\frac{g(g+1)}{2}-3g-3=\frac{(g-2)(g-3)}{2}$ relations among them; and no mistake was found in this, since in the hundred years of existence of the theory.
 
 The goal of the Torelli theorem is to show the first assertion after giving it a precise meaning. 
 \end{quote}

Torelli's theorem was generalized and extended in several ways and there is extensive current research in algebraic geometry motivated by it. Moreover, Torelli's original proof is characteristic free and it lead directly to extensions to curves on other fields. The \emph{Schottky problem} asks for a characterization of Jacobians of curves among all principally polarized Abelian varieties. This is another subject of current research which is in the lineage of the work of Torelli. The last article \cite{T1915} written by Torelli is devoted to that problem.
  
  \section{Siegel}\label{s:Siegel}
 
After Torelli, it is natural to talk about Siegel.\footnote{Carl Ludwig Siegel (1896-1981) and Weil were friends. This friendship started in G\"ottingen and continued in Princeton, cf. Weil's \emph{Collected Papers} \cite{Weil1979} Vol. I p. 521 and 557. In Vol. II (p. 544), Weil writes that ``to comment on Siegel has always appeared to me as one of the most useful tasks that a mathematician of our epoch can undertake." Weil and Siegel also worked on common subjects; for instance, one may recall the famous \emph{Siegel-Weil formula}, and the work of Weil contains various generalizations of  results of Siegel on quadratic forms \cite{1962a} \cite{1964b}.
Krantz recounts in \cite{Krantz} (p. 186) that when Weil was asked who was the greatest mathematician of the twentieth century, the answer, without hesitation, was ``Carl Ludwig Siegel."  In this section, we talk about Siegel's work on the Torelli space, but there are several other contributions to topics which are related to the modern Teichm\"uller theory; for instance, the computation of volumes of fundamental domains of discrete groups, but we cannot expand more on this in the present paper.} He introduced techniques for understanding moduli that are in the tradition of Riemann and of Torelli. In particular, he worked extensively on the question of reconstructing the equations of an algebraic curve from the associated periods of meromorphic Abelian differentials. Siegel went a step further than his predecessors, and he managed to avoid the complications of the method of continuity. In   comments that Weil made on an Air Force final report, as an addendum to his \emph{Collected Papers} edition (\cite{Weil1979}, Vol. II p. 545), he writes:
\begin{quote} \small
The theory of moduli of curves, inaugurated by Riemann, has taken
two decisive steps in our time, first in 1935, with the work of Siegel
(\emph{Ges. Abh. No. 20}, \S 13, Vol. I, pp. 394-405), and then with
the remarkable works of Teichm\"uller; it is true that some concerns were raised about the
latter, but these were finally cleared up by Ahlfors in 1953 (\emph{J. d'An. Math.} 3, pp.
1-58). On the other hand, we realized at the end that Siegel's discovery of automorphic functions belonging to the symplectic group applied first to the moduli of Abelian varieties and then, only by extension, to those of curves, via their Jacobians and using Torelli's theorem. Thus, thanks to Siegel, we have at our disposal the first example of a theory of moduli for varieties of dimension $>1$. We owe to Kodaira and Spencer (\emph{Ann. of Math}. 67 (1958), pp. 328--566) the discovery of the fact that progress on cohomology allows us not only to address a new aspect of the same problem but also, at least from the local point of view, to tackle the general case of complex structures on varieties.
\end{quote}
In the article \cite{Siegel20}  which Weil quotes, Siegel generalizes the classical study of quadratic forms in two variables to the study of forms in $n$ variables and he develops a series of results in class field theory. At the end of this article, Riemann surfaces appear and Siegel gives a construction of the Torelli space. He generalizes the half-plane representing the hyperbolic plane (which was used for binary forms) to a higher-dimensional half-space which is now called the Siegel upper half-space. This is a space of symmetric square matrices whose imaginary part is positive definite. Siegel describes in the Siegel half-space a fundamental domain for the action of the symplectic group (which is in some sense a higher-dimensional analogue of the modular group acting on the two-dimensional upper half-space). He notes that this action is discrete and that this provides a complex structure on a generic part of moduli space.
This result is described in the form of a map  from the Torelli space into the Siegel upper half-space, equipped with a complex structure which is K\"ahler. We shall dwell below on the importance of K\"ahler structures in the  theory of moduli. Combining this with the canonical map from Teichm\"uller space to the Torelli space, one obtains a map from Teichm\"uller space $\Theta$ into the Siegel space. This map is described by Weil in his Bourbaki seminar on Teichm\"uller space (\cite{Weil1958a}, p. 383).  Weil declares that this map is holomorphic when $\Theta$ is provided with its natural complex structure. The image is an analytic subvariety of the Siegel space, whose points are all smooth except those corresponding to hyperelliptic Riemann surfaces.  Weil then says: ``As for $\mathcal{M}$ (Riemann's moduli space) there is virtually no doubt that it can be provided with the structure of an algebraic variety (non-complete of course, and with multiple points), the ``variety of moduli," so the natural mapping of $\Theta$ onto $\mathcal{M}$ is holomorphic."

\section{Teichm\"uller} \label{s:Teich}

In 1958, Weil gave a Bourbaki seminar \cite{Weil1958c} on  the results of Teichm\"uller. The aim of these seminars is to present recent and important advances on a specific topic. In the present case, Teichm\"uller's results were already twenty years old. The fact that these results were presented twenty years after their discovery is meaningful.\footnote{
 Oswald Teichm\"uller (1913-1943) lived and died in the darkest period of German history. Little is known about his life, compared to the ones of other mathematicians of his stature. The reason is that he was a notorious Nazi, and writing about his life has been a delicate matter for historians. Dieudonn\'e says, in a review he wrote on the paper \cite{Hauser} (MR1152479): ``Oswald Teichm\"uller's life is a tragic example of one of the most brilliant minds in his generation of mathematicians, who fell prey to a fanatical doctrine that was bent upon stifling all feelings of decency and compassion." A biographical sketch, based on notes written by his mother, is contained in English translation in Abikoff's article \cite{Abikoff1986}. Teichm\"uller entered the University of  G\"ottingen at 17, with very broad interests. Two years later, in 1932, he joined the Nazi Party. He obtained his doctorate in 1935 under Helmut Hasse. The subject of his thesis was linear operators on Hilbert spaces over the quaternions. Fenchel was one of his teachers, and in a private conversation with Abikoff (which the latter shared with the authors of this article), he described him as a lonely boy who grew up in the Harz mountains. It seems that Teichm\"uler found friends in the Party. In 1935-36, he followed lectures on function theory by Egon Ullrich and Rolf Nevanlinna and became interested in value distribution problems. Nevanlinna was teaching at G\"ottingen, as a visiting professor. Teichm\"uller moved in 1937 to Berlin, where he joined the Bieberbach group and where he earned a qualification to lecture. The Bieberbach conjecture became one of his strong centers of interest. In 1939, he obtained a regular position at the University of Berlin, but the same year he was drafted in the army, and he remained a soldier until his death. He died at the Eastern Front, in the Dnieper region, in 1943. He wrote several of his papers on moduli while he was a soldier. As a mathematician, Teichm\"uller was exceptionally talented. In the review we mentioned, Dieudonn\'e writes: ``the diversity of mathematical problems which he tackled in his 34 published papers, during barely 10 years, is amazing for such a young man: from logic, through algebra, number theory and function theory to topology and differential geometry, most of these papers display an originality that later research entirely justified." Five papers, on various subjects, were completed during the last year of his life, while he was on the front, and they were published in 1944, the year after his death (four of them in  \emph{Deutsche Mathematik}  and one in \emph{Crelle's Journal}). Concerning Teichm\"uller's deeds, Bers quotes in \cite{B1960} a famous sentence by Plutarch (\emph{Perikles} 2.2) which is often repeated: ``It does not of necessity follow that, if the work delights you with its grace, the one who wrought it is worthy of your esteem." We quote Gustave Flaubert: ``L'homme n'est rien, l'\oe uvre est tout." (Letter to George Sand, December 1875).} In the first page of \cite{Weil1958c}, Weil writes: 
\begin{quote}\small
Teichm\"uller proved, first heuristically, then (it is said) rigorously, that we can define a topological space, homeomorphic to an open ball of dimension $6g-6$, whose points are ``naturally" in one-to-one correspondence with the classes of Teichm\"uller surfaces (classes with respect to the equivalence relation defined by isomorphism); Ahlfors gave another proof of the same fact.
\end{quote}
 In the same year, in a report written for an Air Force  contract (\cite{Weil1958b}  p. 91), Weil writes:
\begin{quote} \small
[...] Teichm\"uller's chief contribution was to define on $T$ [Teichm\"uller space] a certain topology, the ``natural" one in a sense described below, and then to prove that $T$, with this topology, is homeomorphic to an open cell of real dimension $6g-6$. So far, I have mainly been concerned with the local properties of Teichm\"uller space and of its ``natural" complex analytic structure. The definition of the latter depends upon ideas introduced by Teichm\"uller himself, but which do not appear to have been fully understood until Kodaira and Spencer attacked similar problems for higher dimensions.
\end{quote}

Teichm\"uller had a new approach to the moduli problem. He defined Teichm\"uller space as a cover of moduli space in which the singularities disappear, and he thoroughly studied its metric and complex structures. This point of view was new as compared to those of Poincar\'e, Klein and others who did not study the moduli space of Riemann surfaces as an intrinsic object, but as a tool for understanding uniformization. We start with a list of the major contributions of Teichm\"uller to the moduli problem.
\begin{enumerate}

\item \label{t1} The formulation of Riemann's moduli problem as a problem of existence of a complex structure on the space of equivalence class of Riemann surfaces and the interpretation of the number of moduli as a complex dimension. 

\item \label{t101} The solution of the moduli problem   described in (\ref{t1}).

\item \label{t102} The introduction of marking as a homeomorphism defined up to homotopy from a fixed topological surface to a varying Riemann surface.

\item  \label{t2} The introduction of quasiconformal mappings as an essential tool in the study of moduli of Riemann surfaces.

\item \label{t23} The introduction of the Teichm\"uller metric. Teichm\"uller showed that this metric is Finsler and he studied its infinitesimal norm and its geodesics. 

\item \label{t25} The  infinitesimal  theory of quasiconformal mappings, and a thorough study of the partial differential equations of which they are the solution.

\item \label{t24} The introduction and the study of Teichm\"uller discs (which Teichm\"uller called ``complex geodesics").

\item \label{t26} The comparison between the hyperbolic length of closed geodesics and the quasiconformal dilatation of a map.

\item \label{t27} A correct use of the continuity method in the moduli problem, on which Poincar\'e, Klein and others had wrestled with for a long period of time. For the first time, the method was applied to objects which were known to be manifolds of the same dimension.  

\item \label{t28} The existence and uniqueness of the universal Teichm\"uller curve and its use in the existence of the complex structure on Teichm\"uller space. At the same time, this introduced the first fibre bundle over Teichm\"uller space,  \emph{Teichm\"uller universal curve}. 

\item \label{t29} The proof that the automorphism group of the universal Teichm\"uller curve is the extended mapping class group.

\item \label{t212}The idea that Teichm\"uller space is defined by a universal property. Grothendieck expressed this later on,  by the fact that Teichm\"uller space represents a functor.

\item \label{t213} The idea of using the period map to define a complex structure on Teichm\"uller space.  

\item \label{t215} The question of whether there is a Hermitian metric on Teichm\"uller space. 
\end{enumerate}

 There has been a lot of confusion about what Teichm\"uller proved, and this is due to several reasons. The first reason is probably non-mathematical, viz., his declared anti-semitism. Because of that, many mathematicians were reluctant to the idea of reading his papers. Another reason is that part of his results on moduli are expressed in the language of algebraic geometry. Since the papers were read essentially by analysts, several points were missed. For example, Ahlfors' review (MR0018762) of Teichm\"uller's 1944 paper \cite{T32} shows that he completely missed all the important results in that paper, probably for lack of knowledge or appreciation of algebraic geometry, and may be also of topology. In fact, it took Ahlfors several years to realize the importance of that paper. The definition of the mapping class group which is given in that review is even wrong.  The confusion that surrounds Teichm\"uller's paper is also due in part to his style. In fact it is difficult for someone who does not go into his papers thoroughly to figure out what he states as a conjecture that he cannot prove, what he first states as a conjecture and then proves later in the same paper, or what he postpones for a subsequent paper. We shall elaborate on all this in the text that follows. Let us first make some quick comments on the above list of results.
   
  The results (\ref{t102})
(\ref{t2}),
(\ref{t23}),
(\ref{t25}),
(\ref{t24}),
(\ref{t26}),
(\ref{t27}), are contained  in the paper \cite{T20} which we shall summarize below. The results (\ref{t1}), (\ref{t101}),  (\ref{t28}), (\ref{t29}), (\ref{t212}), (\ref{t213}) are contained in the paper \cite{T32} on which we shall also comment.
   
A consequence of (\ref{t25}), which also involves the equivalence relation on Beltrami differentials and the complex structure on each tangent space of Teichm\"uller space, is  (using modern terms) the existence of a natural almost-complex structure  on Teichm\"uller space. This structure was shown, later on, using the deformation theory of Kodaira-Spencer or the Newlander-Nirenberg theorem, to be integrable. This approach was at the basis of the complex theory of Teichm\"uler space developed by Weil, Ahlfors and Bers. 

In fact, it seems that the first appearance of the expression ``complex analytic manifold" appears in Teichm\"uller's paper \cite{T32} (under the German name \emph{komplexe analytische Fl\"ache}). Let us quote Reinhold Remmert, one of the founders of  the theory of several complex variables, talking about the history of the theory of complex manifolds of higher dimension (\cite{Remmert1998}, p. 225): 
\begin{quote}\small
It seems difficult to locate the first paper where complex manifolds explicitly occur. In 1944 they appear in Teichm\"uller's work on ``Ver\"anderliche Riemannsche Fl\"achen" \cite{T32} \emph{Collected Papers} p. 714);\footnote{The quote concerns the complex manifolds which are not defined as domains of some complex vector space. (Domains of $\mathbb{C}^N$ had been studied earlier by several authors, in particular, C. L. Siegel, in the context of modular forms or automorphic functions of several variables.) Teichm\"uller space was not known yet to be a complex domain. An embedding of that space in some $\mathbb{C}^N$ was found later; cf. \S \ref{s:AB}.}
 here we find for the first time the German expression ``komplex analytische Mannigfaltigkeit." The English ``complex manifold" occurs in Chern's work (\cite{Chern1946}, p. 103); he recalls the definition (by an atlas) just in passing. And in 1947 we find ``vari\'et\'e analytique complexe" in the title of \cite{Weil1947b}. Overnight complex manifolds blossomed everywhere.
\end{quote}

Concerning Item (\ref{t28}): In his book \emph{Lectures on quasiconformal mappings} \cite{Ahlfors-lectures}, Ahlfors constructs (Section D of Chapter 4) the Teichm\"uller curve, with a reference to Kodaira-Spencer and without mentioning that this object was already defined by Teichm\"uller in his paper \cite{T32}. The second edition of this book contains an update, in the form of a very rich survey by Earle and Kra, titled \emph{A supplement to Ahlfors lectures}, where references to the Teichm\"uller curve are given (in \S 2.7), namely to the papers by Grothendieck \cite{GGG}, Engber \cite{Engber1975} and Earle-Fowler \cite{EF}, but without any mention of Teichm\"uller's paper. The reference is also missing in the three papers quoted.
 
  Concerning Item (\ref{t27}): Teichm\"uller applied the method of continuity in his main paper \cite{T20} and also in \cite{T24}, and he also applied the invariance of domain in his paper \cite{T29} where he proves the general existence result for quasiconformal mappings. In this paper, the surfaces are represented by polygons in the hyperbolic plane, in the tradition of Poincar\'e.

 Concerning Item (\ref{t212}), a few words on the origin of quasiconformal mappings are in order. Ahlfors gives a brief summary of the early use of these mappings  in \cite{A1964}, p. 153. According to this account, these mappings were introduced (by a different name) in the 1928 paper \cite{Greotzsch1928} by Gr\"otzsch. Ahlfors notes that Gr\"otzsch's paper was ``buried in a small journal," that it first remained unnoticed, and that he learned of it by word of mouth in 1931. He adds  that the problem of finding the best mapping between a rectangle and a square was first considered as a mere curiosity, and that the full strength of quasiconformal mappings and their use in the deformation theory of Riemann surfaces was first realized by Teichm\"uller. The expression ``quasiconformal mapping" appears in print in Ahlfors' 1935  fundamental paper on covering surfaces \cite{Ahlfors1935}.\footnote{This is the paper for which Ahlfors earned the Fields medal.}
In the comments he makes in his \emph{Collected Papers} edition (\cite{Ahlfors-collected} Vol. 1, p. 213), Ahlfors recalls that Gr\"otzch  called these mappings ``nichtconforme Abbildungen," and that almost simultaneously with his own 1935 paper, Lavrentiev, in the Soviet Union, published a paper (in French) in which he introduced a notion which is equivalent to quasiconformal mappings, which he termed ``fonctions presque analytiques."  Teichm\"uller in his paper \cite{T20} mentions Gr\"otzsch, Ahlfors and Lavrentiev, and he refers several times to what he calls the ``Gr\"otzsch-Ahlfors method."

Teichm\"uller had the reputation of not writing complete proofs. In fact, it seems that it was considered as a German tradition to write only the main ideas and not the technical details. According to Abikoff, this was the style of the journal in which most of Teichm\"uller's articles were published. He writes in \cite{Abikoff1986}: ``The tradition of \emph{Deutsche Mathematik} is one of heuristic argument and contempt for formal proof. Busemann notes that Teichm\"uller manifested those traits early in his career but when pressed could offer a formal proof." One should also mention that this journal was not easily accessible.\footnote{Abikoff says in \cite{Abikoff1986} that he started his search for some of  Teichm\"uller's papers during a stay in Paris at the end of the 1970s, and that it was very difficult for him to find the ones published in \emph{Deutsche Mathematik}.} Today, his \emph{Collected papers} are easily accessible, but there is still a tendency not to read them. 
 
   A recurrent problem in Teichm\"uller's writing is his use of the word ``conjecture." In fact, he uses it in two different senses: either as a statement which he does not prove at that precise location in the text, but which he  proves later, or as a statement of which he does not give a complete proof, even if he gives a sketch of a proof which is sufficient for an astute reader to work out a complete proof (this happens several times in \cite{T20}). Ahlfors writes in \cite{Ahlfors1954}: ``It requires considerable effort to extricate Teichm\"uller's complete and incontestable proofs from the maze of conjectures in which they are hidden." In fact, most of the times the word ``conjecture" occurs  in the paper \cite{T20}, it concerns the relation between quadratic differentials and extremal quasiconformal mappings, which is  one of the main results in that paper. Teichm\"uller gives a lot of evidence for it, but the complete proof is given in the later paper \cite{T29}.
Another instance occurs  in \S 19 of the same paper, where Teichm\"uller writes, after the definition of his metric: ``We will later be led to the reasonable hypothesis that $\mathfrak{R}^\sigma$ is with regard to our metric a \emph{Finslerian space}." Here, $\mathfrak{R}^\sigma$ is Riemann's moduli space and this also looks like a conjecture, but in fact, it is proved later in the same article, for Teichm\"uller space, which is introduced later on, and the Finsler structure is studied in detail. 
   
  As an example of the confusion in the literature that surrounds Teichm\"uller's statements, we  quote Abikoff from \cite{Abikoff1986} (p. 16), talking about the paper \cite{T20}: ``He shows that he is aware that he has no proof of the key existence theorem. (He offered a proof of this theorem in \cite{T29} but it was never really accepted by the mathematical community although, as Bers notes, Teichm\"uller's proof is correct and complete.)"
 The reader wonders why the theorem was never really accepted since Bers noted that the poof is complete and correct.

As an additional sign of that confusion, we note that the definition of Teichm\"uller space and of the action of the mapping class group on it was very poorly understood, even among eminent geometers. We  mention as examples Narasimhan's Preface to Riemann's collected works \cite{Riemann-Gesammelte} and Ahlfors' review mentioned above. They both contain a wrong definition of the mapping class goup. Dieudonn\'e, in his book \cite{Dieudonne-cours}, in the section on Teichm\"uller's work  (\S IC.12, Vol. I p. 218) also gives an incorrect definition of the mapping class group. A mistake also occurs in Bers' paper \cite{Bers-simultaneous}. All this shows that several eminent analysts and algebraic geometers have some trouble with basic topological notions. Much more recently, F. Kirwan, in her article ``Moduli spaces" published in the \emph{Encyclopedia of mathematical physics} \cite{Enc} has a wrong definition of Teichm\"uller space and a very naive and also wrong description of the complex structure on that space. She attributes the fact that Teichm\"uller space is homeomorphic to a ball to Bers. At the beginning of her article, the author makes a short history of the subject, and her view on this history probably represents that of the majority of the algebraic geometers. In this short history, she passes directly from Riemann to Mumford, with no mention of Teichm\"uller and of Grothendieck, who led the foundations of the analytic and the algebraic structure of moduli space. This article, published in 2006, is a revised version of a survey article she wrote in 1998 and it seems that no algebraic geometer caught these mistakes before the revision.\footnote{Narasimhan, in his Editor's Preface of the new reprint of Riemann's Collected Works edition \cite{Riemann-Gesammelte} (1990), writes (p. 15): ``In the 1930s, O. Teichm\"uller introduced another idea into the study of questions of moduli. He considered \emph{pairs} consisting of a compact Riemann surface of genus $g$ and a fixed choice of homology basis, and showed the equivalence classes of such pairs for a space homeomorphic to $\mathbb{R}^{6g-6}$ ($g>1$)." The space he describes is the Torelli space, instead of Teichm\"uller space. Moreover, Narasimhan talks about the complex structure on Teichm\"uller space without  mentioning Teichm\"uller's 1944 paper \cite{T32}, which is the basic paper on the subject. Kirwan writes in \cite{Enc} p. 453: ``We consider the space of all pairs consisting of a compact Riemann surface
of genus $g$ and a basis  $\gamma_1, ..., \gamma_{2g}$ for $H_1(\Sigma, \mathbb{Z})$ as above such that
$\gamma_i \gamma_{i+g}=1=-\gamma_{i+g}\gamma_i$ 
if $1\leq i\leq g$ and all other intersection pairings $\gamma_i \gamma_j$ are zero.
If $g\geq 2$, this space (called Teichm\"uller space) is naturally homeomorphic to an open ball in $\mathbb{C}^{3g-3}$ (by a theorem of Bers)." This space she describes is again the Torelli space and not Teichm\"uller space. Ahlfors writes in his review: ``Let $H_0$ be a fixed surface, $H$ a variable surface and $T$ a topological map of $H_0$ onto $H$. Then the couple $(H,T)$ is an element of the space $\mathfrak{R}$ if and only if $H$ can be mapped onto $H_0$ by a conformal transformation $S$ such that $TS^{-1}$ can be deformed to the identity." Dieudonn\'e writes that the mapping class group is ``the discrete group $\gamma$ of isomorphism classes of Riemann surfaces of genus $g$, where two isomorphisms belong to the same class if they fix a given point of Teichm\"uller space."}
 
In the rest of this section, we review some papers by Teichm\"uller on moduli, for the sake of making  clearer what results he obtained.
  
  We start with the paper \emph{Extremale quasikonforme Abbildungen und quadra\-tische Differentiale} (Extremal quasiconformal maps and quadratic differentials)  \cite{T20}. The paper is translated in \cite{Theret}, with an extended commentary in \cite{T20C}.
 It contains 170 sections divided into 31 chapters. This theory is developed for the most general surfaces of finite type: orientable or not, with or without boundary, with or without distinguished points which may be in the interior or on the boundary. The author declares (\S 12): ``in later works, we shall study  for example the case where two distinguished points are moved infinitely close together. In the present work, this will be excluded." This remark in itself is a whole program on degeneration of Riemann surfaces. The paper contains several such ideas.
  
Surfaces equipped with
  complex structures and distinguished points are called
   \emph{principal regions}.\footnote{Hauptbereich.} Teichm\"uller denotes (\S 13) the space of conformally equivalent principal regions by $\mathfrak{R}^\sigma$. This is Riemann's moduli space. The author says that the conformal invariants ``are precisely the functions on $\mathfrak{R}^\sigma$." The exact value of
    $\sigma$ is given in \S 14 as:
  \[\sigma -\rho = -6 +6g +3\gamma +2n +3h +k
 \]
 where

 $g$ = number of handles;
  
 $\gamma$ = the number of crosscaps;
 
 $n$ = the number of distinguished interior points;
 
 $h$ = the number of distinguished boundary components;

 $k$ = the number of distinguished points on the boundary;
 
 $\rho$ = the number of parameters of the continuous group of conformal mappings of the principal region onto itself.
 
Teichm\"uller states that ``on a small scale, $\mathfrak{R}^\sigma$ is homeomorphic to the $\sigma$-dimensional Euclidean space." Strictly speaking, this is not correct, because of the existence of singular points. In \S 49, he will consider Teichm\"uller space, which has no singular points. 
  
   
 In \S 15, Teichm\"uller introduces the \emph{dilatation quotient} of a quasiconformal mapping, and in \S 16, he states the following general problem, which concerns the behavior of conformal invariants under quasiconformal mappings:
 \begin{problem}\label{p:conf}
 Let a conformal invariant $J$, seen as a function on $\mathfrak{R}^\sigma$ for a fixed principal region and a number $C$ be given. What values does $J$ assume for those principal  regions onto which the given principal region can be quasiconformally mapped so that the dilatation quotient is everywhere $\leq C$. 
 \end{problem} 
 This is also a broad problem with many ramifications.

In \S 18, the Teichm\"uller distance is introduced as the logarithm of the dilatation quotient. In 
\S 19, the author says that this distance is Finsler. He then gives is a series of examples before stating a theorem which is a wide generalization of a result of Gr\"otzsch on the conformal representation of quadrilaterals.

 In \S 34, Teichm\"uller starts the study of surfaces whose universal cover is the hyperbolic plane. In \S 35, he addresses the question of conjugating a conformal structure by a quasiconformal homeomorphism. He proves an inequality between quasiconformal dilatation and the dilatations of the isometries of the hyperbolic plane. This is an early version of an inequality, rediscovered later on by Sorvali \cite{Sor} and Wolpert \cite{Wolpert79}. He then asks a question concerning the topology defined by the sequence of inequalities associated with the elements of the isometry group acting on the upper half space. This is analyzed in \cite{AT}. 
 
 In \S 36, the study of extremal quasiconformal mappings starts. The author states that  ``this should naturally not be a proof, but rather a heuristic consideration: the maximum dilatation must be constant $=K$ everywhere." In \S 37, he defines the measured foliations on the surface obtained through the consideration of the infinitesimal ellipses that are images of infinitesimal circles by quasiconformal maps.\footnote{Teichm\"uller mentions in a note that ``the idea seams to have been expressed by Lawrentieff."} 
 He addresses the question of which foliation is associated with a given extremal quasiconformal mappings. 
He then defines (\S 38) the notion of a \emph{locally extremal} quasiconformal mapping. 
 \S 41 to  \S 45 concern the Riemann--Roch theorem. In \S 46, he writes:
\begin{quote}\small
\emph{We now conjecture that there is a connection between everywhere finite quadratic differentials and extremal quasiconformal mappings}
\end{quote} 
and he adds the following, which gives an answer to te latter question addressed in \S 37.
\begin{quote}\small
\emph{Here, I arrived in a night in 1938 at the following conjecture:} Let $d\zeta^2$ be an everywhere finite quadratic differential on $\mathcal{F}$ different from 0. One assigns to every point of $\mathcal{F}$ a direction where $d\zeta^2$ is positive. All the extremal quasiconformal mappings are described through the direction fields thus obtained and arbitrary constant dilatation quotient $K\geq 1$.
\end{quote}
He then spends several sections on testing this conjecture on special cases. 

Chapter 15 starts  at \S 49. It is titled
\emph{The topologically determined principal regions}. The mapping clas group, although without a name, appears for the first time in this section. The expression ``topologically determined" means marked. The author introduces the Teichm\"uller space $R^\sigma$ with its topology, equipping the moduli space with the quotient topology. He explains how a marking by homeomorphisms determines homotopy classes of curves. He states as a conjecture that $R^\sigma$ is homeomorphic to a $\sigma$-dimensional Euclidean space. 
  
 Chapter 16 starts at \S 53 and is titled \emph{Definition of principal regions through metrics}. From here on,  the conformal structures are given as structures underlying Riemannian metrics, and they are written, in the tradition of Gauss, as:
\[Edx^2+2Fdxdy+Gdy^2.\]
 This infinitesimal structure leads Teichm\"uller to the approach to conformality and to the moduli problem using partial differential equations.
 In \S 55 and \S 56, he explains how the main problem on extremal quasiconformal mappings is now replaced by a problem between Riemannian metrics. \S 57ff. concern the determination of the quasiconformal dilatation in terms of the coefficients of the Riemannian metric on the surface.  
In Chapter 17, Teichm\"uller starts a classification of infinitesimal quasiconformal mappings and a  thorough development of the theory of the Beltrami equation.  In \S 85, he introduces in this study Poincar\'e series and automorphic functions. 

Chapter 20 starts at \S 88 and is titled \emph{The linear metric space $L^\sigma$ of the classes of infinitesimal quasiconformal mappings}. There, Teichm\"uller describes infinitesimal quasiconformal mappings through quadratic differentials. At the end of \S 91, he asks whether one can define a Hermitian product on the tangent space.  

Chapter 21, called \emph{The duplication}, starts at \S 92. It concerns the operation  of  \emph{doubling} of a Riemann surface with boundary, and also the operation of passing to the orientation-cover for non-orientable surfaces.

In \S\S 108 and 109, he obtains the dimension of the tangent space to Teichm\"uller space. He considers periods of differentials and the period map. He then gives a formula  for the dimension of Teichm\"uller space using the theorem of Riemann--Roch. In \S 112, he announces another proof using the ray structure and the geodesics. 

Chapter 23, which starts at \S 113, is titled \emph{Going to finite mappings: $R^\sigma$ as a Finsler space}. Teichm\"uller returns to the finite (versus the infinitesimal) quasiconformal mappings. He considers that the problem of the infinitesimal extremal quasiconformal mappings is solved. He gives a formula for the finite extremal quasiconformal mappings and he explains how a passage to the limit gives the infinitesimal quasiconformal extremal mapping. He declares: ``The supposition is based on nothing." In \S 114, he reformulates more precisely his main conjecture, in terms of the metric, and he then moves on to the study the Finsler structure and the existence of geodesics connecting arbitrary points. The extremal quasiconformal mappings arise now from the extremal infinitesimal quasiconformal mappings, as one keeps the direction fields and gives to the dilatation quotient a finite constant value $>1$. In  \S 116 and  \S 117, he proves that the Teichm\"uller metric is Finsler. In \S 120, he studies the backward extension of geodesic rays and he proves that the two geodesic rays emanating from the same point corresponding to $d\zeta^2$ and $-d\zeta^2$ join together into a geodesic line. 

In  \S 121, given a point $P$ in $R^\sigma$, a complex number $\varphi$  of modulus 1, and a real number $K\geq 1$, Teichm\"uller considers the point $P(K,\varphi)$ of $R^\sigma$ whose direction is $e^{-i\varphi}d\zeta^2$ and at distance $\log K$ from $P$.
 The target principal region is denoted by $\mathfrak{h}(K,\varphi)$ and is described by the metric
 \[
 \vert d\zeta^2\vert + \frac{K^2-1}{K^2+1}\mathfrak{R}e^{-i\varphi}d\zeta^2.\]
 He considers $\log K$ and $\varphi$ as ``geodesic polar coordinates" and he calls the two-dimensional space constituted by all these Riemann surfaces $P(K,\varphi)$, for variable $K$ and $\phi$, a \emph{complex geodesic}. He states that in the special case where the surface $\mathfrak{h}$ is a topologically fixed ring with no distinguished point (\S 26 to \S 28), the associated space $R^2$ is a ``unique complex geodesic."  Equipped with the metric under study, it is the hyperbolic plane (constant curvature $-1$). More generally, he establishes the following:
 \begin{proposition}
 Any complex geodesic is isometric to the hyperbolic plane.
 \end{proposition}
In  \S 122, he proves the following:
 \begin{proposition}
 Between any two points of $R^\sigma$ passes a unique geodesic, and this geodesic is distance minimizing. 
 \end{proposition}
In \S 123, for any point $P$ in $R^\sigma$, using the geodesic polar coordinates defined in \S 121, he maps  the space $R^\sigma$  homeomorphically onto the Euclidean space $\mathbb{R}^\sigma$, which gives the following:
 \begin{proposition}
 The space $R^\sigma$ of classes of conformally equivalent topologically fixed principal regions of a given type is homeomorphic to a $\sigma$-dimensional Euclidean space.
  \end{proposition}
In    \S 124ff., Teichm\"uller considers special cases of this theorem (about 24 cases). In particular, in \S 129 he considers the case of pentagons, for which he can find explicit parameters, and then the case of hexagons. The space of hexagons is three-dimensional. He studies geodesic lines in this space, and he discusses their behavior at infinity. An application of the Schwarz reflection principle reduces the case of a ring domain to that of the torus (\S 131). In \S 141, he proves: 
\begin{proposition}
A conformal mapping of a closed orientable Riemann surface onto itself   which is homotopic to the identity can be deformed to the identity through conformal maps.
\end{proposition}
In particular, for genus $>1$, since there are no holomorphic continuous deformations, the map itself is the identity. The proof uses the action on quadratic differentials.

In \S 142, he considers the (extended) mapping class group $\mathfrak{F}$. It acts on Teichm\"uller space, and it keeps the metric invariant. He then studies the special cases of the torus and of the closed orientable surface of genus 2.

In  \S 147 he starts the study of convexity in Teichm\"uller space. He makes the following definition:  \emph{Under a geodesic manifold in $R^\sigma$, we understand a nonempty subset of $R^\sigma$ which always contains for two distinct points the whole geodesic running through both points.}
He relates this study to the question of finite subgroup actions on Teichm\"uller space, and he asks whether any finite subgroup of the mapping class group fixes a Riemann surface.\footnote{This is Nielsen's realization problem \cite{Nielsen}.}
In \S 150, he examines an example of a geodesic manifold of dimension $4g-2$. He then addresses the question of what is a geodesic subspace with respect to complex geodesics.

In \S 158, Teichm\"uller studies periods of tori embedded in 3-space. This is carried out again in his paper \cite{T30}, and it turns out to be important
in the work on the complex structure of Teichm\"uler space.
He then considers general problems on quasiconformal mappings, and in particular the following, which we state using his own words:
\begin{problem}
Let a sufficiently regular topological mapping of the circle $\vert z\vert =1$ onto itself be given. One wants to extend it to a quasiconformal topological mapping of the disk $\vert z\vert\leq 1$ onto itself so that the maximum of dilatation quotients becomes the smallest possible. 
\end{problem}
Teichm\"uller adds: ``Unlike what is before, the boundary must hence not only go over itself in totality but for every boundary point the image point is prescribed. This is as if all boundary points were distinguished."
 He makes several guesses and conjectures. He notes: ``I would like to call this statement not a conjecture but only a stimulation for a search for exact conditions." He then passes to a generalization to an arbitrary bordered surface with a distinguished arc on the boundary. 

Chapter 29, titled \emph{A metric}, starts at
\S 160. He defines the following metric on a Riemann surface: Consider a Riemann surface $\mathfrak{h}$ of negative Euler characteristic equipped with a variable interior point $\mathfrak{p}$. The space of pairs $(\mathfrak{h},\mathfrak{p})$ is a two-dimensional manifold. A Finsler metric can be defined on that space by considering the logarithm of the dilatation of the extremal quasiconformal mappings $\mathfrak{h}$ to itself taking $\mathfrak{p}$ to $\mathfrak{p}'$. 
Teichm\"uller addresses several questions about this metric and he says that in the case where $\mathfrak{h}$ is the 3-punctured sphere, the space obtained is the hyperbolic plane. He notes that one may make the same study  with distinguished points on the boundary.

At the end of the paper, Teichm\"uller addresses the question of why one studies quasiconformal mappings.  He makes connections with very early work of Tissot\footnote{Nicolas Auguste Tissot  (1824-1897) is a famous cartographer, who also taught mathematics at the \'Ecole Polytechnique.} on  geographical map drawing. He mentions that  Picard and Ahlfors proved theorems about these mappings which were direct generalizations of theorems on conformal mappings, but that they were nevertheless interested in conformal mappings, and their results, he says, ``taught us in which way the presupposition of conformality restricts the behavior of a mapping." He then adds: ``This view is outdated," and he reviews the new results and theories that quasiconformal mappings may lead to. He announces further results in future papers, in particular, detailed proofs of the results he announces in the present paper.

  The results of \cite{T20} are complemented  in the papers  \cite{T29}, \cite{T23}, \cite{T24} and \cite{T25}.

 In the paper \emph{Bestimmung der extremalen quasikonformen Abbildungen bei geschlossenen orientierten Riemannschen Fl\"achen} (Extremal quasiconformal maps of closed oriented Riemann surfaces) \cite{T29}, Teichm\"uller completes the proofs of the results he obtained in \cite{T20}. A translation of the paper is available \cite{T29A}, and a detailed commentary in \cite{T29C}.
  Let us quote from the introduction of that paper.
  \begin{quote}\small
In a longer article\footnote{Extremale quasikonforme Abbildungen und quadratische Differentiale \cite{T20}.} I gave heuristic arguments for the existence of an extremal quasiconformal mapping under a certain class of constraints, that means a map with minimal upper bound of the dilatation quotient. Additionally, I gave an analytic description of this map. I confirmed this there in the simplest examples, especially for the {\it torus\/}. Later\footnote{Vollst\"andige L\"osung einer Extremalaufgabe der quasikonformen Abbildung \cite{T24}.} I could add the more difficult case of a {\it pentagon\/}. Only now, I succeeded in proving my conjecture on the existence and the analytic form of the extremal quasiconformal maps to its full extent.

In 1939, it was risky to publish a lengthy article entirely built on conjectures. I had studied the topic thoroughly, was convinced of the truth of my conjectures and I did not want to keep back from the public the beautiful connections and perspectives to which I had arrived. Moreover, I wanted to encourage attempts for proofs. I acknowledge the reproaches that have been made to me from various sides, even today, as justifiable but only in the sense that an {\it unscrupulous\/} imitation of my procedure would certainly lead to a proliferation of our mathematical literature. But I never had any doubt about the correctness of my article, and I am glad now to be able to actually prove the main part of it.

At that time I was missing an exact theory of {\it modules\/}, the conformal invariants of closed Riemann surfaces and similar principal domains. In the meantime, particularly with regard to the intended application to quasiconformal maps, I developed such a theory. I will have to briefly report on it elsewhere. The present proof does not depend on this new theory, and instead works with the notion of {\it uniformization\/}. However, I think one will have to combine both to bring the full content of my article \cite{T20} into mathematically exact form.
\end{quote}

This paper was read by Ahlfors and Bers, who after that acknowledged the soundness of all the results stated by Teichm\"uller in his longer paper \cite{T20}. We shall quote Ahlfors and Bers, in the last section of the present paper. Ahlfors wrote in his review [MR0017803]: ``The restriction to closed surfaces is probably inessential. It is proved that every class of topological maps of a closed surface onto another contains an extremal quasiconformal map, i.e., one for which the upper bound of the dilatation is a minimum. [...] The surfaces are represented by means of Fuchsian polygons, the parameters being coefficients in a set of independent automorphisms. The proof is again based on the method of continuity."

The paper \emph{Ein Verschiebungssatz der quasikonformen Abbildung} (A displacement theorem for quasiconformal mapping) \cite{T31} is purely on the quasiconformal theory.  Teichm\"uller solves the following question: Describe the quasiconformal map $f$ from the unit disc to itself such that:
\begin{itemize}
\item $f$ is the identity on the boundary circle;
\item $f(0)=-x$, where $x$ is some given number $0<x<1$;
\item $f$ has the smallest dilatation among the maps satisfying the above two properties.
\end{itemize}
The tools are those introduced in \cite{T20}. By taking two coverings of the disk ramified at $0$ and $-x$ respectively, the author reduces the question to that of finding the extremal map between two ellipses in the Euclidean plane with horizontal and vertical axes, which are obtained from each other by a quarter-circle rotation, such that the map preserves the horizontal axes and the vertical axes and is affine on the boundary of the ellipse. Teichm\"uller shows that the natural affine map between the ellipses is the extremal map with the given boundary conditions. See \cite{Alberge} for a detailed commentary on this paper, and a review of some the developments. The paper is translated in \cite{Karbe}.

Finally, we comment on the paper \emph{Ver\"anderliche Riemannsche Fl\"achen} (Variable Riemannian surfaces) \cite{T32}. This is the last one that Teichm\"uller wrote on the problem of moduli. Together with \cite{T20}, it  is the most important  one that he wrote on the subject. 

In this paper, the author presents a completely new point of view on Teichm\"uller space, through complex analytic geometry. The paper is difficult to read, both for analysts and for low-dimensional topologists and geometers  because of its very concise style, and also because it depends heavily on the language of algebraic geometry (function fields, divisors, valuations, places and so on).  We state the main existence theorem:\begin{theorem} \label{th:T}
There exists an essentially unique globally analytic family of marked Riemann surfaces $\underline{\mathfrak{H}}[\frak{c}]$, where $\mathfrak{c}$ runs over a $\tau$-dimensional complex analytic manifold $\mathfrak{C}$ such that for any marked Riemann surface $\underline{\mathfrak{H}}$ of genus $g$ there is one and only one $\mathfrak{c}$ such that the Riemann surface  $\underline{\mathfrak{H}}$ is conformally equivalent to an $\underline{\mathfrak{H}}[\frak{c}]$ and such that the family $\underline{\mathfrak{H}}[\frak{c}]$ satisfies the following universal property: If $\underline{\mathfrak{H}}[\frak{p}]$ is any globally analytic family of Riemann surfaces with base $\mathfrak{B}$, there is a holomorphic map $f:\mathfrak{B}\to\mathfrak{C}$ such that the family $\underline{\mathfrak{H}}[\frak{p}]$ is the pull-back by $f$ of the family $\underline{\mathfrak{H}}[\frak{c}]$. 
\end{theorem}

  In this statement, $\mathfrak{C}$ is Teichm\"uller space and $\underline{\mathfrak{H}}[\frak{c}]\to  \mathfrak{C}$ is a fiber bundle over $\mathfrak{C}$ where the fiber above each point is a marked Riemann surface representing this point.   The fiber bundle  $\underline{\mathfrak{H}}[\frak{c}]\to  \mathfrak{C}$ is the Teichm\"uller universal curve. The theorem says that any analytic family of Riemann surfaces is obtained from the universal Teichm\"uller curve by pull-back by a certain holomorphic map. It follows easily from the uniqueness property\footnote{This is the statement made in Theorem \ref{th:T} that the family $\underline{\mathfrak{H}}[\frak{c}]\to  \mathfrak{C}$ is ``essentially unique." From the context, this means that it is unique up to composition by elements of the mapping class group.} that the automorphisms group of the universal Teichm\"uller curve is the extended mapping class group. This fact also remained unnoticed.\footnote{A proof of this result is contained in a paper by Andrei Duma (1975) who apparently was not aware of the fact that Teichm\"uller had already proved this result, see \cite{Duma}. Duma's proof does not use Royden's theorem on the complex structure automorphism group of Teichm\"uller space. Earle and Kra had already given a proof of that theorem, using Royden's theorem, without being aware that it was known to Teichm\"uller. We thank Bill Harvey and Cliff Earle for bringing this to our attention. See the MR review of Duma's paper by Earle (MR0407323) and the Zentralblatt review by Abikoff (Zbl 0355.32021).} An English translation \cite{Annette} as well as a commentary \cite{AAJP} of this paper are available. 
 
 The paper \cite{T32} did not attract much attention and it was practically never cited, especially compared with Teichm\"uller's previous papers on the subject, which were analyzed and commented by Ahlfors and Bers.\footnote{Ahlfors, in his 1960 paper \cite{Ahlfors1960},  writing about the complex structure of Teichm\"uller space, does not consider  that Teichm\"uller's approach in \cite{T32} is conclusive. He refers only to the quasiconformal approach in the paper \cite{T20}. He writes: ``In the classical theory of algebraic curves many attempts were made to determine the `modules' of an algebraic curve. The problem was vaguely formulated, and the only tangible result was that the classes of birationally equivalent algebraic curves of genus $g>1$ depend on $6g-6$ real parameters. More recent attempts to go to the bottom of the problem by more powerful algebraic methods have also ended in failure. The corresponding transcendental problem is to study the space of closed Riemann surfaces and, if possible, introduce a complex analytic structure on that space. In this direction, considerable progress has been made. The most important step was taken by Teichm\"uller \cite{T32}."
 But in his 1964 paper \cite{A1964}, Ahlfors emits a different opinion: ``It is only fair to mention, at this point, that the algebraists have also solved the problem of moduli, in some sense even more completely than the analysts. Because of the different language, it is at present difficult to compare the algebraic and analytic methods, but it would seem that both have their own advantages." } Grothendieck is among the first who understood the importance of that paper, about 15 years after its publication. We shall review his work on the subject in \S \ref{s:Gro}.

     \section{Weil}\label{s:Weil}

       Weil got interested very early in the theory of moduli. Together with Grothendieck, Kodaira, Ahlfors and Bers, he was influenced by Teichm\"uller's work, and probably more than the others, even though he did not publish much on this topic. Ahlfors and Bers were analysts, and Weil had a much broader vision, steeped in Riemann, and with a special taste for number theory. In his \emph{Souvenirs d'apprentissage},  he recalls the following, from the year 1925 (he was 19) (\cite{Weil-Souvenirs} p. 46 of the English translation): ``In the course of my walks, I would even stop to open a notebook of calculations on Diophantine equations. The mystery of Fermat's equations attracted me, but I already knew enough about it to realize that the only hope of progress lay in a fresh vantage point. At the same time, reading Riemann and Klein had convinced me that the notion of birational invariance had to be brought to the foreground. My calculations showed me that Fermat's methods, as well as his successors', all rested on one virtually obvious remark, to wit: If $P(x,y)$ and $Q(x,y)$ are homogeneous polynomials algebraically prime to each other, with integer coefficients, and if $(x,y)$ are integers prime to each other, then $P(x,y)Q(x,y)$ are `almost' prime to each other, that is to say, their G.C.D. admits a finite number of possible values; if, then, given $P(x,y)Q(x,y)=z^n$, where $n$ is the sum of the degrees of $P$ and $Q$, $P(x,y)$ and $Q(x,y)$ are `almost' exact $n$-th powers. I attempted to translate this remark into a birationally invariant language, and had no difficulty in doing so. Here already was the embryo of my future thesis." In 1926, Weil spent some time in G\"ottingen, with a grant from the Rockefeller foundation, after a few months spent in Italy. He makes a picturesque description of this visit in his \emph{Souvenirs d'apprentissage} (\cite{Weil-Souvenirs} p. 49ff): ``I chose Courant in G\"ottingen, because of linear functionals [...] Courant extended a cordial welcome to me [...] I  started explaining my ideas on functional calculus. [...] Courant listened patiently. Later on I learned that as of that day he concluded that I would be `\emph{unproduktiv}.' Leaving his house, I met his assistant, Hans Lewy, whose acquaintance I had made the day before. He asked me, `Has Courant given you a topic?' I was thunderstruck: neither in Paris nor in Rome had it occurred to me that one could `be given' a topic to work on [...] I learned little from Courant and his group."\footnote{Weil had nevertheless a great respect for Courant. He writes, in the same pages: ``It has sometimes occurred to me that God, in His Wisdom, one day came to repent for not having had Courant born in America, and He sent Hitler into the world expressly to rectify this error. After the war, when I said as much to Hellinger, he told me, `Weil, you have the meanest tongue I know'."} It seems that Weil learned everything on his own, in books. In fact, this is how he learned to read and write, and this also how he learned Latin and Greek, when he was a child. His thesis, defended in 1928, titled \emph{Les courbes alg\'ebriques} and written under Hadamard, was completely innovative and at the same time rooted in the works of Riemann, Abel and Jacobi. In a letter written to his sister Simone, Andr\'e Weil writes: ``My work consists in deciphering a trilingual text."\footnote{The letter is quoted in \cite{Hindry}, where the author adds: ``Roughly speaking, this consists in the analogies between number fields (arithmetic), function fields (analysis) and Riemann surfaces (geometry and topology).}

In any case, unlike several others mathematicians interested in Teichm\"uller theory, Weil was able to understand Teichm\"uller's papers written in the language of algebraic geometry. He was also a differential geometer and indeed his main contribution in Teichm\"uller theory lies in the local differential geometry of that space. 

It seems that Weil is the first who gave the name \emph{Teichm\"uller space} to the space discovered by Teichm\"uller, and in any case, he gave this name independently of others: in his 
 comments on the paper he wrote for Emil Artin's birthday  (\cite{Weil1979} Volume II p. 546), he writes: ``I was asked for a contribution to a volume of articles to be offered to Emil Artin in March 1958 for his sixtieth birthday; this led me to the decision of writing up my observations, even incomplete, on the moduli of curves and on what I called ``Teichm\"uller space."

 There is another terminology which Weil introduced, which does not survive, but which is interesting. In his 1958 Bourbaki seminar \cite{Weil1958c} and in his paper (\cite{Weil1958a}), Weil called the object today called a \emph{marked Riemann surface} a \emph{Teichm\"uller surface}, and he also introduced the terminology \emph{Torelli surface}.  A \emph{Teichm\"uller surface} is a pair $(S,[f])$ where $S$ is a Riemann surface $S$ homeomorphic to $S_0$ and $[f]$ is the homotopy class of some orientation-preserving homeomorphisms  $f:S_0\to S$. The Teichm\"uller space of $S_0$ is then the space of Teichm\"uller surfaces $(S,[f])$ up to the following relation 
 $(S,[f])\equiv (S',[f'])$ if there is a holomorphic map $h:S\to S'$ such that $h\circ f$ is homotopic to $h'$. A \emph{Torelli surface} is a Riemann surface $S$ homeomorphic to $S_0$, together with a homomorphism from the fundamental group of $S_0$ to the homology group of $S$ that is induced by an orientation-preserving homeomorphism between the two surfaces. 
A Torelli surface is, like a Teichm\"uller surface, a \emph{reinforced} Riemann surface, that is, it carries more structure\footnote{The adjective is used by Weil in \cite{Weil1958a}.} than a Riemann surface in the same sense in which a Teichm\"uller surface is a reinforced Riemann surface. The choice of names reflects clearly the sequence of maps $\mathcal{T}_g\to \mathrm{Tor}_g\to \mathcal{M}_g$ between the Teichm\"uller, the Torelli and the Riemann spaces.

Weil also highlighted the relation with the period map and the Siegel upper half-space. After a choice of a basis for $\pi_1(S)$, a Torelli structure determines the images of the generators in $H_1(S)$, which are generators of this group. With these generators of $H_1(S)$, one can define the so-called ``renormalized differentials of the first kind," that is, those whose matrices are of the form $(I_g,Z)$ where $I_g$ is the unit $g\times g$ matrix and $Z$ is a symmetric matrix whose imaginary part is positive non-degenerate. Such a matrix can be considered as a point of the Siegel space. Torelli's theorem states that the $Z$ determines $S$ up to isomorphism.

Let us say a few words on some of Weil's early works, which will show how he arrived naturally to the subject.
  
  Weil was interested since his early years in complex function theory of several variables, and in fact, he is one of the promoters of that theory, which was rather dormant since some works Poincar\'e did on that subject in the nineteenth century. In an interview he gave to the \emph{Notices of the AMS} in 1999, Henri Cartan, who is the founder of one of the major schools in that field, says the following (\cite{Cartan-AMS} p. 784): ``I believe it was Andr\'e Weil who suggested that it could be interesting. He told me about the work of Carath\'eodory on circled domains. That was the beginning of my interest."\footnote{Cartan completed his doctorate in 1928, the same year as Weil, under Montel, on functions of one complex variable. This conversation with Weil probably occurred soon after. In 1930, Cartan published a paper titled \emph{Les fonctions de deux variables complexes et les domaines cercl\'es de M. Carath\'eodory} \cite{CC}, and several others in the same year and the one that followed, in which he studied and classified bounded domains in $\mathbb{C}^2$ with infinitely many automorphisms having one fixed point. We note that among these papers is one of the two papers that Henri Cartan co-authored with his father (who was not so much involved in these questions), \cite{Cartan-Cartan}.} Soon after, Weil published the two papers \cite{W1932b} and \cite{1935d} in which he introduced a method of integration in the theory of several complex variables, generalizing the Cauchy integral. In the paper \cite{Weil1938a}, he started the study of vector bundles over  algebraic curves, with fibers of any dimension $\geq 1$. In his comments in Vol. I of his \emph{Collected papers}, Weil recalls that in that paper, he ``subconsciously" endeavored to build a moduli space for these bundles, but that at that time he did not have clear ideas on the subject. He states that a special case of the construction was implicit in the classical literature, and it amounts to that of the Jacobian of a curve. He also says that his hope was  to define functions which would generalize the theta functions, but that he did not succeed. The idea, he tells us, even though it was unclear to him, was to conceive a ``moduli space" for the analogue of a Fuchsian group $G$ in dimension $\geq 1$. This would be a quotient of a variety of equivalence classes of representations of $G$ into $\mathrm{GL}(r,\mathbb{C})$, equipped with its natural complex structure. Weil made the observation that there exist equivalence classes of representations that do not contain unitary ones, but that such a unitary representation, if it exists, is essentially unique. He writes in his comments: ``It seems to me that I had the feeling of the important role of unitary representations in the future development of the theory." The famous \emph{Weil conjectures} which he made in 1949 \cite{Weil1949b}, which led to an huge amount of work in algebraic geometry during several decades, with a culmination in the works of Grothendieck and Deligne, are also related to Riemann's ideas. In some sense, they concern the applications of the topological methods of Riemann to the setting of fields of finite characteristic. Many other works of Weil are related in one way or another  to the subject of the present paper, but we shall concentrate only on some of them, which are very closely related to Riemann's and Teichm\"uller's works.   
     
     In the years 1957 and 1958, Weil wrote four papers on the problem of moduli:
 \begin{itemize}
 
 \item A Bourbaki seminar, on the Torelli theorem, titled \emph{Sur le th\'eor\`eme de Torelli} \cite{Weil1957}, which we mentioned in \S \ref{s:Torelli}.
 
\item A Bourbaki seminar titled \emph{Modules des surfaces de Riemann} \cite{Weil1958c}.
 
 \item An unpublished manuscript, in English, carrying the same title as the Bourbaki seminar, dedicated to Emil Artin on his sixtieth birthday \cite{Weil1958a}.
 
\item An unpublished report on an AF contract \cite{Weil1958b}.
\end{itemize}

This collection of papers constitutes a brief account of Riemann's moduli space, of Torelli space  and of Teichm\"uller space, and they contain the essential ideas of several constructions. The two unpublished papers (and not the Bourbaki seminars) appear in Weil's  \emph{Collected Papers}  \cite{Weil1979}. We shall review these four papers, together with the comments that Weil made on them.  
   The papers and the comments provide us with a clear summary of why Weil  was interested in Teichm\"uller theory, and also what he learned from Teichm\"uller's papers and from Kodaira, Spencer, Ahlfors and Bers. The papers also contain open questions. The report \cite{Weil1958a} starts with:
 
\begin{quote}\small
The purpose in the following pages is partly to clarify my own ideas on an interesting topic, at a stage when they are still unripe for publication.  In speaking of these ideas as ``my own," my intention is not to claim originality for them. They are little more than a combination of those of Teichm\"uller with the ideas on the variation of complex structures, recently introduced by Kodaira, Spencer and others into the theory of moduli.
 \end{quote}

Weil describes two approaches to the complex structure on Teichm\"uller space in the paper \cite{Weil1958b}.
The first one uses harmonic  Beltrami differentials and quasiconformal maps, and he attributes it to Bers. The second one uses the period map. With the latter, he recovers a construction of Rauch of the complex structure in the complement of the hyperelliptic locus. He then says that this complex structure can be extended to Teichm\"uller space by ``using well-known general theorems on analytic spaces."\footnote{Such an extension was published by Ahlfors later, in his 1960 paper \cite{Ahlfors1960}. Ahlfors does not mention Weil's Bourbaki seminar. In his comments to his \emph{Collected papers} edition, Vol. 2, p. 122, Ahlfors writes: ``H. Rauch published two notes (On the transcendental moduli, 1955, and On the moduli, 1955) in which he settled the problem for non hyperelliptic surfaces. For my part I was also in possession of a proof in the nonhyperelliptic case, but I had held up publication until I was able to construct the complex structure for the whole space $\mathcal{T}_g$."}

In his 1958 paper dedicated to Emil Artin, Weil writes (\cite{Weil1958a}, p. 413):
\begin{quote} \small  Perhaps the most remarkable of Teichm\"uller's results is the following: when provided with a rather obvious ``natural" topology, the set $\Theta$ of all classes of mutually isomorphic Teichm\"uller surfaces is homeomorphic to an open cell of real dimension $6g-6$. This global result will neither be used nor discussed in the following pages, the chief purpose of which is to consider the local properties of $\Theta$ and to define on it a ``natural" complex analytic structure, of complex dimension $3g-3$, and a ``natural" Hermitian metric.
\end{quote}
In the same paper, he writes:
\begin{quote}\small
[...] In order to justify the statements that we have made so far, we shall make use of the Kodaira--Spencer technique of variation of complex structures. This can be introduced in an elementary  manner in the case of complex dimension 1, which alone concerns us here; this, in fact, had already been done by Teichm\"uller; but he had so mixed it up with his ideas concerning quasiconformal mappings that much of its intrinsic simplicity got lost. Perhaps the worst feature of his treatment, in the eyes of the differential geometer, is that his extremal mappings are destructive of the differentiable structure; this corresponds to the fact that his metric on $\Theta$ is almost certainly not to be defined by a $ds^2$, even though it is presumably a Finsler metric.\footnote{Weil writes ``is presumably a Finsler metric," but Teichm/"uller proved it is in his paper \cite{T20}. }
\end{quote}
 Weil, who was chiefly interested in the complex structure of Teichm\"uller space, does not mention the paper \cite{T32}. His approach was, like the one of Ahlfors and Bers, to extract the complex structure from the almost-complex structure through the theory of the Beltrami equation developed in Teichm\"uller's paper \cite{T20}. In his Bourbaki seminar, he starts with a  differentiable function $\mu$ on the upper half-plane $\mathbb{H}^2$ equipped with a discrete group of isometries $\Gamma$ representing a surface $S$,  such that $\vert \mu\vert <1$ and such that the tensor $\mu d\overline{z}/dz$ is invariant by $\Gamma$ (that is, $\mu$ is of type $(-1,1)$). The differential $\zeta= dz+\mu d\overline{z}$ defines a new complex structure on $\mathbb{H}^2$, and one gets a Riemann surface $S(\mu)$. Given a marked Riemann surface $S$, Weil denotes by $\mathrm{Cl}(S)$ the corresponding element in Teichm\"uller space. We shall use instead the notation $[S]$ which is common today. We have a map $\mu\mapsto [S(\mu)]$. The theory that Weil sketches is summarized in the following
\begin{theorem}\label{Th:T}
If $\mu=\mu(z,u)$ depends continuously (resp. differentiably, resp. holomorphically) on a parameter $u$, then the map $u\mapsto [S(\mu_u)]$ is continuous (resp. differentiable, resp. holomorphic).
\end{theorem}
Weil notes that if $F_\mu$ is the homeomorphism of the plane $\mathbb{H}^2$ which conformally represents $\mathbb{H}^2$ equipped with the structure defined by the form $\zeta$ on $\mathbb{H}^2$ equipped with its usual complex structure (defined by $dz$), and if we normalize $F_\mu$ so that it preserves the points $0,1,\infty$, then, if $\mu=\mu_u$ depends continuously (resp. differentiably) on a parameter $u$, the same holds for $F_\mu$.\footnote{Regarding this result, Weil writes in \cite{Weil1958c} p. 416: ``Some reliable ellipticians guaranteed this to me."} He then adds that the analogous statement for a holomorphic result is false.\footnote{Weil writes in \cite{Weil1958c} p. 416: ``One should not think that if $\mu_u$ depends in a holomorphic manner on $u$, the same is true for $F_{\mu_{u}}$, since this is false."} Contrary to what Weil says, this result is correct, and it is the famous Ahlfors-Bers Riemann mapping theorem for variable metrics \cite{ABers}.

Then, Weil studies the differentiable structure determined by Theorem \ref{Th:T}. He introduces the  notion of an infinitesimal variation of a structure in the sense of Kodaira-Spencer and the notion of a trivial infinitesimal variation, which amounts to the fact that $(\partial \mu /\partial u)_{u=0}$
 is of the form  $\partial \xi/\partial \overline{z}$ where $\zeta$ is a complex-valued form of type $(-1,0)$ on $\mathbb{H}^2$.  He concludes (``heuristically," he says, postponing the proof to later) that one can identify the (real) tangent vectors to Teichm\"uller space at a point $[S_0]$ with a ``Kodaira-Spencer space" defined as the quotient of the (real) vector space of functions of type $(-1,1)$ defined on $\mathbb{H}^2$ quotiented by the subspace of functions $\partial \xi/\partial\overline{z}$ where $\zeta$ is of type $(-1,0)$. Thus, the complex structure is defined using the Kodaira-Spencer theory.

In the paper \cite{Weil1958a}, Weil introduces the Riemannian metric on Teichm\"uller space which has become known as the Weil-Petersson metric. This is defined using a version of an inner product which was defined in 1939 by Petersson\footnote{This is Hans Petersson (1902-1984). Lehner writes about Petersson, in Chapter 1 (Historical Development) of his book \cite{Lehner} (p. 35): ``The most important contributor to the theory of automorphic functions in recent times is H. Petersson, whose investigations begin about 1930. He was a student of Hecke and much of his work consists of extending to more general discontinuous groups what Hecke developed for congruence subgroups of the modular group. Petersson, like Klein and Hecke, is greatly interested in the correspondence between Riemann surface theory and automorphic function theory. Because he insists on considerable generality, his papers are hard to read, but it is also true that he returns again and again to the classical examples from which the general theory sprang, reinterpreting and deepening them with the newest results." After recalling the Petersson product, Lehner writes: ``Petersson's investigations of the new Hilbert spaces revolutionized the theory of automorphic forms of negative dimension. Formerly difficult theorems could now be proved by the methods of liner algebra."} in the context of automorphic forms, cf. \cite{Petersson1} and \cite{Petersson2}. Petersson defined the scalar product of pairs $\phi$ and $\psi$ of automorphic forms of the same class by the formula
\[(\phi,\psi)=\int\int \phi(z)\overline{\psi(z)}y^r\frac{dxdy}{y^2}\]
where $z=x+iy$ in the upper half-plane, with $r\geq 0$ being some integer weight and where the integral is over a fundamental domain.
 In his first papers, Petersson pointed out several applications of this product to the theories of modular forms and of automorphic forms. Many other applications were discovered later by others.  
 Weil noticed that this product gives a metric on Teichm\"uller space. 
Let us quote his Bourbaki seminar in which he defines the Weil-Petersson metric:
 \begin{quote}\small
 One observes immediately that the space of quadratic differentials on $S_0$ is equipped with a well-known Hermitian metric (a particular case of that of Petersson in the theory of automorphic forms), given,
 for
 $w=qdz^2$, $w'=q'dz^2$, by
 \[(w,w')=i\int\int_{S_{0}}q\overline{q'}y^2dzd\overline{z}\]
 with $y=\mathrm{Im}(z)$. After checking some small poor differentiability conditions, this means that this form defines on Teichm\"uller space a ``natural" Hermitian structure, and one checks (by a stupid computation) that it is K\"ahler.
The preceding can also be expressed by saying that any quadratic differential $w$ on $S_0$ defines a covector (with complex values) at the equivalence class $[S_0]$ of $S_0$ in Teichm\"uller space, a covector which, by definition, is of type (or ``bidegree") $(1,0)$, for the quasi-complex structure that we defined above.\footnote{This is the almost-complex structure defined by the equivalence classes of Beltrami differentials.}  Then a trivial computation allows one to check that for such a structure,  the map $u\mapsto [S_u]$ is holomorphic if $\mu_u$ depends holomorphically on a parameter $u$ with values in $\mathbb{C}^n$.
If we take then, to be more precise:
\[\mu_u(z)=\sum_{i=1}^{3g-3}u_i y^2 \overline{q}_i (z),\]
where the $q_i$ are such that the $q_i dz^2$ form a basis of the vector space of quadratic differentials on $S_0$, we deduce that the map $u\mapsto [S_0]$ is an isomorphism from a neighborhood of $0$ in $\mathbb{C}^{3g-3}$ to a neighborhood of $[S_0]$ in Teichm\"uller space, in the sense of quasi-complex structures. This confirms the fact that the latter space is complex analytic of dimension $3g-3$. [N.B. This reasoning is due to Bers (personal communication)].
 
In the same vein, we can now treat by a computation, without difficulty, a variety of questions concerning the differentials, at a given point in Teichm\"uller space, of functions defined on that space, and for instance:

A.  Coefficients for the substitutions of the [Fuchsian] group $\Gamma$. In this way, we find that if this group is normalized, the $6g-6$ coefficients of the first $2g-2$ ``distinguished" generators may be taken as local coordinates (in the sense of real analytic structures) on Teichm\"uller space. Thus, this space is represented as an open set on $\mathbb{R}^{6g-6}$;

B. Periods of normalized integrals of the first kind on $S_{\mu}$. We can thus check the theorem of Rauch, which says the following: In the Siegel space of symmetric matrices with positive non-degenerate imaginary part, the points corresponding to the ``Torelli surfaces" form an analytic subspace of complex dimension $3g-3$, whose singular points are those which correspond to the hyperelliptic curves. The latter form a complex analytic variety (with no singular point) of complex dimension $2g-1$. In the neighborhood of every point of this space corresponding to a non hyperelliptic surface $S_0$, we can write, as local coordinates on that set, every system of $3g-3$ periods $p_{ij}$ such that the quadratic differentials $\zeta_i\zeta_j$ (where $\zeta_1,\ldots,\zeta_g$ are the normalized differentials of the first kind) are linearly independent over $\mathbb{C}$.

 \end{quote}

Weil concludes his paper \cite{Weil1958a} with the following (p. 389):
 \begin{quote} \small
 This raises the most interesting problems of the whole theory: is this a K\"ahler metric?\footnote{In his Bourbaki seminar paper \cite{Weil1958c} which we quoted above, written the same year, Weil writes, on the contrary (on p. 414 and on p. 418), that the space \emph{has} a K\"ahler structure, and that this can be checked ``by a stupid computation." The same existence result is stated in his AF Report \cite{Weil1958b}, p. 392 of Vol. II of Weil's \emph{Collected papers}.} has it an everywhere negative curvature? is the space $\Theta$, provided with its complex structure and with this metric, a homogeneous space?
It would seem premature even to hazard any guess about the answers to these questions.
 \end{quote}

 The first two questions were answered positively by Ahlfors in the two papers \cite{Ahlfors-Some} and \cite{Ahlfors-Curvature}. In the paper \cite{Ahlfors-Some}, Ahlfors writes in a footnote:
``According to an oral communication, the fact has been known to Weil, but his proof has not been published." In the second paper \cite{Ahlfors-Curvature}, Ahlfors writes:
``Intrinsic definitions lead to a 
metric on $T_g$, introduced by Teichm\"uller, to a Riemannian structure whose 
use was suggested by A. Weil, and finally to a complex analytic structure 
of dimension $3g-3$. It was proved by Weil, and again in \cite{Ahlfors-Some}, with very 
little computation, that the Riemannian metric is K\"ahlerian with respect to 
the complex structure." 
  In a note contained in his \emph{Collected papers} edition, Vol. 2, p. 155, Ahlfors writes:
          ``Weil knew that the Petersson metric was K\"ahlerian, but had not published the proof. This turned out to be an almost immediate consequence of the calculations in \cite{Ahlfors-Some}, and in \cite{Ahlfors-Curvature} I showed through hard work that the metric has negative Ricci and sectional curvature." This last information contrasts with what Weil says in his Bourbaki seminar, where he gives the impression that everything was easy for him. In any case,  it is impossible to know with certainty who first proved that the Weil-Petersson metric is K\"ahler. The homogeneity question that Weil addressed concerning the complex structure was answered in the negative by Royden in 1969 \cite{Royden1969}.  The homogeneity question for the Weil-Petersson metric is a consequence of the 2002 result of Masur and Wolf \cite{MW2002}. 

To understand his motivations, it should also be noted that Weil, at the time he was also working and making conjectures on the K\"ahler structure of Teichm\"uller space, was also working on the K\"ahler geometry of K3 surfaces. In fact, in his final report to the AF \cite{Weil1958b}, Weil reports on two problems: the problem of moduli, which we are discussing here, and ``a study of the K\"ahler varieties topologically identical with the non-singular quartics in projective 3-space (henceforward called K3 surfaces)."\footnote{Here also, the name seems to be due to Weil. In his commentary on this paper, Weil writes that the name K3 is in honor of Kummer, K\"ahler, Kodaira and of the beautiful Kashmir mountain K2. A K3 surface is Calabi-Yau of smallest dimension which is not a torus. Most K3 surfaces are not algebraic varieties. A K3 surface is also K\"ahler (a theorem of S.-T. Siu).} It is interesting to note that the same group of mathematicians contributed to both problems, including Siegel, Kodaira, Spencer and Andreotti, and the two theories share many common problems and techniques: the existence of a complex-analytic structure and of a K\"ahler metric, intersection matrices and period maps, the use of polarization of Abelian varieties, the theory of automorphic functions, etc. Finally, we mention that in the same year, Weil published his important book on K\"ahler manifolds, \emph{Introduction \`a l'\'etude des vari\'et\'es k\"ahl\'eriennes} \cite{Weil-Kahler}. He had been working on the subject for at least ten years; see \cite{Weil1947b} and \cite{Weil1949d}.
     
In discussing the theory of the Beltrami equation, Weil several times refers to Bers, and for some other points he refers to ``analysts" and to  some ``reliable ellipticians" (e.g. \cite{Weil1958c} p. 418) and to Kodaira and Spencer, with whom he said he had several discussions. In his report \cite{Weil1958b}, he writes: ``In order to give, in what follows, a coherent account of these topics, it will be necessary to include much of the work of my colleagues, and it would be unpractical to try unravel in detail what may belong to me and what belongs to each one of them. It should be understood that they deserve a large share of credit for the work described in this report."

  It seems that Weil left the subject after he wrote these four papers, and we have two conjectural reasons for that. The first one is that Weil saw he was superseded on these questions  by Ahlfors and Bers.\footnote{From Weil's comments on his paper for Emil Artin's birthday, in his \emph{Collected Works} (\cite{Weil1979} Vol. II p. 546): ``Soon after that, I noticed that on more than one point I was matching Ahlfors and Bers; these continued their research, and soon after they overtook me."} The other reason is related to Alexander Grothendieck's appearance on the scene. Grothendieck completely changed  the foundations of analytic geometry, and in doing this he was greatly influenced by Teichm\"uller's construction of the analytic structure of Teichm\"uller space.\footnote{These are Grothendieck's own words (\cite{G1}, p. 1): ``Chemin faisant, la n\'ecessit\'e deviendra manifeste de revoir les fondements de la G\'eom\'etrie analytique" (On the way, the necessity of revising the foundations of analytic geometry will be manifest).} We shall elaborate on this in the next section. The relation between Grothendieck and Weil was not always friendly. In fact, Weil hardly mentions Grothendieck in the comments he makes in his \emph{Collected Works}. In his paper \cite{MG}, Mumfor writes that ``Weil radiated cynicism about anyone else's abstractions." One certain fact is that Grothendieck's \emph{\'El\'ements de g\'eom\'etrie alg\'ebrique} (1960--1967) almost condemned Weil's \emph{Foundations of algebraic geometry} (1946) to oblivion.

 \section{Grothendieck}\label{s:Gro}  
    
In the academic year 1960-1961, Grothendieck\footnote{Alexander Grothendieck was born in Berlin, in 1928. It is impossible to say anything significant on his life in a few lines. We refer the reader to the collective book \cite{Schneps} and in particular the article \cite{Cartier} by Cartier.} gave a series of 10 talks at Cartan's seminar in which he presented Teichm\"uller's existence and uniqueness result on the complex analytic structure of Teichm\"uller space (Theorem \ref{th:T} above). This result was reformulated by Grothendieck and given a complete proof  in the language of algebraic geometry. The new point of view included the statement that Teichm\"uller space represents a functor.

   A large part of Grothendieck's work is developed in the setting of categories and functors, and his approach to Teichm\"uller theory was in the same spirit.  This abstract categorical setting is a unifying setting, and it is the result of Grothendieck's broad vision and a profound feeling of unity in mathematics, especially algebraic geometry, complex geometry, and topology. The fact that Grothendieck formulated the moduli problems in such a setting came out of a need. With the classical methods, it was possible to construct analytic moduli spaces, but not algebraic ones. In fact, the question of ``what is an algebraic deformation?" was one of the main questions formulated clearly by Grothendieck, using the language he developed. Paraphrasing and expanding the introduction in \cite{G1}, Grothendieck's goal in his series of lectures, which he announced in the first lecture, is the following:

\begin{enumerate}
\item To introduce a general functorial mechanism for a global theory of moduli. Teichm\"uller theory is one example to which this formalism applies, but the theory also applies to other families, like the family of elliptic curves (the case of genus 1), which  so far (according to Grothendieck) had not been made very explicit in the literature. 
\item To give a ``good formulation" of a certain number of moduli problems for analytic spaces. Grothendieck gives a precise formulation of several moduli problems in the setting of analytic spaces, including the moduli of Hilbert schemes of points or Hilbert schemes of subvarieties, or of Hilbert schemes with a given Hilbert polynomial. These were used as step-stones, but they also have an independent interest.  According to Grothendieck, for all these moduli problems, the state of the art in most of the situations is such that one could only conjecture some reasonable existence theorems. He nevertheless gives some existence and uniqueness results, and he proves in particular the existence of Teichm\"uller space as a universal object.

\item Under some ``projectivity  hypothesis" (this is the analogue of compactness in the algebraic geometry setting) for the morphisms that will be considered, to give some existence theorems for the problems in (2). \emph{This includes in particular the existence theorem for Teichm\"uller space.}\footnote{Grothendieck underlines.}
\item Grothendieck says that for that purpose, it will be necessary to reconsider the foundations of analytic geometry, by getting inspiration from the theory of schemes.\footnote{In this setting, an analytic space is defined as a system of local rings, each ring being a quotient of the ring of convergent complex power series in several variables. This is a wide abstraction of the idea of Riemann where the stress is on functions rather than on the space itself. The theory of schemes gives a way of gluing these local fields. One can also think of this in analogy with the manifolds defined by local charts, in the differential category.} In particular, it will be important to admit nilpotent elements\footnote{Recall that points in an analytic space are defined as ideals in certain rings of fields, and Grothendieck says here that in order to get the full strength of the theory, we must admit nilpotent elements in the ideals. For instance, one should be able to distinguish between the function $y^2$ and the function $y$ which have the same zero set, namely, the point $0$.} in the local fields defining analytic spaces,  and also in more general spaces that consist of families of spaces,  so that the theorems are stated with their full strength.
\end{enumerate} 

Let us now quote Grothendieck's statement:

 \begin{theorem}[Theorem 3.1. of \cite{G1}] \label{the:GG} 
 There exists an analytic space $T$ and a $\mathcal{P}$-algebraic curve $V$ above $T$ which are universal in the following sense: For every $\mathcal{P}$-algebraic curve $X$ above an analytic space $S$, there exists a unique analytic morphism $g$ from $S$ to $T$ such that $X$ (together with its $\mathcal{P}$-structure) is isomorphic to the pull-back of $V/T$ by $g$.
\end{theorem}

  In this statement an \emph{algebraic curve over an analytic space} is a family of algebraic curves which depends analytically on a parameter. The reader should notice the intertwining of the words ``algebraic" and  ``analytic" in the statement. This is one of the major ideas in modern algebraic geometry.
The analytic space $T$ which is referred to here is Teichm\"uller space and the $\mathcal{P}$-algebraic curve $V$ above $T$ is the Teichm\"uller universal curve. The term ``$\mathcal{P}$-algebraic" refers to a rigidification of the curves, and $\mathcal{P}$ is a functor, the so-called Teichm\"uller rigidifying functor.

It is probable that Grothendieck did not read any of Teichm\"uller's papers, but that he heard of them by word of mouth.  In the bibliographical references of the first article of the series \cite{G1}, in which he outlines the whole theory,  the only reference he makes (apart from references to his own works) is to Bers' paper \cite{Bers-Tata}. On p. 9 of this paper, Grothendieck writes: ``It is easy to check, using if necessary a paper by Bers \cite{Bers-Tata}, that the space which is introduced here axiomatically (and which we shall prove the existence) is homeomorphic to the Teichm\"uller space of analysts. It follows that Teichm\"uller space is homeomorphic to a ball, and therefore, is contractible, and in particular connected and simply connected. A fortiori, the Jacobi spaces of all levels are connected, and likewise the moduli space introduced in \S 5 as a quotient space of Teichm\"uller space. It seems that at the time being there is no proof, using algebraic geometry, of the fact that the moduli space is connected, (which, in algebraic geometry, is interpreted by saying that two algebraic curves of the same genus $g$ are part of a family of curves parametrized by a connected algebraic variety)."

One of the main tools is the introduction of categories whose objects are fibrations over complex spaces whose fibers are complex spaces and functors between them.  Grothendieck introduces several  functors in the theory. Some of them are \emph{rigidifying}  and some are \emph{representable}. Classical objects like projective spaces, Grassmannians, and Eilenberg-Mac Lane spaces, represent functors. Rigidifying functors are obtained by equipping the fibers with some extra structure.  Marking is a rigidification of the structure and it is transported from fiber to   fiber by analytic continuation. In this process a monodromy appears, which, in the case of the functor built from Riemann surfaces, is encoded by the mapping class group action on the fibers. The rigidifying functor resulting from the rigidification of Riemann surfaces is the so-called \emph{Teichm\"uller functor}, which is an example of a \emph{representable functor}. We recall the definition. Let $C$ be a category, $X$ an object in $C$ and $h_X$  a contravariant functor from $C$ to the category $\mathbf{(Ens)}$ of sets, defined by the formula 
\[h_X(Y)=\mathrm{Hom}(Y,X)\] at the level of objects. At the level of morphisms, if $Y$ and $Z$ are objects in $C$ and  $f:Y\to Z$ a morphism between them, then the image of $f$ is the map $p\mapsto p\circ f$ from the set $\mathrm{Hom}(Z,X)$ to the set $\mathrm{Hom}(Y,X)$. If $C$ is a category and \[F:C\to \mathbf{(Ens)}\] a contravariant functor from $C$ to the category of sets, then
  $F$ is said to be \emph{representable} if there exists an object $X$ in $C$ such that the functor $F$ is isomorphic to the functor $h_X$.
We say that $F$ is \emph{represented} by the object $X$. 
According to Grothendieck in \cite{Gro-Bou2},  the ``solution of a universal problem" always consists in showing that a certain functor from a certain category to the category  $\mathbf{(Ens)}$ is representable. The fact that a functor is representable makes the abstract category on which it is defined more understandable,  since it says that the objects and morphisms of this category can be replaced by objects and morphisms of the category of sets. 

The main result of Grothendieck in his series of lectures is the following:

 \begin{theorem}[\cite{G1} p. 8] \label{the:Functor} 

  The rigidifying Teichm\"uller functor, $\mathcal{P}$,  for curves of genus $g$, is representable.  
  \end{theorem}

In fact,  Teichm\"uller space was the first non-trivial example of a complex space representing a functor. Grothendieck also showed that the Teichm\"uller functor is universal in the sense that any rigidifing functor can be deduced from an operation of fiber product from the Teichm\"uller rigidifying functor (see \cite{G1} p. 7).

Theorems \ref{the:GG}  and \ref{the:Functor}  are equivalent, although there are differences in the formulations and in the  proofs.
Teichm\"uller first constructs Teichm\"uller space, and then moduli space. It seems that he thought it was impossible to construct moduli space directly. 
Grothendieck constructs moduli space directly as a space representing a functor in the algebraic category, so that the space is canonically equipped with an algebraic structure (and a fortiori an analytic structure). Teichm\"uller used a topological marking of Riemann surfaces in order to construct Teichm\"uller space. Grothendieck, in dealing with the singularities of moduli space caused by the non-trivial automorphisms of Riemann surfaces, introduced a marking by cubical differentials and another marking through level structures on the first homology.\footnote{In \cite{uni}, the authors mention the fact that one can deal with the problem of surfaces having nontrivial automorphisms by considering ``level structures," and they refer for this to the book \cite{HaMo}. No reference to the papers of Grothendieck  \cite{G1} \cite{G10} is given, either in \cite{uni} or in \cite{HaMo}.} The algebraic structure on moduli space cannot be lifted to Teichm\"uller space since the latter is an infinite cover of the former. In fact, Teichm\"uller space does not carry a natural algebraic structure.

Inspired by Teichm\"uller's ideas, Grothendieck proposed a general approach to the construction of moduli spaces of algebraic varieties, in particular Hilbert schemes for families of closed subvarieties of a given projective variety. This also led him to his theory of \emph{families} of algebraic varieties. The idea of marking allowed him to remove nontrivial automorphism groups of the varieties and to construct fine moduli spaces. The desingularization of the moduli spaces was obtained by looking at smooth covering varieties, cf. the papers \cite{GGG}, and also \cite{Gro-Bou1},  \cite{Gro-Bou2}.     

 Like Teichm\"uller's paper \cite{T32}, Grothendieck's papers are difficult to read for people not used to the abstract language of algebraic geometry.\footnote{Abikoff writes, in a review published in 1989 in the Bulletin of the AMS \cite{Abikoff-BAMS}, on an book on Teichm\"uller theory by Nag: ``First, algebraic geometers took us, the noble but isolated practitioners of this iconoclastic discipline, under their mighty wings. We learned the joys of providing lemmas solving partial differential and integral equations and various other nuts and bolts results. These served to render provable such theorems as: The ?\%$\sharp$\$! is representable."} On p. 9 of his first lecture \cite{G1}, Grothendieck writes: ``One can hope that we shall be able one day  to eliminate analysis completely from the theory of Teichm\"uller space, which should be purely geometric."  A detailed review of the 10 lectures in Cartan's seminar is contained in the paper \cite{AAPP}.

Grothendieck, in his lectures, raised several questions, among them the question of whether Teichm\"uller space is Stein (\cite{G1} p. 14). Bers and Ehrenpreis answered this question affirmatively in \cite{BE1964}, apparently without being aware that Grothendieck had raised it.

Grothendieck left Teichm\"uller theory just after he gave these lectures. The reason was certainly his complete investment in the foundations of general algebraic geometry. Mumford brought later a new point of view, based on geometric invariant theory, which allowed him to give an intrinsic construction of Riemann's moduli space and equipping it with the structure of an algebraic (quasi-projective) variety, after the work of Baily on the subject. Mumford's work  \cite{Mumford} is mentioned by Grothendieck in his last lecture, as well as Baily's work \cite{Baily}. In the introduction of that lecture, he writes:
\begin{quote}
The method indicated in the text bumps, in the context of schemes, on a difficulty related to taking quotients, which does not exist in the transcendental case. By a similar method, using in a more systematic way Picard schemes and their points of finite order, the lecturer was able to construct the Jacobi modular schemes $M_n$ \emph{of high enough level}, but, for lack of knowledge of whether $M_n$ is quasi-projective, it was not possible to take a quotient by finite groups in order to obtain the moduli spaces of arbitrary level, and in particular the classical moduli space $M_1$. These difficulties have just been overcome by Mumford, using a new theorem of passage to the quotient which can be applied to polarized Abelian schemes, and from there to curves.
\end{quote}

Grothendieck returned to Riemann surfaces several years later, after he had officially put an end to his extraordinarily intense mathematical activity. He introduced the subject which became known as the Grothendieck-Teichm\"uller theory, whose overall goal is to understand the absolute Galois group
$\Gamma_{\mathbb Q}=\mathrm{Gal}(\overline{\mathbb Q}/\mathbb Q)$ through its action on the fundamental group of the tower of moduli space $\mathcal{M}_{g,n}$  of algebraic curves 
of genus $g$ with $n$ punctures. The main reference is  \cite{Gro-esquisse}.
Several papers were written after Grothendieck. See in particular \cite{Loc} \cite{Lochak2012} \cite{NS2000}.
The work is also surveyed in \cite{AJP3}.

 \section{Ahlfors and Bers}\label{s:AB} 
  In the years that followed Teichm\"uller's work, a huge effort was made by  Ahlfors and Bers to provide proofs for the analytic part of his work, and in particular, on his existence and uniqueness theorem for extremal quasiconformal maps. We shall mainly quote some of their writings related to the major results on that subject, explaining their motivation and their goal.
 
 Ahlfors writes in his 1953 paper celebrating the 100 years anniversary of Riemann's  inaugural dissertation:
\begin{quote}\small
In the premature death of Teichm\"uller, geometric function theory, like other branches of mathematics, suffered a grievous loss. He spotted the image of Gr\"otzsch's technique, and made numerous applications of it, which it would take me too long to list. Even more important, he made systematic use of extremal quasiconformal mappings, a concept that Gr\"otzsch had introduced in a very simple special case. Quasiconformal mappings are not only a valuable tool in questions connected with the type problem,  but through the fundamental although difficult work of Teichm\"uller it has become clear that they are instrumental in the study of modules of closed surfaces. Incidentally, this study led to a better understanding of the role of quadratic differentials, which in somewhat mysterious fashion seem to enter in all extremal problems connected with conformal mapping.
\end{quote}
In the commentaries he made on his paper \emph{On quasiconformal mappings} \cite{Ahlfors1954} (1954) in his Collected Papers edition (\cite{Ahlfors-collected} Vol. 1, p. 1), Ahlfors writes:

\begin{quote}\small
More than a decade had passed since Teichm\"uller wrote his remarkable paper \cite{T20} on extremal quasiconformal mappings and quadratic differentials. It has become increasingly evident that Teichm\"uller's ideas would profoundly influence analysis and especially the theory of functions of one complex variable, although nobody at that time could foresee the extent to which this would be true. The foundations of the theory were not commensurate with the loftiness of Teichm\"uller's vision, and I thought it was time to reexamine the basic concepts. My paper has serious shortcomings, but it has nevertheless been very influential and has led to a resurgence of interest in quasiconformal mappings and Teichm\"uller theory. [...]
\\
Based in this definition the first four chapters are a careful and rather detailed discussion of the basic properties of quasiconformal mappings to the extent that they were known at that time. In particular a complete proof of the uniqueness part of Teichm\"uller's theorem was included. Like all other known proofs of the uniqueness it was modeled on Teichm\"uller's own proof, which used uniformization and the length-area method. Where Teichm\"uller was sketchy I tried to be more precise.
\\
In the original paper Teichm\"uller did not prove the existence part of his theorem, but in a following paper  \cite{T29} he gave a proof based on a continuity method. I found his proof rather hard to read and although I did not doubt its validity I thought that a direct variational proof would be preferable. My attempted proof on these lines had a flaw, and even my subsequent correction does not convince me today. In any case my attempt was too complicated and did not deserve to succeed. Later L. Bers \cite{B1960}  published a very clear version of Teichm\"uller's proof. The final credit belongs to R. Hamilton, who gave an amazingly short and direct proof of the existence theorem. The consensus today is that the existence part is easier to prove than the uniqueness.
\end{quote}

 Ahlfors adds in a note in his      
      \emph{Collected papers} edition (Vol. 2, p. 155):        
      \begin{quote}\small
      During the period of roughly 1958-1961, L. Bers and I were busy developing the fundamentals of the theory of Teichm\"uller space. We were in constant touch with each other, and it is sometimes hard to tell, and difficult to remember, who did what. If I were an impartial judge I would give the laurel to Bers for having introduced what has become the standard approach to Teichm\"uller space as an open subset of the complex linear space of quadratic differentials, and I would give myself honorable mention for having helped develop the analytic techniques.
      \end{quote}

The complex structure of Teichm\"uller space was worked out by Ahlfors and Bers, but, as Weil did it,  from the point of view of the Beltrami equation. In a paper published in 1960 on the complex structure \cite{Ahlfors1960}, Ahlfors says that he does not agree with Teichm\"uller's statement made in \cite{T32} that his work in \cite{T20} on the metric structure is of no use for the complex structure of that space (in fact, Ahlfors bases his construction of the complex structure on that work). He then adds the following: 
\begin{quote}\small
The problem [of moduli]  is not a clear cut one, and several formulations seem equally reasonable. In his paper of 1944, Teichm\"uller analyzes the situation and ends by setting his sights extremely high. He does not claim complete success, and due to the sketchy nature of the paper I have been unable to determine just how much he proved. [...] I have not succeeded in showing that the Teichm\"uller topology is the only topology which permits a complex analytic structure of the desired kind. Nevertheless, there are sufficient indications that the Teichm\"uller topology is the only one that is reasonable to consider.
\end{quote}
              
              Let us now talk about Bers.
              
 Bers started his research by studying partial differential equations, and he used them later in uniformization problems, in the tradition of Poincar\'e. In his paper  ICM 1958 paper, 
\emph{Spaces of Riemann surfaces}, \cite{Bers-Spaces1960}, he writes:
 \begin{quote}\small
This address is a report on recent work, partly not yet published, on the classical problem of moduli. Much of this work consists in clarifying and verifying assertions of Teichm\"uller whose bold ideas, though sometimes stated awkwardly and without complete proofs, influenced all recent inverstigations, as well as the work of Kodaira and Spencer on the higher dimensional case. Following Teichm\"uller, we consider not the space of closed Riemann surfaces of a given genus $g$ but rather an appropriate covering space and certain related spaces. [...]
Our main technical tools are uniformization theory and the theory of partial differential equations. The problem of moduli has also an algebro-geometric aspect, but the topological and analytical methods used here are, of course, restricted to the classical case. 
\end{quote}

In the same paper, after presenting Teichm\"uller's proof for the existence of extremal quasiconformal mappings, Bers writes:
 \begin{quote}\small
  The statement that $\mathcal{T}_g$ is a $(6g-6)$-cell is already contained in the work of Fricke. Fricke's proof is quite different and very difficult to follow.
  
 The existence of a ``natural" complex analytic structure in $\mathcal{T}_g$ has been asserted by Teichm\"uller \cite{T32}; the first proof was given by Ahlfors (The complex analytic structure of the space of complex closed Riemann surfaces, to appear) after Rauch (On the transcendental moduli, NAS, 42-49, 1955) showed how to introduce complex-analytic coordinates in the neighborhood of any point of $\mathcal{T}_g$ which is not a hyperelliptic surface. Other proofs are due to Kodaira (to appear) and to Weil (Bourbaki). 

 \end{quote}
 He states as an application: ``Every canonical dissection of a closed Riemann surface $S=U/G$ can be deformed into a dissection which maps into a convex non-euclidean polygon in $U$."
He adds: ``This was stated, with a different and complicated proof, by Fricke"  and he refers to Fricke and Klein \cite{fk}. 
 At the end of the paper, Bers asks: 

\begin{quote}
Is $\mathcal{T}_g$ a subset of $\mathbb{C}^{3g-3}$?
\end{quote}
 \begin{quote}
The space of (unmarked) Riemann surfaces is the factor-space $\mathcal{T}_g/\Gamma_g$, $\Gamma_g$ being the so-called mapping class group. Give a precise description of this space. 
\end{quote}

In his paper \cite{Bers1960a} (1960), he announces, with a sketch of a proof, the following theorem:
        \begin{theorem}\label{th:BBe}
        Teichm\"uller space $\mathcal{T}_{g,n}$ is (holomorphically equivalent to) a bounded domain in a complex vector space.
        \end{theorem}
        
          For the existence of the complex structure,
 Bers mentions in \cite{Bers1960a} five papers (by Ahlfors, by himself, by Kodaira-Spencer, by Rauch and by Weil).\footnote{Regarding Weil, the reference is to his Bourbaki seminar \cite{Weil1958c}.} The proof of Theorem \ref{th:BBe}
 is based on the Ahlfors-Bers \emph{Riemann's mapping theorem theorm for variable metrics} \cite{ABers} (1960). In the next year, in a note \cite{Bers1960ab}, Bers announces a mistake in the proof in \cite{Bers1960a} (which is not in the Ahlfors-Bers Riemann mapping theorem and which nonetheless does not affect the statement of  \ref{th:BBe}, and he gives an outline of a new proof. He announces that in the meanwhile, ``simultaneously and independently," Ahlfors found another proof.

  In 1964, Bers and Ehrenpreis showed that any finite-dimensional Teichm\"uller space can be embedded as  a domain of holomorphy in some $\mathbb{C}^N$  \cite{BE1964}. This is equivalent to the fact that Teichm\"uller space is a Stein manifold, a question which was asked by Grothendieck in \cite{G1}.
The Bers embedding is defined in terms of Kleinian groups, and in that theory, quasiconformal mappings and fine properties of their Schwarzian derivatives play a central role. The embedding is based on Bers' \emph{simultaneous uniformization theorem} \cite{Bers-simultaneous}. Bers also worked on the generalization of the deformation theory to surfaces of infinite type and to the universal Teichm\"uller space which later on appears to be of use in physics. Bers worked on deformations of Kleinian groups and he studied boundaries of Teichm\"uller spaces and of spaces of Kleinian groups. It is not possible to review this work here.  We refer to Bers' collected papers \cite{Bers-Selected}, and the survey by Kra and Maskit \cite{MK}.

 We would like to conclude this article by a quote from Adolphe Buhl, commenting the work of Poincar\'e \cite{Buhl}.  
 \begin{quote}\small
 A brilliant work is not some novel entanglement which is so complex that nobody, until then, had managed to build something so complicated.
 
 On the contrary, this is the clear perception of a very simple harmony which the unrefined eyes of the  contemporaries did not see. It is the rapid sketch of an artist, unexpected though luminous and striking, as soon as it exists. This is true in all domains of thought and even more particularly in the mathematical domain. 
 \end{quote}
 
    In writing this paper, we were interested in the origin of the ideas that are at the basis of the work of Bers. We hope to have communicated this interest to the reader.  We also tried to convey the conviction that mathematics, despite its division into several subfields, is a coherent living organism.  Mathematical work requires personal investment but it is above all a community work. Finally, the work we describe shows that mathematicians are persistent and patient people, and the important constructions are only achieved after several decades, sometimes centuries, of collective work, exactly like Middle Ages cathedrals.

\end{document}